\documentclass[preprint]{elsarticle}

\biboptions{sort}

\usepackage{amsmath}            
\usepackage{amsfonts}
\usepackage{amssymb}		
\usepackage{amsthm}		
\usepackage{eucal}		
\usepackage{mathtools}          
\usepackage{lscape}             

\usepackage{tabularx}           

\usepackage{arydshln}           
\renewcommand\ge\geqslant
\renewcommand\le\leqslant

\newtheorem{theorem}{Theorem}[section]
\newtheorem{corollary}[theorem]{Corollary}
\newtheorem{proposition}[theorem]{Proposition}

\newtheorem{conjecture}[theorem]{Conjecture}

\theoremstyle{definition}
\newtheorem{definition}[theorem]{Definition}
\newtheorem{example}[theorem]{Example}

\theoremstyle{remark}
\newtheorem{remark}[theorem]{Remark}

\numberwithin{equation}{section}

\numberwithin{equation}{section}

\newcommand{\defeq}{\mathrel{:=}}
\newcommand{\eqdef}{\mathrel{=:}}

\newcommand\myatop[2]{\genfrac{}{}{0pt}{}{#1}{#2}}

\newcommand\genvert[2]{\genfrac{|}{|}{0pt}{}{#1}{#2}}
\newcommand\stircyc[2]{\genfrac{[}{]}{0pt}{}{#1}{#2}}
\newcommand\stirsub[2]{\genfrac{\{}{\}}{0pt}{}{#1}{#2}}
\newcommand\eulerian[2]{\genfrac\langle\rangle{0pt}{}{#1}{#2}}

\newcommand\tgenvert[2]{\genfrac{|}{|}{0pt}{1}{#1}{#2}}
\newcommand\tstircyc[2]{\genfrac{[}{]}{0pt}{1}{#1}{#2}}
\newcommand\tstirsub[2]{\genfrac{\{}{\}}{0pt}{1}{#1}{#2}}

\newenvironment{smallarray}[1]
 {\null\,\vcenter\bgroup\scriptsize
  \arraycolsep=.13885em
  \hbox\bgroup$\array{@{}#1@{}}}
 {\endarray$\egroup\egroup\,\null}

\makeatletter
\newenvironment{sizeequation}[1]{%
  \skip@=\baselineskip
  #1%
  \baselineskip=\skip@
  \equation
}{\endequation \ignorespacesafterend}

\newenvironment{sizealign}[1]{%
  \skip@=\baselineskip
  #1%
  \baselineskip=\skip@
  \align
}{\endalign \ignorespacesafterend}

\makeatother

\begin{document}

\begin{frontmatter}
\title{Triangular recurrences, generalized Eulerian numbers, and related number triangles}
\author{Robert S. Maier}
\ead{rsm@math.arizona.edu}
\address{Depts.\ of Mathematics and Physics, University of Arizona, Tucson,
AZ 85721, USA}
\begin{abstract}
  Many combinatorial and other number triangles are solutions of
  recurrences of the Graham--Knuth--Patashnik (GKP) type.  Such
  triangles and their defining recurrences are investigated
  analytically.  They are acted upon by a transformation group generated
  by two involutions: a~left--right reflection and an upper binomial
  transformation, acting row-wise.  The group also acts on the
  bivariate exponential generating function (EGF) of the triangle.  By
  the method of characteristics, the EGF of any GKP triangle has an
  implicit representation in~terms of the Gauss hypergeometric
  function.  There are several parametric cases when this EGF can be
  obtained in closed form.  One is when the triangle elements are the
  generalized Stirling numbers of Hsu and Shiue.  Another is when they
  are generalized Eulerian numbers of a newly defined kind.  These
  numbers are related to the Hsu--Shiue ones by an upper binomial
  transformation, and can be viewed as coefficients of connection
  between polynomial bases, in a manner that generalizes the classical
  Worpitzky identity.  Many identities involving these generalized
  Eulerian numbers and related generalized Narayana numbers are
  derived, including closed-form evaluations in combinatorially
  significant cases.
\end{abstract}
\end{frontmatter}

\section{Introduction}
\label{sec:intro}

\subsection{Notation}
\label{subsec:notation}

Recurrences of the form
\begin{equation}
  \genvert{n+1}{k+1} = 
         [\alpha n + \beta (k+1) + \gamma]\, \genvert{n}{k+1}
         + [\alpha' n + \beta'k + \gamma']\, \genvert{n}{k},
\label{eq:gkp}
\end{equation}
satisfied by an infinite triangle of numbers $\genvert{n}k$, $0\le
k\le n<\infty$, with $\genvert{n}k$ equal to zero by convention if
$k<0$ or~$k>n$, and normalized so that the apex element $\genvert00$
equals unity, occur in pure and applied combinatorics, the analysis of
discrete algorithms, and elsewhere in mathematics.  Graham, Knuth, and
Patashnik \cite[Chapter~6,\ Problem~94]{Graham94} have indicated the
need for a general theory of such triangular recurrences, which are
now said to be of P94 or GKP type~\cite{Salas2021}.

This would include the construction of explicit
solutions~$\genvert{n}k$ of minimal rank, for the broadest choices of
the parameter vectors $\alpha,\beta;\gamma$ and
$\alpha',\beta';\gamma'$, and the identification of especially simple
`fundamental' solutions in~terms of which other solutions can be
expressed.  Here, `rank' refers to the depth to which summations are
nested in any explicit formula, the summand(s) being products and
quotients of factorials and powers~\cite[\S5.7]{Comtet74}.  The
strongest results to~date in these directions include those of
Spivey~\cite{Spivey2011} and Barbero~G. et~al.~\cite{Barbero2014}, who
employed both series manipulations and generating functions.  The
present work builds on theirs.

Many familiar numerical or combinatorial triangles are solutions of
GKP recurrences.  The new, explicitly parametric notations
\begin{subequations}
\label{eq:12}
\begin{align}
  \genvert{n}{k}
  &\eqdef
  \left[
    \begin{array}{ll|l}
      \alpha, & \beta & \gamma \\
      \alpha', & \beta' & \gamma' \\
    \end{array}
    \right]_{n,k},
\\
\label{eq:rowpolys}
\sum_{k=0}^n
  \genvert{n}{k}
  \,t^k 
  &\eqdef
  \left[
    \begin{array}{ll|l}
      \alpha, & \beta & \gamma \\
      \alpha', & \beta' & \gamma' \\
    \end{array}
    \right]_{n}(t),
\\
\label{eq:egf}
\sum_{n=0}^\infty
\sum_{k=0}^n
  \genvert{n}{k}
  \,t^k \frac{z^n}{n!}
&\eqdef
  \left[
    \begin{array}{ll|l}
      \alpha, & \beta & \gamma \\
      \alpha', & \beta' & \gamma' \\
    \end{array}
    \right](t,z)
\end{align}
\end{subequations}
will be employed here.  The first symbolizes the infinite triangle
derived from the recurrence~(\ref{eq:gkp}), and the second is its
$n$'th row polynomial, the univariate ordinary generating function of
its $n$'th row.  The third is the bivariate exponential generating
function (EGF) of the triangle as a whole.  The six GKP parameters
will be allowed to be complex, like the generating function
arguments~$t,z$.  Transformed or modified versions of a GKP triangle
$\genvert{n}k$ or a GKP parameter array
$\left[\begin{smallarray}{cc|c}\alpha&\beta&\gamma\\\alpha'&\beta'&\gamma'\end{smallarray}\right]$
will be indicated by an asterisk, as ${\genvert{n}k}^*$ or
$\left[\begin{smallarray}{cc|c}\alpha&\beta&\gamma\\\alpha'&\beta'&\gamma'\end{smallarray}\right]^*$.
An alternative six-parameter notation will be introduced in
Section~\ref{subsec:newparam}.

Partial lists of GKP recurrences and solutions that have appeared in
the literature can be found in~\cite{Barbero2014,Theoret94,Theoret95}.
In particular,
\begin{equation}
\begin{alignedat}{2}
  \stirsub{n}{k} 
  &
  = \left[
    \begin{array}{cc|c}
      0, & 1 & 0 \\
      0, & 0 & 1
    \end{array}
    \right]_{n,k},
  \qquad\qquad
  &
  \stircyc{n}{k}
  &
  = \left[
    \begin{array}{cc|c}
      1, & 0 & 0 \\
      0, & 0 & 1
    \end{array}
    \right]_{n,k},
  \\
  \eulerian{n}{k} 
  &
  = \left[
    \begin{array}{cc|c}
      0, & 1 & 1 \\
      1, & -1 & 0
    \end{array}
    \right]_{n,k},
  \qquad\qquad
  &
  \binom{n}{k}
  &
  = \left[
    \begin{array}{cc|c}
      0, & 0 & 1 \\
      0, & 0 & 1
    \end{array}
    \right]_{n,k},
\end{alignedat}
\end{equation}
where $\stirsub{n}k$ are the Stirling subset numbers (also denoted by
$S(n,k)$ and called the Stirling numbers of the second kind),
$\stircyc{n}k$ are the Stirling cycle numbers (also denoted by
$(-1)^{n-k}s(n,k)$ and called the unsigned Stirling numbers of the
first kind), and $\eulerian{n}k$ are the Eulerian numbers (also
denoted in the traditional indexing by~$A_{n,k+1}$).  Each of these
triangles has a combinatorial interpretation.  For any
$n$\nobreakdash-set, $\stirsub{n}k$~counts the number of its
partitions into $k$~blocks, and $\stircyc{n}k$~counts the number of
its permutations that have $k$~cycles.  If the $n$\nobreakdash-set is
totally ordered, $\eulerian{n}k$~counts the number of its permuations
that have $k$~descents.  The Stirling and Eulerian numbers are
introduced in~\cite{Comtet74,Graham94,Riordan58}.  The Eulerian
numbers are reviewed in~\cite{Foata2010} (see also~\cite{Takacs79})
and are treated more abstractly in~\cite{Foata70,Petersen2015}.

\subsection{Context and overview}
\label{subsec:overview}

This paper introduces a new, GKP-type generalization of the Eulerian
numbers~$\eulerian{n}k$, in addition to developing further the
analytic theory of GKP recurrences and their solutions.  Two other
parametric GKP triangles are also studied.  The following remarks
place the new Eulerian numbers in context.

Hsu and Shiue~\cite{Hsu98} introduced a parametric, GKP-type
generalization of the Stirling numbers $\stirsub{n}k$
and~$\stircyc{n}k$, which subsumes various previously treated ones.
In a slight modification of their notation, it is
\begin{equation}
\label{eq:hsushiue}
S_{n,k}(a,b;\,r) \defeq
\left[
  \begin{array}{cc|c}
    -a, & b & r \\
    0, & 0 & 1
  \end{array}
  \right]_{n,k}.
\end{equation}
(Examples are listed in~\cite{Benyi2002,He2013,Hsu98}.)  They originally defined
the $S_{n,k}(a,b;r)$ numbers as coefficients of connection between
graded polynomial bases of factorial type, i.e.,
\begin{equation}
\label{eq:hsushiue2}
  (x)^{\underline{n},a} = \sum_{k=0}^n S_{n,k}(a,b;\,r) (x-r)^{\underline k,b}.
\end{equation}
When $(a,b;r)=(0,1;0)$, this reduces to the original definition of~$\stirsub{n}k$,
\begin{equation}
  x^n = \sum_{k=0}^n \stirsub{n}{k} \,x^{\underline k},
\end{equation}
which is the Newton--Gregory expansion of~$x^n$, and similarly when
$(a,b;r)=(-1,0;0)$, it reduces to the original definition
of~$\stircyc{n}{k}$.  (Here $(x)^{\underline{n},a}$ and~$x^{\underline
  k}$ denote the falling factorials $(x)(x-\nobreak
a)\dotsm\allowbreak[x-\nobreak(n-1)a]$ and $x(x-\nobreak
1)\dotsm\allowbreak [x-\nobreak(k-1)]$; rising factorials will be
indicated by an overbar.)  The explicit general formula
\begin{equation}
\label{eq:chakcorcino}
  S_{n,k}(a,b;\,r) = 
\frac{1}{b^kk!}\sum_{j=0}^k (-1)^{k-j} \binom{k}{j}(bj+r)^{\underline{n},a}
\end{equation}
was pointed out by Corcino~\cite{Corcino2001}.
Equation~(\ref{eq:chakcorcino}) is a rank\nobreakdash-1 formula that
subsumes the well-known formula for the Stirling subset
numbers~\cite{Comtet74,Graham94}, which is
\begin{equation}
\label{eq:stirling5}
  \stirsub{n}k = \frac{1}{k!}  \sum_{j=0}^k (-1)^{k-j}
  \binom{k}j j^n.
\end{equation}
But (\ref{eq:chakcorcino}) clearly does not apply when~$b=0$, as is
the case for the Stirling cycle numbers~$\stircyc{n}k$, for which a
rank\nobreakdash-$2$ formula must be used~\cite{Comtet74}.  For a
discussion of the related definitions (\ref{eq:hsushiue}),
(\ref{eq:hsushiue2}), and~(\ref{eq:chakcorcino}), see
\cite{Corcino2001} and~\cite[\S4.2]{Mansour2016}.\footnote{It should be
  noted that numbers equivalent to the $S_{n,k}(a,b;r)$ had been
  introduced previously by Singh Chandel~\cite{SinghChandel77} and
  Charalambides and Koutras~\cite{Charalambides83}.  Their starting
  point was not (\ref{eq:hsushiue}) or~(\ref{eq:hsushiue2}),
  but~(\ref{eq:chakcorcino}).}  When restricted to integer parameter
values, $S_{n,k}(a,b;r)$ has been interpreted
combinatorially~\cite{Benyi2002,Charalambides83,Corcino2001a,Maltenfort2020,Mansour2013}.

The GKP recurrences solved in the present paper are largely ones with
$\beta\beta'\neq\nobreak0$ (they are of `type~I' in the classification
of~\cite{Barbero2014}).  They would seem unrelated to the generalized
Stirling numbers of~(\ref{eq:hsushiue}).  But in fact, if
$\beta'\neq0$ then
\begin{equation}
\label{eq:infact}
  \left[
  \begin{array}{cc|c}
    \alpha, & \beta & \gamma \\
    0, & \beta' & \gamma'
  \end{array}
\right]_{n,k}
=
\left(\frac{\gamma'}{\beta'}\right)^{\overline k}(\beta')^k
  \left[
  \begin{array}{cc|c}
    \alpha, & \beta & \gamma \\
    0, & 0 & 1
  \end{array}
\right]_{n,k}
= (\gamma')^{\overline k, \beta'} S_{n,k}(-\alpha,\beta;\,\gamma).
\end{equation}
This is because if the lower parameter vector
$(\alpha',\beta';\gamma')$ of a GKP recurrence is equal to
$(0,0;c_1)$, replacing it by $(0,c_1;sc_1)$ will multiply the solution
$\genvert{n}k$ by~$s^{\overline k}$~\cite{Spivey2011}.  Also,
multiplying $(\alpha',\beta';\gamma')$ by any common constant
factor~$A$ clearly multiplies $\genvert{n}k$ by~$A^{k}$.  (If~the
upper vector $(\alpha,\beta;\gamma)$ were multiplied by~$A$, the
solution would be multiplied instead by~$A^{n-k}$.)  Thus solutions of
the Hsu--Shiue type~(\ref{eq:hsushiue}) are fundamental ones, in~terms
of which the solution of any GKP recurrence with $\alpha'=0$ and
$\beta'\neq0$ can be expressed.

Formulas for the GKP solution $\genvert{n}k$ when $\beta\beta'\neq0$
have been systematically derived in three cases: (A\,I)~$\alpha'=0$,
(A\,II)~$\alpha+\beta=0$, and (A\,III) $\frac{\alpha}\beta =
1+\frac{\alpha'}{\beta'}$.  In each case, $\genvert{n}k$~can be
expressed as a double sum involving subset and cycle Stirling numbers,
binomial coefficients, and generalized factorials.
(See~\cite{Spivey2011}, and \cite[Proposition~2.5]{Encinas2019} for a
compact restatement.)  The derived expressions for~$\genvert{n}k$ are
therefore of high rank: up~to~$5$, though the nested summations can
sometimes be simplified.  In case (A\,I), with $\alpha'=0$, it is
better to compute the triangle $\genvert{n}k$ by
applying~(\ref{eq:infact}) to reduce it to the generalized Stirling
triangle $S_{n,k}(-\alpha,\beta;\gamma)$, which can be computed from
formula~(\ref{eq:chakcorcino}).  The resulting formula
for~$\genvert{n}k$ is of rank~$1$.  The handling of case (A\,II), when
$\alpha+\beta=0$, can be similarly improved, because one can show that
there is an involution
\begin{equation}
  \left[
  \begin{array}{cc|c}
    \alpha, & \beta & \gamma \\
    0, & \beta' & \gamma'
  \end{array}
\right]_n(t)
  = t^n
  \left[
  \begin{array}{cc|c}
    \beta', & -\beta' & \gamma' \\
    \alpha+\beta, & -\beta & \gamma
  \end{array}
\right]_n\left(\frac1t\right).
\end{equation}
That is, row polynomials in cases (A\,I), (A\,II) are reversed or
`reflected' versions of each other.

A major theme of the present paper is that the analysis of case
(A\,III) leads naturally to a new generalization of the Eulerian
numbers~$\eulerian{n}k$, combinatorial interpretations of which remain
to be explored.  When $\beta\beta'\neq0$, one can assume without loss
of generality that $\beta'=-\beta$, which is an innocuous
normalization.  Case (A\,III) is then the case when
$\alpha+\alpha'=\beta$, which suggests defining the parametric
generalized Eulerian triangle
\begin{equation}
\label{eq:eulerian1}
  E_{n,k}(a,b;\,c_0,c_\infty) \defeq
  \left[
    \begin{array}{cc|c}
      -a, & b & c_0\\
      a+b,& -b & c_\infty
    \end{array}
    \right]_{n,k},
\end{equation}
which reduces to~$\eulerian{n}k$ when $(a,b;c_0,c_\infty)$ is
$(0,1;1,0)$ and to the traditionally indexed numbers $A_{n,k}$ when it
is $(0,1;0,1)$.  

Equation~(\ref{eq:eulerian1}) will be shown to imply that for
all~$n\ge0$,
\begin{equation}
\label{eq:eulerian2}
  (c_0+c_\infty)^{\overline{n},b} (x)^{\underline{n},a}
  =
  \sum_{k=0}^n E_{n,k}(a,b;\,c_0,c_\infty)\, (x-c_0)^{\underline{k},b} (x+c_\infty)^{\overline{n-k},b},
\end{equation}
which if $(c_0+c_\infty)^{\overline{n},b}\neq0$, defines the
generalized numbers $E_{n,k}(a,b;\,c_0,c_\infty)$, $0\le k\le n$, as
expansion coefficients.  (For a proof that $(x-c_0)^{\underline{k},b}
(x+c_\infty)^{\overline{n-k},b}$, $0\le k\le n$, are a basis for the
space of polynomials of degree~${\le n}$ in the indeterminate~$x$,
see~\cite{Carlitz78}.)  When $(a,b;c_0,c_\infty) = (0,1;1,0)$,
eq.~(\ref{eq:eulerian2}) reduces to the celebrated identity of
Worpitzky,
\begin{equation}
\label{eq:eulerian3}
  n!\,x^n = \sum_{k=0}^n  \eulerian{n}k\, (x-k)^{\overline{n}}.
\end{equation}
Thus (\ref{eq:eulerian2}) is a generalized Worpitzky identity.  Just
as it is clear from~(\ref{eq:hsushiue2}) that the generalized Stirling
numbers $S_{n,k}(a,b;r)$ are coefficients that connect a factorial
basis of the ring of polynomials (depending on~$a$) to another one
(depending on~$b$), so it is clear from (\ref{eq:eulerian2}) that for
all $n\ge0$, the generalized Eulerian numbers
$E_{n,k}(a,b;\,c_0,c_\infty)$, $0\le k\le n$, relate the factorial
element $(x)^{{\underline n},a}$ of the
$(n+\nobreak1)$\nobreakdash-dimensional space of polynomials of
degree~${\le n}$ to a $b$\nobreakdash-dependent `bifactorial' basis of
this space.

A rank\nobreakdash-$1$ formula for these numbers will also be derived,
applying when $b\neq0$, namely
\begin{equation}
\label{eq:eulerian4}
\begin{split}
& E_{n,k}(a,b;\,c_0,c_\infty) = 
\\
& \quad\frac{1}{b^k k!}
\sum_{j=0}^k (-1)^{k-j} \binom{k}j
(bn+c_0+c_\infty)^{\underline{k-j},b} (c_0 + c_\infty)^{\overline{j}, b} (bj+c_0)^{\underline{n},a}.
\end{split}
\end{equation}
When $(a,b;c_0,c_\infty) = (0,1;1,0)$ this reduces to the classical
formula~\cite{Comtet74,Graham94}
\begin{equation}
\label{eq:eulerian5}
  \eulerian{n}k = \sum_{j=0}^k  (-1)^{k-j} \binom{n+1}{k-j} (j+1)^n.
\end{equation}
The reader should notice that the five equations
(\ref{eq:eulerian1})--(\ref{eq:eulerian5}), dealing with the new
generalized Eulerian numbers, are bijective (by~design) with
(\ref{eq:hsushiue})--(\ref{eq:stirling5}), which dealt with the
Stirling numbers of Hsu and Shiue.

There are two subcases of the parametric $E_{n,k}(a,b;c_0,c_\infty)$,
extending the standard numbers $E_{n,k}(0,1;1,0)=\eulerian{n}k$, which
have been treated in the literature.  The first is when $c_0+\nobreak
c_\infty=b$.  It could be called the single-progression subcase,
because the expansion functions in~(\ref{eq:eulerian2}) then simplify:
\begin{equation}
  (x-c_0)^{\underline{k},b} (x+c_\infty)^{\overline{n-k},b}
  =
  [x-c_0-(k-1)b]^{\overline{n},b},
\end{equation}
with the two arithmetic progressions combining into one, as seen
in~(\ref{eq:eulerian3}).  This subcase includes the `degenerate'
Eulerian numbers of Carlitz~\cite[\S8]{Carlitz79}, which are of the
form $E_{n,k}(\lambda,1;c_0,1-c_0)$ and have recently been
combinatorially interpreted (when $c_0=1$) by
Herscovici~\cite{Herscovici2020}.  It has been extensively treated
in~\cite{Charalambides82} (see also~\cite{Hsu99}).

The second subcase is when $a=0$, causing the
function~$(x)^{\underline n,a}$ expanded in~(\ref{eq:eulerian2}) to
reduce to the monomial~$x^n$, as seen too in~(\ref{eq:eulerian3}).
This could be called the Carlitz--Scoville subcase.  It can be traced
to~\cite{Carlitz74}, in an equivalent symmetric formulation, and also
to the solution of a stochastic model of habitat selection by certain
insect larvae~\cite{Charalambides91,Janardan88}.

Generalized Eulerian numbers of the $a=0$ type have appeared in many
applications.  The numbers $E_{n,k}(0,1;u,v)$ are the $(u,v)$-Eulerian
numbers, which when $u,v\in\mathbb{N}$ have a combinatorial
interpretation~\cite{Barbero2015}.  When $u=0$ they reduce to the
order\nobreakdash-$v$ Eulerian numbers $A_{n,k}^{(v)}$
of~\cite{Dillon68}.  When $v=0$ they are related by scaling to the
$1/K$\nobreakdash-Eulerian numbers~\cite{Ma2015}, which are of the
form $E_{n,k}(0,K;1,0)$.

The numbers $E_{n,k}(0,1;r,1-r)$ count the permutations of an ordered
$n$\nobreakdash-set that have
$k$~$r$\nobreakdash-descents~\cite{Foata70}, and are identical to the
numbers $E_{n,k}(0,1;\delta,1-\delta)$ studied in~\cite{Harris94}.
They satisfy both of the preceding conditions: $c_0+\nobreak
c_\infty=b$ as well as $a=0$.  The numbers $E_{n,k}(0,2;1,1)\eqdef
{\eulerian{n}{k}}_B$ count the number of signed permutations of an
ordered $n$\nobreakdash-set which have $k$~`signed descents.'
(See~\cite{Bagno2019,Brenti94} and \cite[$\texttt{A060187}$]{OEIS2022}.)  They
are called the type\nobreakdash-$B$ Eulerian numbers or the MacMahon
numbers, and when $(a,b;c_0,c_\infty)=\allowbreak(0,2;1,1)$,
eq.~(\ref{eq:eulerian2}) accordingly reduces to the Worpitzky identity of
type~$B$~\cite{Bagno2022,Brenti94}.  These numbers also satisfy both
conditions, as do the single-progression numbers of~\cite{Xiong2013}.

Restricted to the subcases $c_0+c_\infty=b$ and/or $a=0$, the
identities (\ref{eq:eulerian2}) and~(\ref{eq:eulerian4}) are known.
(For instance, see \cite[eq.~(31)]{Hsu99},
\cite[eq.~(3.1)]{Janardan88}, and~\cite[eq.~18]{Barbero2015}.)  But in
their full generality, (\ref{eq:eulerian2}) and~(\ref{eq:eulerian4})
appear to be new.

Many additional identities involving the new numbers
$E_{n,k}(a,b;c_0,c_\infty)$ and the Hsu--Shiue numbers
$S_{n,k}(a,b;r)$ are derived below.  They include explicit formulas
holding for certain choices of parameter, including ones of
combinatorial significance.  A~fundamental tool is the method of
characteristics, applied as a solution technique to the partial
differential equation (PDE) satisfied by the bivariate generating
function of any GKP triangle~$\genvert{n}k$.

This parametric PDE is acted upon by an order\nobreakdash-$6$
transformation group isomorphic to~$S_3$, the group of permutations of
$3$~letters.  On the triangle level, this group is generated by two
involutive transformations $\genvert{n}k\mapsto{\genvert{n}k}^*$ that
act row-wise: the reflection transformation~(RT) $k\leftarrow n-k$,
and a so-called upper binomial transformation~(UBT\null).  They map
between (respectively) case\nobreakdash-(A\,I) and
case\nobreakdash-(A\,II) triangles, and case\nobreakdash-(A\,II) and
case\nobreakdash-(A\,III) ones.  The two transformations are
illustrated by (see Theorem~\ref{thm:ubtpair})
\begin{equation}
\label{eq:seesin}
E_{n,k}(a,b;\,c_0,c_\infty) = \sum_{j=k}^n (-1)^{j-k} \binom{j}{k}
(c_0+c_\infty)^{\overline{n-j},b} S_{n,n-j}(-a,b;\,c_\infty).
\end{equation}
This relates the new numbers $E_{n,k}$ to those of Hsu and Shiue, and
is an alternative to~(\ref{eq:eulerian4}).  One sees
in~(\ref{eq:seesin}) an initial reflection $j\leftarrow\nobreak
n-\nobreak j$, performed on a case\nobreakdash-(A\,I) triangle
$(c_0+\nobreak c_\infty)^{\overline{j},b}\allowbreak
S_{n,j}(-a,b;\nobreak c_\infty)$ and yielding a
case\nobreakdash-(A\,II) one; and a subsequent binomial
transformation, yielding the triangle $E_{n,k}(a,b;\allowbreak
c_0,c_\infty)$, which belongs to case~(A\,III\null).

By applying the method of characteristics and well-known facts about
the Gauss hypergeometric function~${}_2F_1$, it is possible to derive
an explicit formula for the bivariate exponential generating function
$G(t,z)$ of a GKP triangle not merely in the generalized
Stirling--Eulerian case~(A), which subsumes cases
(A\,I),\allowbreak(A\,II),\allowbreak(A\,III), but also in two others:
(B),~called here the generalized Narayana case, and (C),~called the
generalized secant--tangent case.  If in a GKP recurrence
$\beta\beta'\neq0$ and the pair $(\beta,\beta')$ is adjusted to
equal~$(2,-2)$, which can be done without loss of generality, then the
corresponding cases (B\,I),\allowbreak(B\,II),\allowbreak(B\,III)
become $(\alpha,\alpha')=(1,-2),\allowbreak(-4,3),\allowbreak(1,3)$,
and cases (C\,I),\allowbreak(C\,II),\allowbreak(C\,III) become
$(\alpha,\alpha')=(-1,2),\allowbreak(0,1),\allowbreak(-1,1)$.  Many
explicit formulas for the triangle elements~$\genvert{n}k$ in these
cases are derived; especially, in case~(B\null).  More than a dozen
GKP triangles of the generalized Narayana kind for which
$\genvert{n}k$ can be expressed as a hypergeometric term have been
identified in the OEIS~\cite{OEIS2022}, and are tabulated below.  (See
Tables \ref{tab:2} and~\ref{tab:3}.)

By following a context-free grammar approach to exponential
structures~\cite{Chen93,Dumont96}, some formulas for
case\nobreakdash-(C) GKP triangles are also derived, which justify the
`secant--tangent' description.  In~fact, the approach leads to
interesting formulas in all the just-mentioned cases.  It must be said
that although cases (A),\allowbreak(B),(C) are treated in isolation
here, they are related: by quadratic changes of variable, certain
case\nobreakdash-(A) generating functions can be reduced to
case\nobreakdash-(C) ones, and certain case\nobreakdash-(B) ones to
case\nobreakdash-(A) ones.  This is illustrated by the combinatorics
of polytopes~\cite{Fomin2007}.  The \emph{f}\nobreakdash-vectors of
$n$\nobreakdash-dimensional permutohedra (of~either type~$A_n$ or
type~$B_n$) are the rows of a certain case\nobreakdash-(A) GKP
triangle, and the corresponding $\gamma$\nobreakdash-vectors
(quadratically reduced) are the rows of a case\nobreakdash-(C) one.
For associahedra, there is a similar reduction from case~(B) to
case~(A\null).  But quadratic transformations of GKP triangles are
left to another paper.

The body of this paper is structured as follows.  GKP triangle EGF's,
and the $S_3$ transformation group acting on EGF's or row-wise on
triangles, are introduced in Section~\ref{sec:gfs}.  In
Section~\ref{sec:characteristics} the method of characteristics is
applied to the EGF PDE, and a new GKP parametrization adapted to the
$S_3$\nobreakdash-group and the construction of
${}_2F_1$\nobreakdash-based solutions is introduced.  The generalized
Stirling--Eulerian case~(A) is treated in the multi-part
Section~\ref{sec:stirlingeulerian}.  Many identities, including
contiguous function relations and explicit formulas, are derived.
Results on the generalized Narayana and secant--tangent triangles
(cases (B),(C)) are in Sections \ref{sec:narayana}
and~\ref{sec:sectang}, the latter including some grammar-based
identities.

\section{Generating functions}
\label{sec:gfs}

Suppose that a number triangle $\genvert{n}k$, $0\le k\le n<\infty$,
is a GKP triangle\footnote{As is clear from~(\ref{eq:gkp}), the
  convention of Spivey~\cite{Spivey2011} on the definition of the GKP
  parameters $\gamma,\gamma'$ is adhered~to here, rather than that of
  Barbero~G. et~al.~\cite{Barbero2014}.  In present notation, the
  parameters $\gamma,\gamma'$ used in~\cite{Barbero2014} would be
  written respectively as $\gamma-\alpha$ and
  $\gamma'-\alpha'-\beta'$.}: it satisfies a GKP-type
recurrence~(\ref{eq:gkp}) with parameters
$\alpha,\beta,\gamma;\alpha',\beta',\gamma'\in\mathbb{C}$, and the
initial condition $\genvert{0}0=1$.  (By convention, $\genvert{n}k=0$
if $k<0$ or~$k>n$.)

The EGF (exponential generating function)
$G(t,z)=\left[\begin{smallarray}{cc|c}\alpha&\beta&\gamma\\\alpha'&\beta'&\gamma'\end{smallarray}\right](t,z)$
defined in~(\ref{eq:egf}) equals unity at~$(0,0)$ and is defined and
analytic in a neighborhood of $(0,0)$
in~$\mathbb{C}\times\nobreak\mathbb{C}$; in~fact, in a neighborhood of
$\{t=0\} \cup \{z=0\}$.  It may be possible to compute the EGF in
closed form.  The computation may involve solving from scratch the PDE
(partial differential equation) satisfied by the EGF, or showing that
the EGF is related to that of another GKP triangle.  The latter
relationship may follow from the existence of a PDE-to-PDE
transformation.  `First degree' transformation among the PDE's
satisfied by the EGF's of GKP triangles form a group, as is summarized
in Theorem~\ref{thm:s3} below; quadratic transformations will be
explored elsewhere.

The following triangle-to-triangle transformations are elementary but
useful.  They `trim' a GKP triangle by removing its left or right
edge, provided that the edge consists only of zeroes (with the
exception of the apex element~$\genvert00=1$).  This phenomenon occurs
when $\gamma=0$, resp.\ $\gamma'=0$.
\begin{theorem}
\label{thm:trimming}
  Suppose that a number triangle\/ $\genvert{n}k$, $0\le k\le n<\infty$,
  satisfies a GKP recurrence with parameter array\/
  $\left[\begin{smallarray}{cc|c}\alpha,&\beta&\gamma\\\alpha',&\beta'&\gamma'\end{smallarray}\right]$
  and EGF~$G(t,z)$.  Then,
  \begin{enumerate}
  \item[{\rm(i)}] If\/ $\gamma=0$, $\gamma'\neq0$, the left-trimmed triangle\/
    ${\genvert{n}k}^* \defeq (\gamma')^{-1} \genvert{n+1}{k+1}$, ${0\le k\le n}$, is a GKP triangle with parameter array\/ $\left[\begin{smallarray}{cc|c}\alpha,&\beta&\alpha+\beta\\\alpha',&\beta'&\alpha'+\beta'+\gamma'\end{smallarray}\right]$.  Its EGF\/ equals\/ $(\gamma')^{-1}(\partial/\partial z) G(t,z)/t$.
  \item[{\rm(ii)}] If\/ $\gamma'=0$, $\gamma\neq0$, the right-trimmed
    triangle\/ ${\genvert{n}k}^{*} \defeq (\gamma)^{-1}
    \genvert{n+1}{k}$, ${0\le k\le n}$, is a GKP triangle with
    parameter array\/ $\left[\begin{smallarray}{cc|c}\alpha,&\beta&\alpha+\gamma\\\alpha',&\beta'&\alpha'\end{smallarray}\right]$.  Its EGF\/ equals\/ $(\gamma)^{-1}(\partial/\partial z)G(t,z)$.
  \end{enumerate}
\end{theorem}
\begin{proof}
  Both statements follow by elementary series manipulations.  The
  factors $(\gamma')^{-1}$, $(\gamma)^{-1}$ are included to satisfy
  the normalization ${\genvert00}^*=1$.
\end{proof}

\begin{remark}
\label{rmk:mid}
  A third type of trimming will be encountered below, in the proof of
  Theorem~\ref{thm:2contigsE} and elsewhere.  Let the row polynomials
  $\left[\begin{smallarray}{cc|c}\alpha&\beta&\gamma\\\alpha'&\beta'&\gamma'\end{smallarray}\right]_n(t)$
  of a GKP triangle be denoted by $G_n(t)$, $n\ge0$.  The first two
  row polynomials are $G_0(t)=1$ and
  $G_1(t)=\allowbreak\gamma+\nobreak\gamma't$, and it can be shown by
  induction that if $\gamma+\nobreak\gamma't$ equals
  $A(\beta+\nobreak\beta't)$ for some constant~$A$, each $G_{n+1}(t)$,
  $n\ge0$, will be a multiple of $\beta+\nobreak\beta't$; and if
  $A\neq0$ and $\beta\beta'\neq0$, the quotients
  $G_{n+1}(t)/\left[A(\beta+\nobreak\beta't)\right]$, $n\ge0$, will be
  the row polynomials of a new, `mid-trimmed' GKP triangle, with
  parameter array
  $\left[\begin{smallarray}{cc|c}\alpha,&\beta&\alpha+\gamma\\\alpha',&\beta'&\alpha'+\beta'+\gamma'\end{smallarray}\right]$.
\end{remark}

\begin{example}
\label{ex:22}
  Three distinct EGF's for the standard Eulerian numbers appear in the
  literature~\cite{Comtet74}, but Theorem~\ref{thm:trimming} relates
  them.\footnote{Traditionally the Eulerian numbers were denoted
    by~$A_{n,k}$, and were defined and nonzero for $1\le\nobreak
    k\le\nobreak n$, with $A_{n,k}=0$ if $k>n$ by convention.  To~fit
    them into a GKP framework, one must also set $A_{0,0}=1$ and
    $A_{n,0}=\nobreak0$, $n\ge1$.  In the modern indexing
    $\eulerian{n}k$~signifies $A_{n,k+1}$, except that
    $\eulerian00=1$; note also that when $0\le k\le n$,
    $\eulerian{n+1}k$ equals $A_{n+1,k+1}$ without exception.  The
    occasionally used notation $\eulerian{n}{k-1}$ should be
    understood as signifying $A_{n,k}$.}  They are
  \begin{subequations}
  \begin{alignat}{2}
    \sum_{n=0}^\infty \sum_{k=0}^n \eulerian{n}{k-1}\,t^k\frac{z^n}{n!} &= 
    \left[ 
      \begin{array}{cc|c}
        0,&1&0\\
        1,&-1&1
      \end{array}
      \right](t,z)
    &{}={}&
    \frac{1-t}{1-t{\rm e}^{(1-t)z}},
    \label{eq:threea}
    \\
    \sum_{n=0}^\infty \sum_{k=0}^n \eulerian{n}{k}\,t^k\frac{z^n}{n!} &= 
    \left[ 
      \begin{array}{cc|c}
        0,&1&1\\
        1,&-1&0
      \end{array}
      \right](t,z)
    &{}={}&
    \frac{1-t}{{\rm e}^{(t-1)z}-t}
    \label{eq:threeb}
    \\
    \sum_{n=0}^\infty \sum_{k=0}^n \eulerian{n+1}{k}\,t^k\frac{z^n}{n!} &= 
    \left[ 
      \begin{array}{cc|c}
        0,&1&1\\
        1,&-1&1
      \end{array}
      \right](t,z)
    &{}={}&
    \frac{(1-t)^2 {\rm e}^{(1+t)z}}{({\rm e}^{tz} - t{\rm e}^z)^2}.
    \label{eq:threec}
  \end{alignat}
  \end{subequations}
Left-trimming the triangle in~(\ref{eq:threea}) and right-trimming the
triangle in~(\ref{eq:threeb}) both yield the triangle
in~(\ref{eq:threec}).  From an analytic rather than a combinatorial
point of view, there is little to choose between the classical
definition $A_{n,k}=\eulerian{n}{k-1}$ and the modern shifted
definition~$\eulerian{n}{k}$ of these numbers.
\end{example}

The following theorem was mentioned in the introduction and is also
easily proved.
\begin{theorem}
\hangindent\leftmargini         
\label{thm:thm3}
      {\rm(i)} Let\/ $\genvert{n}k$ denote the GKP triangle with parameters\/
      $(\alpha,\beta,\gamma;\allowbreak\alpha',\beta',\gamma')$, and let\/
      $G(t,z)$ denote its~EGF\null.  Then, the GKP triangle with
      parameters\/ $(A\alpha,A\beta,A\gamma;\allowbreak
      B\alpha',B\beta',B\gamma')$ will be\/ $A^{n-k}B^k\genvert{n}k$, and if\/
      ${A\neq0}$ the EGF of this triangle will be\/ $G(Bt/A,Az)$.
      \begin{itemize}
      \item[{\rm(ii)}] Let\/ $\genvert{n}k$ denote the GKP triangle with
        parameters\/ $(\alpha,\beta,\gamma;\allowbreak0,0,\gamma')$.  Then,
        the GKP triangle with parameters\/ $(\alpha,\beta,\gamma;\allowbreak
        0,\gamma',s\gamma')$ will be\/ $s^{\overline k}\genvert{n}k$.
      \end{itemize}
\end{theorem}

\begin{example}
\label{ex:demorgan}
  An illustration of part~(ii) of the theorem, with $s=1$, is provided
  by the De~Morgan numbers $\textrm{Surj}(n,k)$, $0\le k\le n$, which
  count the number of maps from an $n$\nobreakdash-set onto a
  $k$\nobreakdash-set.  By examination, they satisfy a recurrence of
  GKP type and equal
  $\left[
    \begin{smallarray}{cc|c}0,&1\,&\,0\\
      0,&1\,&\,1
    \end{smallarray}\right]_{n,k}$.
  As
  $\stirsub{n}k=\left[
    \begin{smallarray}{cc|c}0,&1\,&\,0\\
      0,&0\,&\,1
    \end{smallarray}
    \right]_{n,k}$,
  one has $\textrm{Surj}(n,k) = (1)^{\overline k}\stirsub{n}k = k!\stirsub{n}k$.
\end{example}

Many sophisticated transformations of GKP triangles or their EGF's are
based on transformations of the parametric PDE satisfied by the
latter.

\begin{theorem}
\hangindent\leftmargini         
\label{thm:pde}
{\rm(i)} The EGF\/
  $G(t,z)=\left[\begin{smallarray}{cc|c}\alpha,&\beta&\gamma\\\alpha',&\beta'&\gamma'\end{smallarray}\right](t,z)$
  satisfies the first-order PDE
  \begin{equation}
    \label{eq:thepde}
    \left[\mathcal{A}(t)z-1\right]\frac{\partial G}{\partial z}
    +\mathcal{B}(t) \frac{\partial G}{\partial t} + \mathcal{C}(t)G=0,
  \end{equation}
  where\/ $\mathcal{A}(t)=\alpha+\alpha't$, $\mathcal{B}(t)=(\beta
  +\beta't)t$, and\/ $\mathcal{C}(t)=\gamma+\gamma't$, with the
  initial condition\/ $G(t,0)\equiv1$.
  \begin{itemize}
  \item[{\rm(ii)}] The row polynomials\/
    $G_n(t)=\left[\begin{smallarray}{cc|c}\alpha,&\beta&\gamma\\\alpha',&\beta'&\gamma'\end{smallarray}\right]_n\!(t)$ satisfy the differential recurrence
    \begin{equation}
      \label{eq:thediffrec}
      G_{n+1} = \left[\mathcal{A}(t)n + \mathcal{C}(t)\right] G_n + \mathcal{B}(t)G_n'
    \end{equation}
    and the initial condition\/ $G_0(t)\equiv1$.
  \item[{\rm(iii)}] If\/ $\beta\beta'\neq0$, the row polynomials\/
    $G_n(t)$ can be computed from
    \begin{equation}
      \label{eq:thediffrecsoln}
      \bigl[
        t^{1+\hat\alpha}(\beta+\beta't)^{1-\hat\alpha-\hat\alpha'} D_t
        \bigr]^n
      \frac{t^{\hat\gamma}}{(\beta+\beta't)^{\hat\gamma+\hat\gamma'}}
      =
      \frac{t^{\hat\alpha n + \hat\gamma} \:  G_n(t)}
      {(\beta+\beta't)^{(\hat\alpha+\hat\alpha')n + \hat\gamma + \hat\gamma'}},
    \end{equation}
    where\/ $\hat\alpha,\hat\alpha',\hat\gamma,\hat\gamma'$ signify\/
    $\alpha/\beta$, $-\alpha'/\beta'$, $\gamma/\beta$,
    $-\gamma'/\beta'$, and\/ $D_t = {\rm d}/{\rm d}t$.
\end{itemize}
\end{theorem}
\begin{proof}
  Substitute the definitions\/ (\ref{eq:egf}) and\/
  (\ref{eq:rowpolys}) of the EGF and row polynomials into\/
  (\ref{eq:thepde}) and\/~(\ref{eq:thediffrec}); and in both, use the
  triangular recurrence\/~(\ref{eq:gkp}).  By examination, the formula
  for $G_n(t)$ provided by~(\ref{eq:thediffrecsoln})
  satisfies~(\ref{eq:thediffrec}).
\end{proof}

Consider the transformation induced by a lifting map or change of
variables $(t^*,z^*)\mapsto(t,z)$, which is of the form $(t,z) =
(R(t^*),S(t^*) z^*)$ where $R,S$ are rational functions of their
argument.  ($R$~will be taken to be nonconstant.)  Substitution
into~(\ref{eq:thepde}) yields the following.
\begin{theorem}
  The transformed EGF
  \begin{displaymath}
    G^*(t^*,z^*) \defeq G(R(t^*),S(t^*)z^*) = G(t,z),
  \end{displaymath}
lifted from the EGF\/ $G(t,z)$ by the map\/ $(t^*,z^*)\mapsto(t,z)$
specified by\/ $(R,S)$, satisfies the first-order PDE
  \begin{equation}
    \label{eq:liftedpde}
    \left[\mathcal{A^*}(t^*)z^*-1\right]\frac{\partial G^*}{\partial z^*}
    +\mathcal{B^*}(t^*) \frac{\partial G^*}{\partial t^*} + \mathcal{C^*}(t^*)G^*=0,
  \end{equation}
  where
  \begin{subequations}
    \begin{align}
      \mathcal{A}^*(t^*) &= (\alpha+\alpha'R)S - (\beta+\beta'R)R\dot S/\dot R,
      \\
      \mathcal{B}^*(t^*) &= (\beta+\beta'R)RS/\dot R,
      \\
      \mathcal{C}^*(t^*) &= (\gamma+\gamma'R)S,
    \end{align}
  \end{subequations}
  or more compactly, 
  $\mathcal{A}^*=(\mathcal{A}\circ R)S-(\mathcal{B}\circ R)\dot S/\dot R$, $\mathcal{B}^* = (\mathcal{B}\circ R) S/\dot R$, 
  $\mathcal{C}^* = (\mathcal{C}\circ R)S$.  An overdot indicates
  differentiation with respect to\/~$t^*$.
\end{theorem}
For the lifted PDE~(\ref{eq:liftedpde}) to be of the GKP type, like
the original PDE (\ref{eq:thepde}), its coefficient functions
$\mathcal{A}^*(t^*)$, $\mathcal{B}^*(t^*)/t^*$, $\mathcal{C}^*(t^*)$,
must be degree\nobreakdash-$1$ polynomials in~$t^*$.  It should be
noted that any two liftings can be composed.  If the pairs $(R,S)$,
$(R^*,S^*)$ specify successive liftings, i.e.,
$(t,z)=(R(t^*),S(t^*)z^*)$ and
$(t^*,z^*)=(R^*(t^{**}),S^*(t^{**})z^{**})$, their composition
$(R,S)\circ(R^*,S^*)$ is the pair
\begin{equation}
\label{eq:comp}
  ({\bf R},{\bf S}) \defeq \left(R\circ R^*, (S\circ R^*)S^*\right),
\end{equation}
which specifies the composite map $(t^{**}\!,z^{**})\mapsto(t,z) =
({\bf R}(t^{**}),{\bf S}(t^{**})z^{**})$.  Also, liftings $(R,S)$ in
which $R$~has a compositional inverse~$\bar R$ (which will be the case
if $R(t^*)=\frac{\lambda + \mu t^*}{\rho + \sigma t^*}$ with
$\lambda\sigma\neq \mu\rho$, i.e., if $t^*\mapsto t$ is a
degree\nobreakdash-$1$ rational map) have compositional inverses of
the same form, i.e.,
\begin{equation}
\label{eq:inverse}
  (R,S)^{-1} = \left(\bar R, \frac1{S\circ \bar R}\right),
\end{equation}
which specifies the inverse map $(t,z)\mapsto (t^*,z^*)$.  Hence such
liftings form a group under composition.

Suppose that in a GKP recurrence, $\beta'=-\beta$, which if
$\beta\beta'\neq0$ is a mere matter of normalization.  Such a
restriction facilitates the study of the pairs $(R,S)$ that yield a
lifted PDE which is of the GKP type, like the original.  The case when
$R,S$ are rational of at~most degree~$1$ in~$t^*$ is especially easy
to treat.  By direct calculation, one finds that if $(R,S)$ equals
$(\frac1{t^*},t^*)$, resp.\ $(1-\nobreak t^*,-1)$, the lifted PDE will
indeed be of the GKP type, having a lifted or transformed parameter
array
$\left[\begin{smallarray}{cc|c}\alpha,&\beta&\gamma\\\alpha',&\beta'&\gamma'\end{smallarray}\right]^*$
equal to
$\left[\begin{smallarray}{cc|c}\alpha'-\beta,&\beta&\gamma'\\\alpha+\beta,&-\beta&\gamma\end{smallarray}\right]$,
resp.\ $\left[\begin{smallarray}{cc|c}-\alpha-\alpha',&\beta&-\gamma-\gamma'\\\alpha',&-\beta&\gamma'\end{smallarray}\right]$.
Both these liftings preserve the property $\beta'=-\beta$ and in~fact
leave $\beta,\beta'$ unchanged, though they transform in an
affine-linear way the vector comprising the other four parameters.
Both liftings are involutions, as follows from either (\ref{eq:comp})
or~(\ref{eq:inverse}).  By~(\ref{eq:comp}), the composition of either
$(\frac1{t^*},t^*)$ or $(1-\nobreak t^*,-1)$ with itself is the pair
$(t^*,1)$, which specifies the identity transformation.

When acting on a GKP-type EGF, row polynomial, and number triangle,
these two liftings yield the involutive transformation identities
\begin{subequations}
  \begin{align}
    \label{eq:firsta}
    \left[
    \begin{array}{cc|c}
      \alpha'-\beta, & \beta & \gamma' \\
      \alpha+\beta, & -\beta & \gamma
    \end{array}
    \right]
    \left(t^*,z^*\right)
    &=
    \left[
    \begin{array}{cc|c}
      \alpha, & \beta & \gamma \\
      \alpha', & -\beta & \gamma'
    \end{array}
    \right]
    \left(\frac{1}{t^*},t^*z^*\right)
    ,
    \\
    \label{eq:firstb}
    \left[
    \begin{array}{cc|c}
      \alpha'-\beta, & \beta & \gamma' \\
      \alpha+\beta, & -\beta & \gamma
    \end{array}
    \right]_n
    \left(t^*\right)
    &=
    (t^*)^n
    \left[
    \begin{array}{cc|c}
      \alpha, & \beta & \gamma \\
      \alpha', & -\beta & \gamma'
    \end{array}
    \right]_n
    \left(\frac{1}{t^*}\right)
    ,
    \\
    \label{eq:firstc}
    \left[
    \begin{array}{cc|c}
      \alpha'-\beta, & \beta & \gamma' \\
      \alpha+\beta, & -\beta & \gamma
    \end{array}
    \right]_{n,k}
    &=
    \left[
    \begin{array}{cc|c}
      \alpha, & \beta & \gamma \\
      \alpha', & -\beta & \gamma'
    \end{array}
    \right]_{n,n-k}
    ,
\end{align}
\end{subequations}
resp.
\begin{subequations}
  \begin{align}
    \label{eq:seconda}
    \left[
    \begin{array}{cc|c}
      -\alpha-\alpha', & \beta & -\gamma-\gamma' \\
      \alpha', & -\beta & \gamma'
    \end{array}
    \right]
    \left(t^*,z^*\right)
    &=
    \left[
    \begin{array}{cc|c}
      \alpha, & \beta & \gamma \\
      \alpha', & -\beta & \gamma'
    \end{array}
    \right]
    \left(1-t^*, -z^*\right)
    ,
    \\
    \label{eq:secondb}
    \left[
    \begin{array}{cc|c}
      -\alpha-\alpha', & \beta & -\gamma-\gamma' \\
      \alpha', & -\beta & \gamma'
    \end{array}
    \right]_n
    \left(t^*\right)
    &=
    (-1)^n
    \left[
    \begin{array}{cc|c}
      \alpha, & \beta & \gamma \\
      \alpha', & -\beta & \gamma'
    \end{array}
    \right]_n
    (1-t^*)
    ,
    \\
    \label{eq:secondc}
    \left[
    \begin{array}{cc|c}
      -\alpha-\alpha', & \beta & -\gamma-\gamma' \\
      \alpha', & -\beta & \gamma'
    \end{array}
    \right]_{n,k}
    &=
    (-1)^{n-k}
    \sum_{j=k}^n 
    \binom{j}{k}
    \left[
    \begin{array}{cc|c}
      \alpha, & \beta & \gamma \\
      \alpha', & -\beta & \gamma'
    \end{array}
    \right]_{n,j}
    .
\end{align}
\end{subequations}
Each of (\ref{eq:firsta}) and~(\ref{eq:seconda}) expresses a
transformed EGF $G^*(t^*,z^*)$, on the left, in~terms of an original
EGF $G(t,z)$, on the right; and is valid in a neighborhood of
$(t^*,z^*)=(0,0)$ in~$\mathbb{C}\times\nobreak\mathbb{C}$, the trivial
case $t^*=0$ not being covered by~(\ref{eq:firsta}). The
row-polynomial identities (\ref{eq:firstb}) and~(\ref{eq:secondb})
come by expanding (\ref{eq:firsta}) and~(\ref{eq:seconda}).  In
general, if $\genvert{n}{\cdot}(t)$ denotes the $n$'th row polynomial
and ${\genvert{n}{\cdot}}^*(t^*)$ its transform, one can write
\begin{equation}
  {\genvert{n}{\cdot}}^*\!(t^*)
  =
  S(t^*)\:\genvert{n}{\cdot}\left(R(t^*)\right)
  ,
\end{equation}
which (\ref{eq:firstb}) and~(\ref{eq:secondb}) exemplify.
Identity~(\ref{eq:secondc}) comes by expanding binomially the factor
$(1-t^*)^k$ in the summation that defines the row polynomial on the
right-hand side of~(\ref{eq:secondb}).

Equation (\ref{eq:firstc}) defines a \emph{reflection} transformation
(RT): it reverses each row of a GKP triangle, yielding a triangle with
altered GKP parameters.  (A~special case was mentioned in the
introduction.)  Equation~(\ref{eq:secondc}) defines a (signed,
involutive) \emph{upper binomial} transformation (UBT), which also
acts row-wise on any GKP triangle with $\beta'=-\beta$.  The
appearance of a UBT in this context was first pointed~out
in~\cite{Spivey2011}.  There is a literature on binomial transforms of
finite or infinite sequences, but most of it deals with lower rather
than upper transforms~\cite{Boyadzhiev2018}.  A~lower binomial
transform would be based not on the operator
$\sum_{j=k}^n\binom{j}k\times$ but on $\sum_{j=0}^k\binom{k}j\times$,
as in (\ref{eq:chakcorcino}) and~(\ref{eq:stirling5}).

The RT and UBT generate a group of transformations, all of which come
from liftings.  By examination, this group is of order~$6$ and is
isomorphic to~$S_3$, the symmetric group on three letters.  The action
of the six transformations, including the identity, RT, and~UBT, can
be summarized as follows.

{
\setlength{\tabcolsep}{4.4pt}
\begin{table}
  \begin{center}
  \begin{tabular}{@{}lllll}
    $g\in S_3$ & $\left(R(t^*),S(t^*)\right)$  & $\left[\begin{smallarray}{cc|c}\alpha,&\beta&\gamma\\ \alpha',&\beta'&\gamma'\end{smallarray}\right]^*$ & ${\genvert{n}k}^*$ \\
    \noalign{\smallskip}\hline
    $(0)(1)(\infty)$ & $(t^*,1)$ & $\left[\begin{smallarray}{cc|c}\alpha,&\beta&\gamma\\ \alpha',&-\beta&\gamma'\end{smallarray}\right]^{\vphantom{\overline n}}$ &  $\genvert{n}k$ \\[6pt]
    $(0\infty)(1)$ & $(\frac{1}{t^*},t^*)$ & $\left[\begin{smallarray}{cc|c}\alpha'-\beta,&\beta&\gamma'\\ \alpha+\beta,&-\beta&\gamma\end{smallarray}\right]$ &  $\genvert{n}{n-k}$\\[6pt]
    $(01)(\infty)$ & $(1-t^*,-1)$ & $\left[\begin{smallarray}{cc|c}-\alpha-\alpha',&\beta&-\gamma-\gamma'\\ \alpha',&-\beta&\gamma'\end{smallarray}\right]$ & $(-1)^{n-k} \sum_{j=k}^n \binom{j}{k} \genvert{n}j $ \\[6pt]
    $(1\infty)(0)$ & $\left(\frac{-t^*}{1-t^*},1-t^*\right)$ & $\left[\begin{smallarray}{cc|c}\alpha,&\beta&\gamma\\ \beta-\alpha-\alpha',&-\beta&-\gamma-\gamma'\end{smallarray}\right]$ & $(-1)^{k} \sum_{j=0}^k \binom{n-j}{n-k} \genvert{n}j $ \\[6pt]
    $(0\infty1)$ & $\left(\frac{-(1-t^*)}{t^*}, -t^*\right)$ & $\left[\begin{smallarray}{cc|c}\alpha'-\beta,&\beta&\gamma'\\ \beta-\alpha-\alpha',&-\beta&-\gamma-\gamma'\end{smallarray}\right]$ & $(-1)^{k} \sum_{j=n-k}^n \binom{j}{n-k} \genvert{n}j $ \\[6pt]
    $(01\infty)$ & $\left(\frac1{1-t^*}, -(1-t^*)\right)$ & $\left[\begin{smallarray}{cc|c}-\alpha-\alpha',&\beta&-\gamma-\gamma'\\ \alpha+\beta,&-\beta&\gamma\end{smallarray}\right]$ & $(-1)^{n-k} \sum_{j=0}^{n-k} \binom{n-j}{k} \genvert{n}j $
  \end{tabular}
  \end{center}
  \caption{Six degree-$1$ transformations of GKP triangles with
    $\beta'=-\beta$, which preserve $\beta$ and~$\beta'$.  The first
    four are the identity transformation and the involutions RT, UBT,
    and $\textrm{RT}\circ\textrm{UBT}\circ\textrm{RT}$.  The last two
    are the order\nobreakdash-$3$ transformations
    $\textrm{UBT}\circ\textrm{RT}$ and
    $\textrm{RT}\circ\textrm{UBT}$.}
\label{tab:1}
\end{table}
}

\begin{theorem}
  For each of the six rows in Table\/~{\rm\ref{tab:1}}, there is a
  transformation of the GKP triangle\/ $\genvert{n}k =
  \left[\begin{smallarray}{cc|c}\alpha,&\beta&\gamma\\\alpha',&\beta'&\gamma'\end{smallarray}\right]_{n,k}$
  with\/ $\beta'=-\beta$ to a new GKP triangle\/ ${\genvert{n}k}^* =
  \left[\begin{smallarray}{cc|c}\alpha,&\beta&\gamma\\\alpha',&\beta'&\gamma'\end{smallarray}\right]^*_{n,k}$
  with\/ 
  $\beta^*=\beta$ and $\beta'^*=\beta'$,
  performed thus:\/ ${\genvert{n}k}^*$,
  the parameter array of which is given in the third column, equals
  the expression given in the fourth.  For each\/ $n\ge0$, the new\/
  n'th row polynomial\/ ${\genvert{n}{\cdot}}^*(t^*)$ equals\/
  $S(t^*)\,\genvert{n}{\cdot}\left(R(t^*)\right)$, and the new EGF\/
  $G^*(t^*,z^*)$ comes from the old EGF\/ $G(t,z)$ as
  $G(R(t^*),S(t^*)z^*)$, in a neighborhood of\/ $(t^*,z^*)=(0,0)$ in\/
  $\mathbb{C}\times\nobreak\mathbb{C}$.
\label{thm:s3}
\end{theorem}

For each of these six maps $(t^*,z^*)\mapsto(t,z) =
(R(t^*),S(t^*)z^*)$, the map $t^*\mapsto R(t^*)$ is a degree-$1$
rational map that stabilizes the subset $\{0,1,\infty\}$ of the
projective $t$\nobreakdash-line, or equivalently permutes the points
$0,1,\infty$, which makes concrete the isomorphism to~$S_3$.  These
permutations (elements $g\in S_3$) are given in cycle notation in the
first column of the table.  The permutation $(0\infty)(1)$ specifies
the~RT, and $(01)(\infty)$ the UBT\null.  In all cases the function
$S(t^*)$ is equal to the denominator of~$R(t^*)$, up~to sign.  This
order\nobreakdash-$6$ group is quite different from the known group of
lower binomial transforms~\cite{Galuzzi98}.

The transformation $\genvert{n}k \mapsto {\genvert{n}k}^*$ specified
by the third involution $(1\infty)(0)$ is conjugated to the UBT by
the~RT: it is the composition $\textrm{RT} \circ \textrm{UBT} \circ
\textrm{RT}$.  It is a variant form of a sequence transformation of
Stanton and Sprott~\cite{Stanton62}.  The transformations specified by
the cyclic permutations $(0\infty1)$ and~$(01\infty)$, when acting on
(the rows~of) any GKP triangle with $\beta'=-\beta$, are not
involutive: they are sequence transformations of order~$3$, each being
both the inverse and the square of the other.  They are the
compositions $\textrm{UBT}\circ\textrm{RT}$ and
$\textrm{RT}\circ\textrm{UBT}$.

This $S_3$ transformation group can be extended in various ways.
(Compare Salas and Sokal~\cite{Salas2021}.)  One can append
$6$~additional elements, in each of which $S(t^*)$ is negated,
relative to what appears in the table.  This negation (in~effect, a
negation of~$z^*$) will multiply ${\genvert{n}k}^*$ by~$(-1)^n$, and
by Theorem~\ref{thm:thm3}(i), negate the entire array
$\left[\begin{smallarray}{cc|c}\alpha,&\beta&\gamma\\ \alpha',&\beta'&\gamma'\end{smallarray}\right]^*$
of transformed parameters.  So, although the additional $6$~elements
will preserve the property that $\beta'=-\beta$, they will negate both
$\beta$ and~$\beta'$.  The extended group $S_3\times\mathbb{Z}_2$ is
isomorphic to the dihedral group with $12$~elements.

One could also relax or alter the condition that the GKP triangle
being transformed satisfy $\beta'=-\beta$.  (The convention adopted
here that the condition $\beta'/\beta=-1$ is fundamental is largely
due to its holding for the Eulerian triangles; recall
Example~\ref{ex:22}.)  For any specified $\beta,\beta'$ with
$\beta\beta'\neq0$, there is a transformation group isomorphic
to~$S_3$ which depends only on the ratio~$\beta:\beta'$ and leaves
$\beta,\beta'$ invariant.  Each of its elements comes from a pair
$(R,S)$ in which $t^*\mapsto R(t^*)$ permutes the points
$0,-\beta'/\beta, \infty$ of the projective line.

An example of this is the original Stanton--Sprott
transformation~\cite[Theorem~3]{Stanton62}, which in present notation
is the involution
\begin{equation}
{\genvert{n}k}^* = 
\sum_{j=0}^k \binom{n-j}{n-k} (-1)^j \genvert{n}j.
\end{equation}
It is similar but not identical to the transformation specified
by~$(1\infty)(0)$.  By examination, it can be viewed as acting on (the
rows~of) any GKP triangle with $\beta'=\beta$, and preserves both
$\beta$ and~$\beta'$.  Its effect is summarized by
$\left[\begin{smallarray}{cc|c}\alpha,&\beta&\gamma\\ \alpha',&\beta'&\gamma'\end{smallarray}\right]^*
=
\left[\begin{smallarray}{cc|c}\alpha,&\beta&\gamma\\ -\beta+\alpha-\alpha',&\beta&\gamma-\gamma'\end{smallarray}\right]$,
and it comes from the pair $(R,S)=\left( \frac{-t^*}{1+t^*}, 1+t^*
\right)$.  The map $t^*\mapsto\frac{-t^*}{1+t^*}$ stabilizes not
$\{0,1,\infty\}$ but $\{0,-1,\infty\}$: in cycle notation, it is the
permutation~$(-1,\infty)(0)$.

\section{The method of characteristics}
\label{sec:characteristics}

The first-order PDE satisfied by the EGF $G(t,z)$ of any GKP triangle
having been derived (see Theorem~\ref{thm:pde}), it will be shown how
in several interesting cases, the PDE can be solved in closed form.
The method of characteristics, which has been applied previously to
the problem of GKP triangles~\cite{Barbero2014,Wilf2004}, will be
exploited to the full.  The key result is Theorem~\ref{thm:2f1formula}
below.

It was noted by Wilf~\cite{Wilf2004} that this method leads to special
functions, in particular the Gauss hypergeometric function~${}_2F_1$.
It will be seen that in three cases, the EGF is nonetheless an
elementary function of its arguments.  The three could be called
(A)~the generalized Stirling--Eulerian case, (B)~the generalized
Narayana case, and (C)~the generalized secant--tangent case.  Case~(A)
was introduced in Section~\ref{sec:intro} and is relatively familiar.
Like~(A), cases (B) and~(C) have three subcases: (I),~(II), and~(III),
which are related by transformations that belong to the
$S_3$\nobreakdash-group of the last section.  (Recall
Table~\ref{tab:1} and Theorem~\ref{thm:s3}; also see
Theorem~\ref{thm:unification} below.)  In the present section only the
EGF for subcase~(I) of each is computed, in Section~3.3.  Cases
(A),~(B), and~(C) are treated in greater generality in the respective
Sections 4,~5, and~6.

\subsection{A new GKP parametrization}
\label{subsec:newparam}

The six transformations of the $S_3$ transformation group, in particular
the GKP parameter maps
$\left[\begin{smallarray}{cc|c}\alpha,&\beta&\gamma\\ \alpha',&\beta'&\gamma'\end{smallarray}\right]
\mapsto
\left[\begin{smallarray}{cc|c}\alpha,&\beta&\gamma\\ \alpha',&\beta'&\gamma'\end{smallarray}\right]^*$,
can be written in a unified and manifestly symmetric form, given in
Theorem~\ref{thm:unification}.

For this a new notation is needed, which as a matter of convention
will be centered on the putatively fundamental case when
$(\beta,\beta')=(1,-1)$.  The parameter array
$\left[\begin{smallarray}{cc|c}\alpha,&1&\gamma\\ \alpha',&-1&\gamma'\end{smallarray}\right]$
can be written alternatively as the tableau
\begin{equation}
  \left[
    \begin{array}{ccc}
      0, & 1, & \infty \\
      \hline
      r_0, & r_1, & r_\infty \\
      g_0, & g_1, & g_\infty
    \end{array}
    \right],
\end{equation}
where
\begin{equation}
\label{eq:rewrite}
\begin{alignedat}{3}
  r_0 &= -\alpha,\qquad  & r_1 &= \alpha+\alpha',\qquad  &r_\infty &= 1-\alpha',\\
  g_0 &= \gamma,\qquad  & g_1 &= -\gamma-\gamma',\qquad  &g_\infty &= \gamma',
\end{alignedat}
\end{equation}
so that $r_0+r_1+r_\infty=1$ and $g_0+g_1+g_\infty=0$.  (The ordering
of the columns is arbitrary: if the parameter-pair $r_0,g_0$
lie in that order below~$0$, etc., an array of this kind has an
unambiguous meaning.)  Such new-style parameter arrays can be used as
specifications of GKP triangles, row polynomials, and EGF's, much as
parameter arrays
$\left[\begin{smallarray}{cc|c}\alpha,&1&\gamma\\ \alpha',&-1&\gamma'\end{smallarray}\right]$
or
$\left[\begin{smallarray}{cc|c}\alpha,&\beta&\gamma\\ \alpha',&\beta'&\gamma'\end{smallarray}\right]$
are used.  (Recall eqs.~(\ref{eq:12}).)

By Theorem~\ref{thm:thm3}(i), whenever $\beta\beta'\neq0$, one can write
\begin{equation}
\label{eq:recalled1}
\begin{aligned}
  \left[
    \begin{array}{cc|c}
      \alpha, & \beta &  \gamma\\
      \alpha', & \beta' &  \gamma'
    \end{array}
    \right]_{n,k}
  &= \beta^{n-k}(-\beta')^k 
  \left[
    \begin{array}{cc|c}
      \alpha/\beta, & 1 &  \gamma/\beta\\
      -\alpha'/\beta', & -1 &  -\gamma'/\beta'
    \end{array}
    \right]_{n,k}\\
  &= \beta^{n-k}(-\beta')^k 
  \left[
    \begin{array}{ccc}
      0, & 1, & \infty \\
      \hline
      r_0, & r_1, & r_\infty \\
      g_0, & g_1, & g_\infty
    \end{array}
    \right]_{n,k},
\end{aligned}
\end{equation}
where (in an extension of (\ref{eq:rewrite}) to arbitrary
nonzero $\beta,\beta'$)
\begin{equation}
\label{eq:old2new}
\begin{alignedat}{3}
  r_0 &= -\alpha/\beta,\qquad  & r_1 &= \alpha/\beta-\alpha'/\beta',\qquad  &r_\infty &= 1+\alpha'/\beta' = (\alpha'+\beta')/\beta',\\
  g_0 &= \gamma/\beta,\qquad  & g_1 &= -\gamma/\beta+\gamma'/\beta',\qquad  &g_\infty &= -\gamma'/\beta'.
\end{alignedat}
\end{equation}
Inverting these, one has that for any $(r_0,r_1,r_\infty)$ and
$(g_0,g_1,g_\infty)$ satisfying the conditions $r_0+r_1+r_\infty=1$
and $g_0+g_1+g_\infty=0$, and $\beta,\beta'$ satisfying
$\beta\beta'\neq0$, it is the case that
\begin{equation}
    \left[
    \begin{array}{ccc}
      0, & 1, & \infty \\
      \hline
      r_0, & r_1, & r_\infty \\
      g_0, & g_1, & g_\infty
    \end{array}
    \right]_{n,k}
    =
    \beta^{k-n}(-\beta')^{-k} 
  \left[
    \begin{array}{cc|c}
      \alpha, & \beta &  \gamma\\
      \alpha', & \beta' &  \gamma'
    \end{array}
    \right]_{n,k},
\end{equation}
and therefore
\begin{equation}
\label{eq:new2old}
    \left[
    \begin{array}{ccc}
      0, & 1, & \infty \\
      \hline
      r_0, & r_1, & r_\infty \\
      g_0, & g_1, & g_\infty
    \end{array}
    \right](t,z)
    =
  \left[
    \begin{array}{cc|c}
      \alpha, & \beta &  \gamma\\
      \alpha', & \beta' &  \gamma'
    \end{array}
    \right](-\beta t/\beta', z/\beta),
\end{equation}
where
\begin{equation}
\label{eq:new2olda}
  \begin{alignedat}{3}
    \alpha&=-\beta r_0,\qquad\quad &\alpha' &= -\beta'(r_0+r_1) = \beta'(r_\infty-1), \\
    \gamma&=\beta g_0,\qquad\quad &\gamma' &= \beta'(g_0+g_1) = -\beta'g_\infty.
  \end{alignedat}
\end{equation}
Formula~(\ref{eq:new2old}) expresses any new-style parametric EGF
in~terms of an old-style one.  The following theorem employs the new
notation but is equivalent to Theorem~\ref{thm:s3}.  That it does away
with the intricate parameter transformations of Table~\ref{tab:1}
justifies the new notation.

\begin{theorem}
\label{thm:unification}
  For each of the six rows in Table\/~{\rm\ref{tab:1}}, the corresponding
  lifting-based transformation of GKP triangles acts as follows on
  EGF's: the equality\/ $G^*(t^*,z^*) = G(t,z)$, where\/ $(t,z) =
  (R(t^*),S(t^*)z^*)$, can be written as
  \begin{equation}
    \left[
    \begin{array}{ccc}
      0, & 1, & \infty \\
      \hline
      r_0, & r_1, & r_\infty \\
      g_0, & g_1, & g_\infty
    \end{array}
    \right](R(t^*),S(t^*)z^*)
    =
    \left[
    \begin{array}{ccc}
      R^{-1}(0), & R^{-1}(1), & R^{-1}(\infty) \\
      \hline
      r_0, & r_1, & r_\infty \\
      g_0, & g_1, & g_\infty
    \end{array}
    \right](t^*,z^*),
  \end{equation}
where\/ $R^{-1}(0),R^{-1}(1),R^{-1}(\infty)$ is a permutation
of\/~$0,1,\infty$.  That is, in the new notation, each element of
the\/ $S_3$\nobreakdash-group acts as a permutation of the
parameter-pairs\/ $(r_0,g_0)$, $(r_1,g_1)$, $(r_\infty,g_\infty)$.
\end{theorem}

This is proved by rewriting the $(\beta,\beta')=(1,-1)$ case of each
of the parameter maps
$\left[\begin{smallarray}{cc|c}\alpha,&\beta&\gamma\\ \alpha',&\beta'&\gamma'\end{smallarray}\right]
\mapsto
\left[\begin{smallarray}{cc|c}\alpha,&\beta&\gamma\\ \alpha',&\beta'&\gamma'\end{smallarray}\right]^*$
of Table~\ref{tab:1} as a map from original parameter-pairs
$(r_0,g_0)$, $(r_1,g_1)$, $(r_\infty,g_\infty)$, to transformed (in
effect, lifted) parameter-pairs $(r_0,g_0)^*$, $(r_1,g_1)^*$,
$(r_\infty,g_\infty)^*$, with the aid of~(\ref{eq:rewrite}).

For instance, the RT (reflection transformation), specified by
$(0\infty)(1)\in S_3$ and previously written as~(\ref{eq:firsta}), can be
rewritten as
\begin{equation}
\label{eq:rtnew}
    \left[
    \begin{array}{ccc}
      0, & 1, & \infty \\
      \hline
      r_0, & r_1, & r_\infty \\
      g_0, & g_1, & g_\infty
    \end{array}
    \right]
    \left(\frac{1}{t^*},t^*z^*\right)
    =
    \left[
    \begin{array}{ccc}
      0, & 1, & \infty \\
      \hline
      r_\infty, & r_1, & r_0 \\
      g_\infty, & g_1, & g_0
    \end{array}
    \right]
    \left(t^*,z^*\right).
\end{equation}
The UBT\/ {\rm(}upper binomial transformation\/{\rm)}, specified by $(01)(\infty)\in
S_3$ and previously written as~(\ref{eq:seconda}), can be rewritten as
\begin{equation}
\label{eq:ubtnew}
    \left[
    \begin{array}{ccc}
      0, & 1, & \infty \\
      \hline
      r_0, & r_1, & r_\infty \\
      g_0, & g_1, & g_\infty
    \end{array}
    \right]
    \left(1-t^*,-z^*\right)
    =
    \left[
    \begin{array}{ccc}
      0, & 1, & \infty \\
      \hline
      r_1, & r_0, & r_\infty \\
      g_1, & g_0, & g_\infty
    \end{array}
    \right]
    \left(t^*,z^*\right).
\end{equation}
As before, these identities are valid in a neighborhood of
$(t^*,z^*)=(0,0)$ in~$\mathbb{C}\times\nobreak\mathbb{C}$, the trivial
case $t^*=0$ not being covered by~(\ref{eq:rtnew}).  Both
(\ref{eq:rtnew}) and~(\ref{eq:ubtnew}) are consistent with the
theorem, and because $(0\infty)(1)$ and $(01)(\infty)$ generate~$S_3$,
the theorem follows.

\subsection{Integrating the PDE}
\label{subsec:integrating}

The following theorem applies to any GKP triangle, parametrized as
explained in Section~\ref{subsec:newparam} by $r_0,r_1,r_\infty$ and
$g_0,g_1,g_\infty$ satisfying $r_0+r_1+r_\infty=1$ and
$g_0+g_1+g_\infty=0$.  It supplies a formula for the EGF which is
implicit rather than explicit, and is based upon a special function:
the Gauss hypergeometric function~${}_2F_1(w)$.  But in several cases
(see Section~\ref{subsec:cases}), the EGF can nonetheless be computed
in closed form.

The function ${}_2F_1(w)$ is parametric, with one lower and two upper
parameters.  Its Maclaurin series is
\begin{equation}
  {}_2F_1\left({\myatop{A,\, B}{C}}\biggm|w\right)   = 
\sum_{k=0}^\infty
    \frac{A^{\overline k}\,B^{\overline k}}{1^{\overline k}\,C^{\overline k}}\,
    w^k,
\end{equation}
and it is defined and analytic in a neighborhood of $w=0$, provided
that $C$~is not a non-positive integer.  In this and the following
subsection the alternative in-line notation ${}_2F_1\left(A,B;C\mid
w\right)$ will be used, for compactness of expressions.

\begin{theorem}
\label{thm:2f1formula}
  In a neighborhood of\/ $(t,z)=(0,0)$, at which it is analytic and
  equals unity, the EGF
\begin{equation}
G(t,z) \defeq
  \left[
    \begin{array}{ccc}
      0, & 1, & \infty \\
      \hline
      r_0, & r_1, & r_\infty \\
      g_0, & g_1, & g_\infty
    \end{array}
    \right](t,z)
\end{equation}
of a GKP triangle is given by the formula
\begin{equation}
  \label{eq:st}
  G(t,z) = \left(\frac{s}t\right)^{g_0} \left( \frac{1-s}{1-t} \right)^{g_1},
\end{equation}
where\/ $s=s(t,z)=t\left(1+zO(t,z)\right)$, with\/ $s(t,0)=t$, is defined
implicitly by
\begin{equation}
\label{eq:pleasant}
  \left(\frac{s}t\right)^{r_0} \left( \frac{1-s}{1-t} \right)^{r_1}
  \!=\,
  \frac{r_0z + {}_2F_1\left(r_0+r_1,1;\,1+r_0\bigm| t\right)}
       {{}_2F_1\left(r_0+r_1,1;\,1+r_0\bigm| s\right)}.
\end{equation}
This formula applies when\/ $r_0$ is not a non-positive integer.
\end{theorem}
\begin{remark}
  Though the parametric function ${}_2F_1(A,B;C;\cdot)$ is not defined
  when $C$~is a non-positive integer, the EGF of any triangle with
  $r_0=0,-1,-2,\dots$ can be computed by taking a limit.
  Alternatively, owing to the fact at~least one of $r_0,r_1,r_\infty$
  must not be a non-positive integer (as $r_0+r_1+r_\infty=1$), one
  can handle any case when $r_0=0,-1,-2,\dots$ by permuting the
  parameter-pairs $(r_0,g_0)$, $(r_1,g_1)$, $(r_\infty,g_\infty)$ with
  the aid of Theorem~\ref{thm:unification}, to obtain a triangle EGF
  which is covered by Theorem~\ref{thm:2f1formula}.
\end{remark}

\begin{proof}
  By~(\ref{eq:new2old}),
  \begin{equation}
    G(t,z) =
    \left[
      \begin{array}{ccc}
        0, & 1, & \infty \\
        \hline
        r_0, & r_1, & r_\infty \\
        g_0, & g_1, & g_\infty
      \end{array}
      \right](t,z)
    =
    \left[
      \begin{array}{ll|l}
        -r_0, & 1 & g_0 \\
        1-r_\infty, & -1 & g_\infty \\
      \end{array}
      \right](t,z),
  \end{equation}
  and by Theorem~\ref{thm:pde}, $G(t,z)$ satisfies the PDE
  \begin{equation}
    \label{eq:newpde}
    \Bigl\{\left[-r_0 + (1-r_\infty)t\right]z-1\Bigr\} \frac{\partial G}{\partial z} + (1-t)t \frac{\partial G}{\partial t} + (g_0 + g_\infty t) G = 0,
  \end{equation}
with the initial condition $G(t,0)=1$.  

To this first-order PDE, the method of characteristics can be applied.
For all $(t,z)$ in a neighborhood of~$(0,0)$
in~$\mathbb{C}\times\nobreak\mathbb{C}$, $G(t,z)$~can be computed by
flowing the initial condition $G=1$ from the $z=0$ line to the
point~$(t,z)$, along the characteristic curve extending to~$(t,z)$.
The Lagrange--Charpit equations coming from~(\ref{eq:newpde}) are
\begin{equation}
\label{eq:lc}
\left[
  \frac1{t(t-1)} - \left(
  \frac{r_0}t + \frac{r_1}{t-1}
  \right)z
  \right]^{-1}
     {\rm d}z
     ={\rm d}t
     = -\left(
     \frac{g_0}t + \frac{g_1}{t-1}
     \right)^{-1}\frac{{\rm d}G}G,
\end{equation}
and it is convenient to parametrize each characteristic by~$t$.  By
the equality between the second and third members, $G$~satisfies
$G\propto t^{-g_0}(1-t)^{-g_1}$ along each characteristic.  If
$s=s(t,z)$ denotes the value of~$t$ at which the characteristic
extending to~$(t,z)$ leaves the $z=0$ line, the initial condition
becomes $G(s(t,z),0)=1$ and can be imposed by expressing $G(t,z)$ as
in~(\ref{eq:st}).

It remains to find $s=s(t,z)$.  It will turn~out that $u\defeq s/t =
1+\nobreak zO(t,z)$, by an application of the implicit function theorem.  By
the equality between the first and second members of~(\ref{eq:lc}),
$z$~as a function of~$t$ along any characteristic satisfies the
inhomogeneous first-order~ODE
\begin{equation}
\label{eq:inhomg}
\frac{{\rm d}{z}}{{\rm d}t}
+
\left(
\frac{r_0}t + \frac{r_1}{t-1}\right)z = \frac1{t(t-1)}.
\end{equation}
The homogeneous solutions of~(\ref{eq:inhomg}) are of the form
$Kt^{-r_0}(1-\nobreak t)^{-r_1}$, where $K$~is arbitrary.
A~particular solution of~(\ref{eq:inhomg}) can be found by
differentiating (after multiplying by $t(t-1)$).  This yields the
homogeneous second-order~ODE
\begin{equation}
\label{eq:hypergode}
  \frac{{\rm d}^2 z}{{\rm d}t^2}
  + \left[ \frac{1+r_0}t + \frac{1+r_1}{t-1}\right]
  \frac{{\rm d} z}{{\rm d}t}
  +\, \frac{r_0+r_1}{t(t-1)}\,z = 0,
\end{equation}
which is a version of the Gauss hypergeometric equation.  Provided
that $r_0$~is not a negative integer, which holds by hypothesis, its
solutions in a neighborhood of $t=0$ include the analytic functions
\begin{equation}
\label{eq:2F1}
C\, {}_2F_1(1-\nobreak
r_\infty,1;\allowbreak1+\nobreak r_0\mid t),
\end{equation}
$C$~being arbitary.  Substitution reveals that $C$~must equal~$-1/r_0$
for (\ref{eq:2F1}) to be a solution of~(\ref{eq:inhomg}).  Combining
the particular solution~(\ref{eq:2F1}) with the homogeneous solution,
one finds that any characteristic must be of the form
\begin{equation}
  z(t) = K t^{-r_0}(1-t)^{-r_1} + (-1/r_0)\, {}_2F_1(r_0+r_1,1;1+r_0\mid
t)
\end{equation}
where $K$~specifies the characteristic.  If the characteristic is to
extend to~$(t,z)$ from~$(s,0)$, it must be the case that
\begin{equation}
\left\{
\begin{aligned}
  z &= K t^{-r_0}(1-t)^{-r_1} + (-1/r_0)\, {}_2F_1(r_0+r_1,1;1+r_0\mid t),\\
  0 &= K s^{-r_0}(1-s)^{-r_1} + (-1/r_0)\, {}_2F_1(r_0+r_1,1;1+r_0\mid s).
\end{aligned}
\right.
\end{equation}
This system determines the function $s=s(t,z)$.  Eliminating~$K$ yields
\begin{equation}
\label{eq:preimplicit}
\begin{aligned}
  &s^{r_0}(1-s)^{r_1}\, {}_2F_1\left(r_0+r_1,1;\,1+r_0\bigm| s\right)
  =\\
  &\qquad t^{r_0}(1-t)^{r_1}\left[r_0\, z +  {}_2F_1\left(r_0+r_1,1;\,1+r_0\bigm| t\right)\right],
\end{aligned}
\end{equation}
which as an equation for $u\defeq s/t$ can be written as
$\mathcal{F}(t,z;u)=0$, where
\begin{equation}
\label{eq:lesspleasant}
\begin{aligned}
  \mathcal{F}(t,z;\,u) &=u^{r_0}(1-tu)^{r_1}\, {}_2F_1\left(r_0+r_1,1;\,1+r_0\bigm| tu\right)\\
&\qquad {}- (1-t)^{r_1}\left[r_0\, z +  {}_2F_1\left(r_0+r_1,1;\,1+r_0\bigm| t\right)\right].
\end{aligned}
\end{equation}
Clearly $\mathcal{F}(0,0;1)=0$, and the question is whether
(\ref{eq:lesspleasant}) defines a function $u=u(t,z)$ which is
analytic in a neighborhood of $(t,z)=(0,0)$, at which point it equals
unity.  By direct computation $(\partial\mathcal{F}/\partial u)$
equals~$r_0$ when $(t,z;u)=(0,0;1)$, and $r_0\neq0$ by hypothesis; so
this follows by the analytic version of the implicit function theorem.
Equation~(\ref{eq:pleasant}) is a rewritten version
of~(\ref{eq:preimplicit}).
\end{proof}

\begin{remark}
  The Gauss hypergeometric ODE~(\ref{eq:hypergode}) would be written
  in Riemann's P-symbol notation as
  \begin{equation}
    \label{eq:Psymbol}
    {\rm P}
    \left[
      \begin{array}{ccc}
        0, & 1, & \infty\\
        \hline
        0, & 0, & 1 \\
        -r_0, & -r_1, & 1-r_\infty
      \end{array}
      \Biggm|
      t
      \right],
  \end{equation}
  which lists the two characteristic exponents of each of its singular
  points (the points $0,1,\infty$ on the projective
  $t$\nobreakdash-line).  (The exponent \emph{differences}, important
  in the construction of local series solutions, are
  $-r_0,-r_1,-r_\infty$.)  By Fuchs's relation, the sum of the six
  exponents equals unity.  The new GKP parametrization introduced in
  Section~\ref{subsec:newparam} was suggested by~(\ref{eq:Psymbol}).
  However, it~also lists the parameters $g_0,g_1,g_\infty$, which
  do~not appear in the ODE though they appear in the
  PDE~(\ref{eq:newpde}), from which the ODE was derived.
\end{remark}

\subsection{Some special cases}
\label{subsec:cases}

The EGF of a GKP triangle, parametrized by $r_0,r_1,r_\infty$ and
$g_0,g_1,g_\infty$ with respective sums $1$~and~$0$, can in some
interesting cases be computed in closed form from
Theorem~\ref{thm:2f1formula}.  These include (A)~when one of
$r_0,r_1,r_\infty$ equals~$1$; (B)~when
$\{r_0,r_1,r_\infty\}=\{-\frac12,-\frac12,2\}$; and (C)~when
$\{r_0,r_1,r_\infty\}=\{\frac12,\frac12,0\}$.  In the following
sections these are related to the generalized Stirling--Eulerian,
Narayana, and secant--tangent triangles.  In each of (A),(B),(C),
$r_0,r_1,r_\infty$ can be permuted with the aid of
Theorem~\ref{thm:unification}.  Hence without loss of generality, it
suffices to examine the cases (A\,I)~when $r_\infty$ equals~$1$;
(B\,I)~when $(r_0,r_1,r_\infty)=(-\frac12,-\frac12,2)$; and
(C\,I)~when $(r_0,r_1,r_\infty)=(\frac12,\frac12,0)$.

The following three theorems evaluate
  \begin{equation}
    G(t,z) =
    \left[
      \begin{array}{ccc}
        0, & 1, & \infty \\
        \hline
        r_0, & r_1, & r_\infty \\
        g_0, & g_1, & g_\infty
      \end{array}
      \right](t,z)
\end{equation}
in cases (A\,I), (B\,I), and~(C\,I).  The first two are especially
easy to treat because when an upper parameter of the hypergeometric
function~${}_2F_1$ equals a non-positive integer~$-N$, the power series
defining the function terminates and becomes a degree\nobreakdash-$N$
polynomial.

\begin{theorem} 
\label{thm:AI}
$\mathrm{(A\,I)}$ If\/ $r_\infty=1$ {\rm(}so that\/ $r_0+r_1=0${\rm)},
the EGF is given in a neighborhood of\/ $(0,0)$ by
  \begin{equation}
    \label{eq:AI}
    \begin{aligned}
    G(t,z)&=
    \left[t+(1-t)(1+r_0z)^{-1/r_0}
      \right]^{-g_0}
    \left[(1-t) + t(1+r_0z)^{1/r_0}
      \right]^{-g_1}\\
    &= \left[(1+r_0 z)^{1/r_0}\right]^{g_0}\left[(1-t)+t (1+r_0z)^{1/r_0}\right]^{g_\infty}
    \end{aligned}
  \end{equation}
when\/ $r_0\neq0$.
\end{theorem}
\begin{proof}
  If $r_0$~is not a non-positive integer, formula~(\ref{eq:pleasant})
  of Theorem~\ref{thm:2f1formula} applies, and if $r_\infty=1$, each
  ${}_2F_1$ in the formula degenerates to the unit (constant)
  function.  Some algebra then yields
  \begin{equation}
    s=s(t,z) = \frac{t(1+r_0z)^{1/r_0}}{(1-t) + t(1+r_0z)^{1/r_0}},
  \end{equation}
and~(\ref{eq:AI})  follows from~(\ref{eq:st}).  If $r_0$~is a
negative integer (though not if $r_0=0$), (\ref{eq:AI}) holds by a
limit argument.
\end{proof}

The value $r_0=0$ not covered by the theorem is handled thus: If
$(r_0,r_1,r_\infty)=(0,0,1)$ then
\begin{equation}
\begin{aligned}
  G(t,z) &= 
    \left[t+(1-t){\rm e}^{-z}
      \right]^{-g_0}
    \left[(1-t) + t{\rm e}^z
      \right]^{-g_1}\\
    &= {\rm e}^{g_0 z} \left[(1-t)+t{\rm e}^z\right]^{g_\infty}
    ,
\end{aligned}
\end{equation}
by taking $r_0\to0$.

\begin{theorem} 
\label{thm:BI}
$\mathrm{(B\,I)}$
  If\/ $(r_0,r_1,r_\infty)=(-\frac12,-\frac12,2)$, the EGF is given in a neighborhood of\/ $(0,0)$ by
  \begin{equation}
    G(t,z) = \left(\frac{s_+}{t_+}\right)^{g_0}
    \left(\frac{s_-}{t_-}\right)^{g_1},
    \label{eq:lastminG}
  \end{equation}
  where
  \begin{displaymath}
    s_\pm = \frac12 \pm \frac{4(t-\tfrac12)+z}{2\sqrt{4+8(t-\tfrac12)z+z^2}},
  \end{displaymath}
  with\/ $s_++s_-=1$, and $t_+=t$ and $t_-=1-t$.
\end{theorem}
\begin{proof}
  As $r_0+r_1=-1$, each ${}_2F_1$ in~(\ref{eq:pleasant}) is a
  degree\nobreakdash-$1$ polynomial function of its argument, and
  (\ref{eq:pleasant}) becomes a quadratic equation for~$s$, the
  solutions of which are $s_+$ and~$s_-$.  The one with the correct
  behavior as $t,z\to0$, satisfying $s=t(1+\nobreak zO(t,z))$,
  is~$s_+$, so in~(\ref{eq:st}), $s$~and $1-\nobreak s$ are
  respectively equal to $s_+$ and~$s_-$.  One can confirm that
  (\ref{eq:lastminG})~satisfies the PDE~(\ref{eq:newpde}).
\end{proof}

It is clear from formula~(\ref{eq:pleasant}) that if $r_\infty$~is a
positive integer and $r_0,r_1$ are nonzero rational numbers,
$s=s(t,z)$ will be an algebraic function, and if moreover
$g_0,g_1,g_\infty\in\mathbb{Q}$, the same will be true of~$G(t,z)$.
But the polynomial of which $s$~is a root will typically be of higher
degree than quadratic.  The following theorem deals with an inherently
non-algebraic case~($r_\infty=0$).

\begin{theorem} 
\label{thm:CI}
$\mathrm{(C\,I)}$
  If\/ $(r_0,r_1,r_\infty)=(\frac12,\frac12,0)$, the EGF is given in a neighborhood of\/ $(0,0)$ by
  \begin{equation}
    \label{eq:lastminuteG}
    G(t,z) = \left(\frac{s_+}{t_+}\right)^{g_0}
    \left(\frac{s_-}{t_-}\right)^{g_1},
  \end{equation}
  where
  \begin{displaymath}
    s_\pm = \left[\sqrt{t_\pm}\cos\left( \frac{z}2 \sqrt{t_+t_-}\right)
    \pm
    \sqrt{t_\mp}\sin\left( \frac{z}2 \sqrt{t_+t_-}\right)\right]^2,
  \end{displaymath}
  with\/ $s_++s_-=1$, and $t_+=t$ and $t_-=1-t$.
\end{theorem}
\begin{proof}
  This comes from the known fact that in a neighborhood of~$t=0$,
  ${}_2F_1\left(1,1;\nobreak\tfrac32;\nobreak t\right)$ equals
  $\sin^{-1}(\sqrt{t}\,)/\sqrt{t(1-t)}$.  Squaring both sides
  of~(\ref{eq:pleasant}), one sees that $s=s(t,z)$ is defined
  implicitly by
  \begin{equation}
    \frac{s(1-s)}{t(1-t)} =
    \left[
      \frac{(z/2) + \sin^{-1}(\sqrt{t}\,) / \sqrt{t(1-t)}}
           {\sin^{-1}(\sqrt{s}\,) / \sqrt{s(1-s)}}
      \right]^2,
  \end{equation}
  which simplifies to the statement that
  \begin{equation}
    s = \sin^2 \left(\frac{z}2\sqrt{t(1-t)} + \sin^{-1}(\sqrt{t}\,)\right),
  \end{equation}
  or equivalently to~$s=s_+$.  By examination, $1-s$ equals $s_-$, and
  (\ref{eq:lastminuteG}) follows from~(\ref{eq:st}).  One can confirm
  that (\ref{eq:lastminuteG}) satisfies the PDE~(\ref{eq:newpde}).
\end{proof}

Cases (A\,I),(B\,I),(C\,I) can be converted to what will be called
(A\,II),\allowbreak(B\,II),\allowbreak(C\,II) and
(A\,III),\allowbreak(B\,III),\allowbreak(C\,III) by applying
respectively the elements $(0\infty)(1)$ and $(1\infty)(0)$ of the
$S_3$\nobreakdash-group, i.e., the sequence transformations
$\textrm{RT}$ and $\textrm{RT}\circ\textrm{UBT}\circ\textrm{RT}$.  The
resulting EGF formulas will appear in the following three sections but
can be summarized as~follows.
\begin{theorem}
\label{thm:38}
  The statements of Theorems\/ {\rm\ref{thm:AI}}, {\rm\ref{thm:BI}},
  and {\rm\ref{thm:CI}} remain valid if the parameter pairs\/
  $(r_0,g_0)$ and\/ $(r_\infty,g_\infty)$ are interchanged, with\/ $t$
  replaced by\/ $\frac1t$ and\/ $z$ by\/ $tz$; and similarly if\/
  $(r_1,g_1)$ and\/ $(r_\infty,g_\infty)$ are interchanged, with\/ $t$
  replaced by\/ $\frac{-t}{1-t}$ and\/ $z$ by\/ $(1-t)z$.
\end{theorem}

However, each of the three theorems is unchanged (or is unchanged
up~to parametrization, in the case of Theorem~\ref{thm:AI}) by the
remaining involution $(01)(\infty)$, i.e., the~UBT\null.

Barbero~G. et~al.~\cite{Barbero2014}, besides computing the bivariate
EGF $G(t,z)$ in the three subcases of case~(A), have treated the case
when (in~present notation) the unordered set $\{r_0,r_1,r_\infty\}$
equals $\{N,1-N,0\}$, for some $N\in\mathbb{Z}\setminus\{0,1\}$.
(See~\cite{Barbero2014}, \S\S{A.1.3}, {A.1.5}, {A.1.6}.)  This case
also has three subcases, which are related by the
$S_3$\nobreakdash-group.  But in each, the EGF turns~out not to be an
elementary function, but rather to be expressible in~terms of an
implicitly defined tree function (of~combinatorial significance).
Analytically, this can be attributed in~part to the ${}_2F_1$'s
in~(\ref{eq:pleasant}) not being elementary functions.

\section{Generalized Stirling--Eulerian triangles}
\label{sec:stirlingeulerian}

In GKP case (A\,I) of the last section, when $r_\infty=1$ or
equivalently $\alpha'=0$, Theorem~\ref{thm:AI} supplies a closed-form
expression for the EGF $G(t,z)$.  This case leads naturally to the
definition of the generalized Stirling and generalized Eulerian
numbers, $S_{n,k}(a,b;r)$ and $E_{n,k}(a,b;c_0,c_\infty)$.  The latter
are characterized by the similar condition $r_1=1$, and
in~\S\ref{subsec:32}, they are alternatively interpreted as connection
coefficients.  (The key result is Theorem~\ref{thm:addfactsE}, which
also includes a rank\nobreakdash-$1$ formula for
$E_{n,k}(a,b;c_0,c_\infty)$.)  In
\S\S\ref{subsec:33}~and~\ref{subsec:34}, some important parameter
choices when `rank\nobreakdash-$0$' formulas for $S_{n,k}(a,b;r)$ and
$E_{n,k}(a,b;c_0,c_\infty)$ exist are examined.  These may contain
binomial coefficients and generalized factorials, without much
summation.

\subsection{Basic formulas}
\label{subsec:31}

Applying the transformations $\textrm{RT}$ and
$\textrm{RT}\circ\textrm{UBT}\circ\textrm{RT}$ to case (A\,I), as
summarized in Theorem~\ref{thm:38}, yields the two additional EGF
formulas that appear in the following three theorems.
(Case~(A\,II) is when $\alpha=-\beta$ and case (A\,III) is when
$\frac\alpha\beta =\allowbreak \frac{\alpha'}{\beta'}+\nobreak1$; as
above, $\beta\beta'\neq0$ is assumed.)  The parameters
$(\alpha,\beta,\gamma;\allowbreak\alpha',\beta',\gamma')$ of each EGF
have been computed from the new parameters
$(r_0,r_1,r_\infty;\allowbreak g_0,g_1,g_\infty)$ with the aid of
(\ref{eq:new2old}) and~(\ref{eq:new2olda}).

These three EGF formulas have been derived previously.
(See~\cite[(A.8),\allowbreak(A.4),\allowbreak(A.2)]{Barbero2014}, and
also~\cite{Neuwirth2001,Theoret94,Theoret95}.)  But the present
derivation, making explicit use of the $S_3$\nobreakdash-group of
sequence transformations to derive the latter two from the first,
seems the most efficient.

\begin{theorem}[(A\,I), $r_\infty=1$, $\alpha'=0$: generalized Stirling]
\label{thm:AInew}
  If\/ $\alpha'=0$ then\/ $G(t,z)$ equals\/ {\rm(}when\/ $\alpha\neq0${\rm)}
  \begin{displaymath}
    (1-\alpha z)^{-\gamma/\alpha}
    \left\{
    1 + (\beta'/\beta)t\left[1- (1-\alpha z)^{-\beta/\alpha} \right]
    \right\}^{-\gamma'/\beta'},
  \end{displaymath}
  which in the\/ $\alpha\to0$ limit becomes
  \begin{displaymath}
    {\rm e}^{\gamma z}
    \left\{
    1 + (\beta'/\beta)t\left[1- {\rm e}^{\beta z} \right]
    \right\}^{-\gamma'/\beta'}.
  \end{displaymath}
\end{theorem}

\begin{theorem}[(A\,II), $r_0=1$, $\alpha=-\beta$: generalized Stirling, reflected]
\label{thm:AII}
  If\/ $\alpha=-\beta$ then\/ $G(t,z)$ equals\/ {\rm(}when $\alpha'+\beta'\neq0${\rm)}
  \begin{displaymath}
    \left[
      1-(\alpha'+\beta')zt
      \right]^{-\gamma'/(\alpha'+\beta')}
    \left\{
    1 + (\beta/\beta')t^{-1}
    \left[
      1-\left(
      1-   (\alpha'+\beta')zt \right)^{\beta'/(\alpha'+\beta')}
      \right]
    \right\}^{\gamma/\beta},
  \end{displaymath}
  which in the\/ $\alpha'\to-\beta'$ limit becomes
  \begin{displaymath}
    {\rm e}^{\gamma'zt}
    \left\{
    1 + (\beta/\beta')t^{-1}
    \left[
      1-
      {\rm e}^{-\beta' zt}
      \right]
    \right\}^{\gamma/\beta}.
  \end{displaymath}
\end{theorem}

\begin{theorem}[(A\,III), $r_1=1$, $\frac{\alpha}{\beta}=\frac{\alpha'}{\beta'}+1$: generalized Eulerian]
\label{thm:AIII}
  If\/ $\frac{\alpha}{\beta}=\frac{\alpha'}{\beta'}+1$ then\/ $G(t,z)$ equals\/
  {\rm(}when $\alpha\neq0$, i.e., $\alpha'+\beta'\neq0${\rm)}
  \begin{align*}
    &\left\{ 
    1-[\beta/(\beta+\beta't)]\left[1-\bigl(1-\alpha z\bigm/[\beta/(\beta+\beta't)]\bigr)^{\beta/\alpha}\right]
    \right\}^{-\gamma/\beta}
    \\
    &\qquad{}\times\left\{ 
    1-[\beta't/(\beta+\beta't)]\left[1-\bigl(1-\alpha z\bigm/[\beta/(\beta+\beta't)]\bigr)^{-\beta/\alpha}\right]
    \right\}^{\gamma'/\beta'},
  \end{align*}
  which in the\/ $\alpha\to0$ or equivalently\/ $\alpha'\to-\beta'$ limit becomes
  \begin{align*}
    &\left\{ 
    1-[\beta/(\beta+\beta't)]\left[1- {\rm e}^{-z(\beta+\beta't)}\right]
    \right\}^{-\gamma/\beta}
    \\
    &\qquad{}\times
    \left\{ 
    1-[\beta't/(\beta+\beta't)]\left[1- {\rm e}^{z(\beta+\beta't)}\right]
    \right\}^{\gamma'/\beta'}.
  \end{align*}
\end{theorem}

Cases (A\,I) and~(A\,II) are related by~$\textrm{RT}$, the action
of which is straightforward.  (It reverses each row of a GKP triangle;
see the second row in Table~\ref{tab:1}.)  The focus will therefore be
on (A\,I) and~(A\,III)\null.  In both, it is natural to reduce the
number of free parameters by defining a specialized or normalized
version of the GKP triangle.  Case (A\,I) will be treated first.
\begin{definition}
\label{def:Sdef}
The $3$-parameter generalized Stirling triangle $S_{n,k}(a,b;r)$ is defined by
\begin{equation}
\label{eq:hsushiue3}
S_{n,k} = S_{n,k}(a,b;\,r) \defeq
\left[
  \begin{array}{cc|c}
    -a, & b & r \\
    0, & 0 & 1
  \end{array}
  \right]_{n,k}
\end{equation}
and satisfies $S_{n+1,k+1} = [-an+b(k+1)+r]S_{n,k+1} + S_{n,k}$.
\end{definition}
The corresponding denormalization is
\begin{equation}
\label{eq:denorm}
  \left[
  \begin{array}{cc|c}
    \alpha, & \beta & \gamma \\
    0, & \beta' & \gamma'
  \end{array}
\right]_{n,k}
=
\left(\frac{\gamma'}{\beta'}\right)^{\overline k}(\beta')^k
  \left[
  \begin{array}{cc|c}
    \alpha, & \beta & \gamma \\
    0, & 0 & 1
  \end{array}
\right]_{n,k}
= (\gamma')^{\overline k, \beta'} S_{n,k}(-\alpha,\beta;\,\gamma)
\end{equation}
and a homogeneity property is
\begin{equation}
\label{eq:Shomog}
S_{n,k}(\lambda a,\lambda b;\, \lambda r) =  \lambda ^{n-k}S_{n,k}(a,b;\,r).  
\end{equation}
The numbers $S_{n,k}(a,b;r)$ are of~course the generalized Stirling
numbers of Hsu and Shiue, which reduce to $\tstirsub{n}k$
and~$\tstircyc{n}k$ when $(a,b;r)=(0,1;0)$, resp.\ $(-1,0;0)$.  When
restricted to integer parameter values, they have been interpreted
combinatorially \cite{Corcino2001a,Maltenfort2020}.  If $b\neq0$ and
$a\neq0$, these numbers have the bivariate~EGF
\begin{equation}
\label{eq:elem}
  \sum_{n=0}^\infty \sum_{k=0}^n k!\,S_{n,k}(a,b;\,r)\, t^k\,\frac{z^n}{n!}
  = (1+az)^{r/a} \left\{1-\frac{t}b\left[(1+az)^{b/a}-1\right]\right\}^{-1},
\end{equation}
which follows from (\ref{eq:denorm}) and
Theorem~\ref{thm:AInew} by choosing $\beta'=\gamma'=1$, and the
equivalent but perhaps less useful~EGF
\begin{equation}
  \sum_{n=0}^\infty \sum_{k=0}^n S_{n,k}(a,b;\,r)\,t^k\,\frac{z^n}{n!}
  = (1+az)^{r/a}\exp \left\{\frac{t}b\left[(1+az)^{b/a}-1\right]\right\},
\end{equation}
which follows from Theorem~\ref{thm:AInew} by taking $\beta'\to0$.
Taking account of the elementary $t$-dependence in~(\ref{eq:elem}),
one has for all $k\ge0$ the `vertical' univariate~EGF
\begin{equation}
\label{eq:vertical}
  \sum_{n=0}^\infty S_{n,k}(a,b;\,r) \frac{k!}{n!}\,z^n
  =(1+az)^{r/a} \left[ \frac{(1+az)^{b/a}-1}b \right]^k,
\end{equation}
which is of the simple form $d(z) h(z)^k$; so irrespective of the
choice of parameters, $\frac{k!}{n!}S_{n,k}$~is a Riordan array, and
$S_{n,k}$ itself is a so-called exponential Riordan
array~\cite{Barry2016}.  Taking the $a\to0$ limit in any of the three
preceding formulas is straightforward (the $b\to0$ limit is not
considered here).

The following additional facts are well known but are proved for
completeness.  In this, $\Delta_x$~is the forward first difference
operator with respect to~$x$, defined by $\Delta_xf(x)=\allowbreak
f(x+\nobreak 1)-\nobreak f(x)$.

\begin{theorem}
\label{thm:addfacts}
  {\rm(i)} When\/ $b\neq0$, $S_{n,k}(a,b;r)$ is given by the
  rank\/\nobreakdash-$1$ formula
  \begin{subequations}
    \begin{align}
      \label{eq:expanded}
      S_{n,k}(a,b;\,r) &= 
      \frac{1}{b^kk!}\sum_{j=0}^k (-1)^{k-j} \binom{k}{j}(bj+r)^{\underline{n},a}\\
      \label{eq:preexpanded}
      &= \frac1{b^kk!}\, \Delta_x^k\bigl[ (bx+r)^{\underline{n},a}\bigr]\bigm|_{x=0}.
    \end{align}
  \end{subequations}
  \begin{itemize}
  \item[{\rm(ii)}] In general the numbers\/ $S_{n,k}(a,b;r)$, $0\le k\le n<\infty$, satisfy
    \begin{equation}
      \label{eq:certainpolybases}
      (x)^{\underline{n},a} = \sum_{k=0}^n S_{n,k}(a,b;\,r) (x-r)^{\underline k,b},
    \end{equation}
  \end{itemize}
  defining them as coefficients of connection between certain graded
  bases of the space of polynomials in an indeterminate\/~$x$.
\end{theorem}
\begin{proof}
  The exponential Riordan array $\tgenvert{n}k$ specified by a pair of
  formal power series $d,h$ in~$z$ of respective orders~$0,1$, i.e.,
  $d(z)=d_0+d_1z+\dots$ with $d_0\neq0$ and $h(z)=h_1z+h_2z^2+\dots$
  with $h_1\neq0$, denoted by $\mathrm{e}\mathcal{R}(d,h)$ or~$[d,h]$,
  is the infinite lower-triangular matrix defined by
  \begin{equation}
    \label{eq:mvprod}
    [d,h]_{n,k}=\genvert{n}k =\frac{n!}{k!} [z^n] d(z)h(z)^k,
  \end{equation}
  where $[z^n]$ extracts the coefficient of~$z^n$.  It is a
  fundamental fact~\cite{Barry2016} that if a column vector~$u_k$ has
  EGF $u(t)=\sum_{k=0}^\infty u_kt^k/k!$, the matrix--vector product
  $v_n=\sum_{k=0}^\infty \tgenvert{n}k u_k$ will have EGF
  $v(z)=\sum_{n=0}^\infty v_nz^n/n!$ equal to $d(z)u(h(z))$.
  In~(\ref{eq:certainpolybases}), the column vector
  $u_k=(x-r)^{\underline k,b}$ has EGF $u(t)=(1+bt)^{(x-r)/b}$ and
  $v_n=(x)^{\underline n,a}$ has EGF $v(z)=(1+az)^{x/a}$, consistent
  with this fact and the functions $d(z),h(z)$ appearing
  in~(\ref{eq:vertical}).  This proves~(\ref{eq:certainpolybases}).

  From~(\ref{eq:certainpolybases}) with $x$~replaced by $bx+r$, it
  follows by Newton's interpolation formula~\cite{Miller66} that when
  $b\neq0$, (\ref{eq:preexpanded})~holds.
  Equation~(\ref{eq:expanded}) is an expanded version
  of~(\ref{eq:preexpanded}).
\end{proof}

Case (A\,III), when $\frac{\alpha}{\beta}=\frac{\alpha'}{\beta'}+1$, can also be specialized or normalized to reduce the number of free parameters.
\begin{definition}
\label{def:Edef}
The $4$-parameter generalized Eulerian triangle
$E_{n,k}(a,b;c_0,c_\infty)$ is defined by
\begin{equation}
\label{eq:euleriannew1}
  E_{n,k}=E_{n,k}(a,b;\,c_0,c_\infty) \defeq
  \left[
    \begin{array}{cc|c}
      -a, & b & c_0\\
      a+b,& -b & c_\infty
    \end{array}
    \right]_{n,k},
\end{equation}
which if $b\neq0$ comes equivalently from a tableau with $r_1=1$, as
\begin{equation}
  b^n
  \left[
    \begin{array}{ccc}
      0, & 1, & \infty \\
      \hline
      a/b, & 1, & -a/b \\
      c_0/b, & \,-(c_0+c_\infty)/b,\, & c_\infty/b
    \end{array}
    \right]_{n,k}.
\end{equation}
It satisfies $E_{n+1,k+1} = [-an+b(k+1)+c_0]E_{n,k+1} + [(a+b)n - bk +
  c_\infty]E_{n,k}$.  
\end{definition}
The corresponding denormalization is
\begin{equation}
  \left[
    \begin{array}{cc|c}
      \alpha, & \beta & \gamma \\
      \left(\frac{\alpha}{\beta}-1\right)\beta', & \beta' & \gamma' \\
    \end{array}
    \right]_{n,k}  
  =
  \left(
  \frac{-\beta'}{\beta}
  \right)^k
  E_{n,k}(-\alpha,\beta;\,\gamma,-\beta\gamma'/\beta')
\end{equation}
and a homogeneity property is
\begin{equation}
\label{eq:Ehomog}
E_{n,k}( \lambda a, \lambda b;\, \lambda c_0, \lambda c_\infty) =  \lambda ^nE_{n,k}(a,b;\,c_0,c_\infty).  
\end{equation}
The generalized Eulerian numbers $E_{n,k}(a,b;c_0,c_\infty)$ reduce
to the standard numbers $\eulerian{n}{k-1}\defeq A_{n,k}$,
$\eulerian{n}{k}$, and $\eulerian{n+1}k = A_{n+1,k+1}$ when
$(a,b;c_0,c_\infty)$ is respectively equal to $(0,1;0,1)$,
$(0,1;1,0)$, and $(0,1;1,1)$; recall Example~\ref{ex:22}.

It is noteworthy that when $b\neq0$, the $E_{n,k}(a,b;c_0,c_\infty)$
or the corresponding row polynomials $G_n(t)=\allowbreak \sum_{k=0}^n
E_{n,k}(a,b;c_0,c_\infty) t^k$, $n\ge0$, can be computed by repeated
differentiation.  According to the formula of
Theorem~\ref{thm:pde}(iii),
\begin{equation}
\label{eq:euleriandiffrec}
  (bt^{1-\hat a}D_t)^n
  \left\{
  \frac{t^{\hat c_0}}{(1-t)^{\hat c_0+\hat c_\infty}}
  \right\}
  =
  \frac{t^{\hat c_0-\hat a n}\, G_n(t)}{(1-t)^{\hat c_0+\hat c_\infty  + n}},
\end{equation}
or equivalently
\begin{equation}
\label{eq:euleriandiffrec2}
b^n\sum_{k=0}^\infty (\hat c_0 + k)^{\underline n,\hat a} (\hat c_0+\hat
c_\infty)^{\overline k}\,\frac{t^k}{k!}
=
\frac{G_n(t)}{(1-t)^{\hat c_0 + \hat c_\infty + n}},
\end{equation}
where $\hat a=a/b$, $\hat c_0= c_0/b$, $\hat c_\infty= c_\infty/b$.
These formulas subsume Euler's ones for the numbers $A_{n,k}$, to
which they reduce when $(a,b;c_0,c_\infty)$ equals $(0,1;0,1)$.

Combinatorial interpretations of the numbers
$E_{n,k}(a,b;\,c_0,c_\infty)$ for general parameter values remain to
be explored, though as mentioned in the introduction, some special
cases (such as the Carlitz--Scoville $a=0$ case) have appeared in the
literature.  When $b\neq0$ and $a\neq0$ these numbers have the
bivariate~EGF
\begin{gather}
\label{eq:speck}
\begin{aligned}
  &\sum_{n=0}^\infty \sum_{k=0}^n E_{n,k}(a,b;\,c_0,c_\infty)
  \,t^k\,\frac{z^n}{n!}=\\
  &\qquad\qquad
  \left\{
  1-(1-t)^{-1}\left[1-(1+az-atz)^{-b/a}\right]
  \right\}^{-c_0/b}\\
  &\qquad\qquad{}\times
  \left\{
  1+t(1-t)^{-1}\left[1-(1+az-atz)^{b/a}\right]
  \right\}^{-c_\infty/b},
\end{aligned}
\\
\intertext{which becomes}
\label{eq:speck2}
\begin{aligned}
  &\sum_{n=0}^\infty \sum_{k=0}^n E_{n,k}(0,b;\,c_0,c_\infty)
  \,t^k\,\frac{z^n}{n!}=\\
  &\quad
  \left\{
  1-(1-t)^{-1}
  \left[
    1-{\rm e}^{-b(1-t)z}
    \right]
  \right\}^{-c_0/b}
  \left\{
  1+t(1-t)^{-1}\left[
    1-{\rm e}^{b(1-t)z}
    \right]
  \right\}^{-c_\infty/b}
\end{aligned}
\end{gather}
when $a\to0$ (the $b\to0$ limit is not considered here).  The EGF's
(\ref{eq:speck}),(\ref{eq:speck2}) are specializations of the EGF's of
Theorem~\ref{thm:AIII}.  Note that by setting $b=1$ with
$(c_0,c_\infty)=(0,1)$, $(1,0)$, and $(1,1)$, one recovers
from~(\ref{eq:speck2}) the classical Eulerian EGF's of
Example~\ref{ex:22}.  More generally, by setting $b=1$ with
$(c_0,c_\infty)=(u,v)$, one obtains an EGF for the $(u,v)$-Eulerian
numbers of Barbero~G. et~al.~\cite{Barbero2015}.  These include the
traditional order\nobreakdash-$v$ Eulerian numbers $A_{n,k}^{(v)}$
defined in~\cite{Dillon68}, for which $(c_0,c_\infty)=(0,v)$.

The analogue of Theorem~\ref{thm:addfacts} for the new numbers
$E_{n,k}(a,b;\,c_0,c_\infty)$, including a rank\nobreakdash-$1$
formula, will appear as Theorem~\ref{thm:addfactsE} below.  These
numbers can also be expressed in~terms of the Hsu--Shiue numbers
$S_{n,k}(a,b;\,r)$, and vice versa, by applying appropriate elements
of the $S_3$\nobreakdash-group of sequence transformations acting
row-wise:
\begin{theorem}
\label{thm:ubtpair}
  For all\/ $n\ge0$ and choices of parameters\/ $(a,b;c_0,c_\infty)$, one
  has the UBT\/ {\rm(}upper binomial transform\/{\rm)} pair
  \begin{subequations}
    \begin{align}
      & E_{n,k}(a,b;\,c_0,c_\infty) = \sum_{j=k}^n (-1)^{j-k} \binom{j}{k}
       (c_0+c_\infty)^{\overline{n-j},b} \,S_{n,n-j}(-a,b;\,c_\infty),\label{eq:ubt1}\\
        &(c_0+c_\infty)^{\overline{n-k},b}\,S_{n,n-k}(-a,b;\,c_\infty)
          = \sum_{j=k}^n\binom{j}{k} E_{n,j}(a,b;\,c_0,c_\infty).\label{eq:ubt2}
    \end{align}
  \end{subequations}
\end{theorem}
\begin{proof}
The first identity is an application of the composite sequence
transformation $\textrm{UBT}\circ\textrm{RT}$, which is the
order\nobreakdash-$3$ transformation listed on the fifth line of
Table~\ref{tab:1}.  Equation~(\ref{eq:ubt1}) can be written as
\begin{equation}
    \left[
    \begin{array}{ll|l}
      -a, & b & c_0 \\
      a+b, & -b & c_\infty
    \end{array}
    \right]_{n,k} = (-1)^{n-k}
    \sum_{j=k}^n \binom{j}{k} \hat{S}_{n,j},
\end{equation}
where
\begin{equation}
  \begin{aligned}
    \hat{S}_{n,j} &= 
    \left\{
    (-1)^j (c_0+c_\infty)^{\underline{j},b}
    \left[
      \begin{array}{ll|l}
        a, & b & c_\infty \\
        0, & 0 & 1
      \end{array}
      \right]_{n,j}
    \right\}_{j\leftarrow n-j}    
    \\
    &=
    \left\{
    \left[
      \begin{array}{ll|c}
        a, & b & c_\infty \\
        0, & -b & -c_0-c_\infty
      \end{array}
      \right]_{n,j}
    \right\}_{j\leftarrow n-j}
    =
    \left[
      \begin{array}{ll|c}
        -b, & b & -c_0-c_\infty \\
        a+b, & -b & c_\infty
      \end{array}
      \right]_{n,j},
  \end{aligned}
\end{equation}
in which the final equality comes from the RT
formula~(\ref{eq:firstc}).  The identity now follows by applying the
UBT formula~(\ref{eq:secondc}).  The second identity is inverse to the
first, and comes from the inverted group element
$\textrm{RT}\circ\textrm{UBT}$.
\end{proof}

When $(a,b;c_0,c_\infty)=(0,1;1,0)$, the UBT pair in this theorem
(incorporating an initial $j\leftarrow n-\nobreak j$ reflection,
resp.\ a final $k\leftarrow n-\nobreak k$ reflection) reduces to the
classical UBT pair~\cite{Graham94} relating the Eulerian numbers
$E_{n,k}(0,1;1,0)=\eulerian{n}k$ and the Stirling subset numbers
$S_{n,k}(0,1;0)=\stirsub{n}k$:
\begin{subequations}
\label{eq:overallpair}
  \begin{align}
    &\eulerian{n}k = \sum_{j=k}^n (-1)^{j-k}\binom{j}k (n-j)!\, \stirsub{n}{n-j},\label{eq:ubtspec1}\\
    &(n-k)!\,\stirsub{n}{n-k} = \sum_{j=k}^n \binom{j}k  \eulerian{n}j.\label{eq:ubtspec2}
  \end{align}
\end{subequations}
The sequence transformation performed in (\ref{eq:ubtspec2}) is an
(unsigned and non-involutive) UBT, and the transformation
in~(\ref{eq:ubtspec1}) is its inverse.

Slightly modified versions of the UBT pair~(\ref{eq:overallpair}) can
be obtained by setting $(a,b;\allowbreak c_0,c_\infty)$ equal to
$(0,1;0,1)$ and $(0,1;1,1)$ in the theorem; again, recall
Example~\ref{ex:22}.  Also, setting $(a,b;\allowbreak c_0,c_\infty)$
equal to $(0,2;1,1)$ yields an additional known UBT pair with a
combinatorial interpretation, which relates the type\nobreakdash-$B$
Eulerian numbers $E_{n,k}(0,2;1,1)\eqdef{\eulerian{n}k}_B$ (see
\cite[$\texttt{A060187}$]{OEIS2022}) and the type\nobreakdash-$B$
Stirling subset numbers $S_{n,k}(0,2;1)\eqdef{\stirsub{n}k}_B$ (see
\cite{Bagno2019} and \cite[$\texttt{A039755}$]{OEIS2022}).

For each~$n$, sequences $(2)^{\overline
  k}S_{n,k}(0,1;1)=(k+1)!\,\stirsub{n+1}{k+1}$ and
$E_{n,k}(0,1;1,1)=\eulerian{n+1}k$, resp.\ $2^{k}k!\,S_{n,k}(0,2;1)$
and $E_{n,k}(0,2;1,1)={\eulerian{n}k}_B$, where $0\le k\le n$, arise
combinatorially as the \emph{f}\nobreakdash-vector and
\emph{h}\nobreakdash-vector of a simplicial complex dual to the
permutohedron of type~$A_n$, resp.~$B_n$~\cite{Fomin2007}.  In that
context, the composite transformation $\textrm{UBT}\circ\textrm{RT}$
implicit in~(\ref{eq:ubt1}) can be identified with the `reverse
Pascal's triangle' construction that maps the
\emph{f}\nobreakdash-vector of a dual simplicial complex to its
\emph{h}\nobreakdash-vector~\cite[Example~5.6]{Fomin2007}.

\subsection{Connection coefficient interpretation}
\label{subsec:32}

Now that the generalized Eulerian numbers $E_{n,k}(a,b;c_0,c_\infty)$
have been introduced as the elements of a parametric GKP triangle, how
they can be efficiently computed will be explained.  For general
parameter values, it is difficult to extract a useful closed-form
expression from~(\ref{eq:speck}), their bivariate~EGF\null.  They can
be computed alternatively by~(\ref{eq:ubt1}) from the Hsu--Shiue
generalized Stirling numbers, which in~turn are given by the
rank\nobreakdash-$1$ formula~(\ref{eq:expanded}).  But the resulting
formula for $E_{n,k}(a,b;c_0,c_\infty)$ is of rank~$2$: it involves a
double summation.

It may be possible to simplify this, but a formula without a multi-sum
can be worked~out by another technique, which is of independent
interest.  Besides yielding the following rank\nobreakdash-$1$
formula, the technique provides a Worpitzky-like interpretation of the
$E_{n,k}(a,b;c_0,c_\infty)$ as connection coefficients.

\begin{theorem}
\label{thm:addfactsE}
\hangindent\leftmargini         
  {\rm(i)} When\/ $b\neq0$, $E_{n,k}(a,b;c_0,c_\infty)$ is given by the
  rank\/\nobreakdash-$1$ formula
  \begin{equation}
    \label{eq:certainpolybasesE}
    \begin{split}
      & E_{n,k}(a,b;\,c_0,c_\infty) = 
      \\
      & \quad\frac{1}{b^k k!}
      \sum_{j=0}^k (-1)^{k-j} \binom{k}j
      (bn+c_0+c_\infty)^{\underline{k-j},b} (c_0 + c_\infty)^{\overline{j}, b} (bj+c_0)^{\underline{n},a}.
    \end{split}
  \end{equation}
  \begin{itemize}
  \item[{\rm(ii)}] In general the numbers\/ $E_{n,k}(a,b;c_0,c_\infty)$, $0\le k\le n$, satisfy
  \begin{equation}
    (c_0+c_\infty)^{\overline{n},b} (x)^{\underline{n},a}
    =
    \sum_{k=0}^n E_{n,k}(a,b;\,c_0,c_\infty)\, (x-c_0)^{\underline{k},b} (x+c_\infty)^{\overline{n-k},b},
  \label{eq:conncoeffinterp}
  \end{equation}
  which if\/ $(c_0+c_\infty)^{\overline{n},b}\neq0$, defines them as
  coefficients that express the factorial polynomial\/
  $(x)^{\underline{n},a}$ with respect to a bifactorial basis of the\/
  $(n+\nobreak1)$-dimensional space of polynomials of degree~${\le
    n}$.
  \end{itemize}
\end{theorem}
\begin{proof}
  Let $G(t,z)$ denote the EGF of the GKP triangle defining $E_{n,k}(a,b;c_0,c_\infty)$, and let $G_n(t)$, $n\ge0$, be its row polynomials.  That is,
  \begin{equation}
    G(t,z)=
    \left[
      \begin{array}{ll|l}
        -a, & b & c_0 \\
        a+b, & -b & c_\infty \\
      \end{array}
      \right](t,z),
    \quad
    G_n(t) = \left[
      \begin{array}{ll|l}
        -a, & b & c_0 \\
        a+b, & -b & c_\infty \\
      \end{array}
      \right]_{n}(t),
  \end{equation}
  and the PDE satisfied by $G(t,z)$ and the differential recurrence
  satisfied by $G_n(t)$ are given in Theorem~\ref{thm:pde}.  The
  lifting transformations previously considered were of the form
  $G^*(t^*,z^*)=G(t,z)$ where $(t,z)=(R(t^*),S(t^*)z^*)$.  Consider
  the more general transformation
  \begin{equation}
    G^*(t^*,z^*) = Q(t^*)\,G\!\left(R(t^*),S(t^*)z^*\right)
  \end{equation}
specified by $(Q,R,S)$, where $Q$~satisfies $Q(0)=1$.  It may be
possible to choose $(Q,R,S)$ so that $G^*(t^*,z^*)$ is the EGF of an
infinite array ${\genvert{n}{k}}^*$ that satisfies a triangular
recurrence of GKP type, though the array will not be lower-triangular
if~$Q\not\equiv1$.

Consider in particular the case when $(R(t^*),S(t^*) =
\left(t^*,(1-t^*)^{-1}\right)$ and the prefactor $Q(t^*)$ equals
$(1-\nobreak t^*)^{-(c_0+c_\infty)/b}$, so that $t=t^*$ and
\begin{equation}
\label{eq:wprefactor}
  G_n^*(t) = (1-t)^{-n-(c_0+c_\infty)/b} G_n(t).
\end{equation}
The unlifted row polynomials $G_n(t)$, $n\ge0$ satisfy the recurrence
\begin{equation}
  G_{n+1} = \left\{
  \left[-a+(a+b)t\right]n
  +(c_0+c_\infty t)
  \right\} G_n + b(1-t)t\,G_n',
\end{equation}
with $G_0(t)\equiv1$.  By direct computation the lifted $G_n^*(t)$,
$n\ge0$, which by~(\ref{eq:wprefactor}) are not polynomials, satisfy
the recurrence
\begin{equation}
\label{eq:preclift}
  G^*_{n+1} = (
  -a n
  +c_0
  ) G^*_n + bt(G_n^*)',
\end{equation}
with $G^*_0(t) =\allowbreak (1-\nobreak t)^{-(c_0+c_\infty)/b}$.
Substituting $G_n^*(t) = \sum_{k=0}^\infty {\genvert{n}{k}}^*t^k$
into~(\ref{eq:preclift}) reveals that the lifted array coefficients
${\genvert{n}{k}}^*$ satisfy a degenerate recurrence of the GKP type,
namely
\begin{equation}
  {\genvert{n+1}{k}}^* = (-an + bk + c_0) {\genvert{n}{k}}^*,\qquad \textrm{$(n,k)\in\mathbb{N}^2$},
\end{equation}
with the non-GKP initial condition
\begin{equation}
  {\genvert0k}^* =
  \binom{(c_0+c_\infty)/b -1 + k}{k}, \qquad k\in\mathbb{N}.
\end{equation}
The explicit formula
\begin{equation}
\label{eq:explicitbyinspec}
    {\genvert{n}k}^* =
  \binom{(c_0+c_\infty)/b -1 + k}{k} (bk+c_0)^{\underline{n},a}, \quad \textrm{$(n,k)\in\mathbb{N}^2$},
\end{equation}
follows by inspection.

Now consider (\ref{eq:wprefactor}) above.  Expanding the prefactor and
likewise its reciprocal in geometric series, and equating like powers
of~$t$, yields the pair
\begin{subequations}
\begin{align}
\label{eq:newpaira}
  {\genvert{n}{k}}^* &= \sum_{j=0}^k
  \binom{n+(c_0+c_\infty)/b -1 + k-j}{k-j} \genvert{n}j,\\
  {\genvert{n}{k}} &= \sum_{j=0}^k (-1)^{k-j}
  \binom{n+(c_0+c_\infty)/b}{k-j} {\genvert{n}{j}}^*,
\label{eq:newpairb}
\end{align}
\end{subequations}
holding for all $(n,k)\in\mathbb{N}^2$.  Here, $\genvert{n}k$
signifies $E_{n,k}(a,b;c_0,c_\infty)$.  Taking
(\ref{eq:explicitbyinspec}) into account, one sees that
(\ref{eq:newpairb}) is equivalent to the claimed
formula~(\ref{eq:certainpolybasesE}) for $E_{n,k}(a,b;c_0,c_\infty)$.

Continuing, let $\delta\defeq (c_0+c_\infty)/b-1$, so that
(\ref{eq:newpaira}) says that
\begin{equation}
  \binom{\delta+k}k (bk+c_0)^{\underline{n},a}
  = \sum_{j=0}^k \binom{n+\delta+k-j}{k-j} E_{n,j}(a,b;\,c_0,c_\infty).
\end{equation}
If $k\le n$ the summation $\sum_{j=0}^k$ can obviously be replaced by
$\sum_{j=0}^n$, and multiplying both sides by $k!\,(\delta+\nobreak
n)^{\underline{n-k}}$ then yields
\begin{equation}
  (\delta+1)^{\overline{n}} (bk+c_0)^{\underline{n},a} =
  \sum_{j=0}^n (k)^{\underline{j}} (k+\delta +1)^{\overline{n-j}} 
  E_{n,j}(a,b;c_0,c_\infty).
\end{equation}
By the formal substitution $k=(x-c_0)/b$ this becomes a statement
that the equality
\begin{equation}
  (\delta+1)^{\overline{n}} (x)^{\underline{n},a} =
  \sum_{j=0}^n \left((x-c_0)/b\right)^{\underline{j}} ((x-c_0)/b+\delta +1)^{\overline{n-j}} 
  E_{n,j}(a,b;c_0,c_\infty)
\end{equation}
holds when $x$~takes on any of the $n+1$ distinct values $bk+c_0$,
$k=0,1,\dots,n$.  But both sides are degree\nobreakdash-$n$
polynomials in~$x$, so the equality must hold for all~$x$.  This is
equivalent to the claim~(\ref{eq:conncoeffinterp}), the $b=0$ case
following by taking the $b\to0$ limit.
\end{proof}

An immediate corollary of the connection
formula~(\ref{eq:conncoeffinterp}) is the involutive identity
\begin{equation}
\label{eq:needinstance}
E_{n,n-k}(a,b;\,c_0,c_\infty) = E_{n,k}(-a,b;\,c_\infty,c_0).  
\end{equation}
This also follows from the RT formula~(\ref{eq:firstc}), if one uses
the definition (\ref{eq:euleriannew1}) of $E_{n,k}$ as a GKP triangle.
An instance of~(\ref{eq:needinstance}), coming from the choice
$(a,b;c_0,c_\infty)=(0,1;1,1)$, or equivalently
$(\alpha,\beta,\gamma;\allowbreak
\alpha',\beta',\gamma')=(0,1,1;\allowbreak 1,-1,1)$ as
in~(\ref{eq:threec}), is the classical Eulerian identity
$\eulerian{n+1}{n-k} = \eulerian{n+1}k$, $0\le k\le n$.  In the
traditional indexing this is $A_{n+1,n-k+1} = A_{n+1, k+1}$, $0\le
k\le n$.

Theorem~\ref{thm:addfactsE}(i) and its inverse, which comes from
(\ref{eq:newpaira}) rather than (\ref{eq:newpairb}), can be rephrased
in the following way.

\begin{theorem}
\label{thm:rephrased}
  For all\/ $n\ge0$ and choices of parameters\/ $(a,b;c_0,c_\infty)$, one
  has the binomial transform pair
\begin{subequations}
\begin{sizealign}{\small}
  &b^kk!\,E_{n,k}(a,b;\,c_0,c_\infty)
  = \sum_{j=0}^k (-1)^{k-j} \binom{k}{j} (bn+c_0+c_\infty)^{\underline{k-j},b} (c_0+c_\infty)^{\overline{j},b}\, (bj+c_0)^{\underline{n},a},\label{eq:rephraseda}\\
  &(c_0+c_\infty)^{\overline{k},b} (bk+c_0)^{\underline{n},a}
  = \sum_{j=0}^k \binom{k}j (bn+c_0+c_\infty)^{\overline{k-j},b}\, b^jj! \,E_{n,j}(a,b;\,c_0,c_\infty).\label{eq:rephrasedb}
\end{sizealign}
\end{subequations}
\end{theorem}

When $(a,b;c_0,c_\infty)=(0,1;1,0)$, the formula (\ref{eq:rephraseda})
reduces to the classical formula~(\ref{eq:eulerian5})
for~$\eulerian{n}k$, which is well known~\cite{Comtet74,Graham94}.
Slightly modified versions of this formula can be obtained by setting
$(a,b;c_0,c_\infty)=(0,1;0,1)$ and $(0,1;1,1)$.  It must be mentioned
that He and Shiue~\cite[eqs.\ (15),(17)]{He2020} recently obtained the
$(a,b;c_0,c_\infty)=(\theta,1;1,0)$ specialization not only
of~(\ref{eq:rephraseda}), but of the inverse identity
(\ref{eq:rephrasedb}) as~well.

\smallskip
Theorem~\ref{thm:rephrased} should be compared with
Theorem~\ref{thm:ubtpair}, which revealed that for all $n\ge0$, the
row sequences $E_{n,k}(a,b;c_0,c_\infty)$, $k=0,\dots,n$ and
$(c_0+\nobreak c_\infty)^{\overline{n-k},b}\allowbreak
S_{n,n-k}(-a,b;c_\infty)$, $k=0,\dots,n$, form a UBT pair; and also
with the following.

\begin{theorem}
\label{thm:lbtpair}
  For all\/ $n\ge0$ and choices of parameters\/ $(a,b;r)$, one
  has the LBT\/ {\rm(}lower binomial transform\/{\rm)} pair
\begin{subequations}
\begin{align}
  b^kk!\,S_{n,k}(a,b;\,r)
  &= \sum_{j=0}^k (-1)^{k-j} \binom{k}{j} (bj+r)^{\underline{n},a},\\
  (bk+r)^{\underline{n},a}
  &= \sum_{j=0}^k \binom{k}j  b^jj! \,S_{n,j}(a,b;\,r).
\end{align}
\end{subequations}
\end{theorem}
\begin{proof}
  The first identity is from Theorem\/~{\rm\ref{thm:addfacts}}, and
  the LBT is inverted in the second.  For the theory of LBT's, see
  Boyadzhiev~\cite{Boyadzhiev2018}.
\end{proof}

For any $n\ge0$, the UBT pair in Theorem~\ref{thm:ubtpair} and the LBT
pair in Theorem~\ref{thm:lbtpair} are of the classical forms
\begin{align}
\label{eq:438}
  &
  v_k = \sum_{j=k}^n (-1)^{j-k} \binom{j}k u_j
  \;\Longleftrightarrow\;
  u_k = \sum_{j=k}^n  \binom{j}k v_j\\
\shortintertext{and}
  &
  v_k = \sum_{j=0}^k (-1)^{k-j} \binom{k}j u_j
  \;\Longleftrightarrow\;
  u_k = \sum_{j=0}^k  \binom{k}j v_j,
\label{eq:lessgen}
\end{align}  
respectively.  In each, two sequences $(u_k)_{k=0}^n$, $(v_k)_{k=0}^n$
are related by a non-involutive binomial transform, or equivalently by
its inverse.  But the LBT pair in Theorem~\ref{thm:rephrased} is of a
more general form than~(\ref{eq:lessgen}), namely
\begin{equation}
\label{eq:moregen}
  v_k = \sum_{j=0}^k (-1)^{k-j} \binom{k}j (A)^{\underline{k-j},b}\,u_j
  \,\Longleftrightarrow\,
  u_k = \sum_{j=0}^k  \binom{k}j (A)^{\overline{k-j},b}\,v_j,
\end{equation}
where $A\neq0$ and $b$ are constants.  (Actually, in the theorem
$A$~depends affine-linearly on the row index~$n$, but because these
transformations act row-wise, that is not a major matter.)  When $A=1$
and $b=0$, (\ref{eq:moregen}) reduces to~(\ref{eq:lessgen}).  Unusual
LBT pairs of the type~(\ref{eq:moregen}), which have not appeared
widely in the literature, are briefly discussed at the end of the next
subsection.

\subsection{Generalized Stirling formulas}
\label{subsec:33}

The generalized Eulerian numbers $E_{n,k}(a,b;c_0,c_\infty)$ can be
computed by Theorem~\ref{thm:ubtpair} from the Hsu--Shiue Stirling
numbers $S_{n,k}$, as well as by the rank\nobreakdash-$1$ summation
formula in Theorem~\ref{thm:addfactsE} (when $b\neq0$).  Because of
the existence of the former method, some identities satisfied by the
$S_{n,k}$ in general, and explicit formulas for certain~$S_{n,k}$,
will now be given.  The formulas will incidentally lead to some
sequence transformations that generalize the classical Stirling
transform (for which see~\cite{Boyadzhiev2018}).  By exploiting the
previously mentioned homogeneity property
\begin{gather}
S_{n,k}(\lambda a,\lambda b;\, \lambda r) =  \lambda ^{n-k}S_{n,k}(a,b;\,r),
\end{gather}
additional explicit formulas can be generated.

Identities (i) and~(ii) in the following theorem could be called
`contiguous function relations' in the spirit of Gauss.
\begin{theorem}
  \label{thm:2contigsK}
  For all\/ $n\ge0$,
  the numbers\/ $S_{n,k}$ satisfy
  \begin{enumerate}
    \item[{\rm (i)}] $S_{n,k}(a,b;r+a) = S_{n,k}(a,b;r) + an\,S_{n-1,k}(a,b;r)$,
    \item[{\rm (ii)}] $S_{n,k}(a,b;r+b) = S_{n,k}(a,b;r) + b(k+1)\,S_{n,k+1}(a,b;r)$,
    \item[{\rm (iii)}] $S_{n,k}(-a,b;r) = S_{n,k}\left(a,b;r+a(n-1)\right)$,
    \item[{\rm (iv)}] $S_{n,k}(a,-b;r) = S_{n,k}\left(a,b;r-bk\right)$,
    \item[{\rm (v)}] $S_{n,k}(a,b;b-a) = S_{n+1,k+1}\left(a,b;0\right)$,
  \end{enumerate}
  when\/ $0\le k\le n$, with\/ $k=-1$ also allowed in\/ {\rm(ii)}
  and\/~{\rm(v)}.  In\/ {\rm(i)}, {\rm(ii)}, and\/~{\rm(v)}, the
  convention that\/ $S_{n,k}=0$ if\/ $k<0$, $k>n$, or\/ $n<0$, is
  adhered to.
\end{theorem}
\begin{proof}
  Each of (i),(ii),(iii) follows by elementary manipulation of the
  bivariate EGF~(\ref{eq:elem}), the cases $a=0$ and $b=0$ holding by
  continuity; though (iv) follows more easily from the
  finite-difference representation~(\ref{eq:preexpanded}).
  Identity~(v) is an example of the left-trimming of a GKP triangle,
  as in Theorem~\ref{thm:trimming}(i).
\end{proof}

\begin{example}
The $r$\nobreakdash-Stirling subset numbers ${\stirsub{n}{k}}_r =
S_{n,k}(0,1;r)$ and cycle numbers ${\stircyc{n}{k}}_r =
S_{n,k}(-1,0;r)$ (where usually $r\in\mathbb{N}$,
see~\cite{Broder84}), generalize $\stirsub{n}{k}$
and~$\stircyc{n}{k}$.  For an $(n+\nobreak r)$\nobreakdash-set, they
count restricted partitions with $k+\nobreak r$~blocks,
resp.\ permutations with $k+\nobreak r$ cycles, the restriction being
that $r$~distinguished elements of the set must be placed in distinct
blocks, resp.\ cycles.\footnote{In present notation, the
  $r$\nobreakdash-Stirling numbers of~\cite{Broder84} would be written
  as ${\stirsub{n-r}{k-r}}_r$ and ${\stircyc{n-r}{k-r}}_r$.}  The
$r$\nobreakdash-Stirling subset numbers can be computed by the
rank\nobreakdash-$1$ formula of Theorem~\ref{thm:addfacts}.
Additionally, setting $(a,b)=(0,1)$ in part~(ii) and $(-1,0)$ in
part~(i) of Theorem~\ref{thm:2contigsK} yields respectively that when
$0\le k\le n$,
\begin{subequations}
\begin{gather}
{\stirsub{n}{k}}_r = {\stirsub{n}{k}}_{r+1}\! - (k+1) {\stirsub{n}{k+1}}_r , \\
\shortintertext{and}
{\stircyc{n+1}{k}}_r = {\stircyc{n+1}{k}}_{r+1}\! - (n+1) {\stircyc{n}{k}}_{r+1}.
\end{gather}
\end{subequations}
These `cross' recurrences are equivalent to known ones with
combinatorial interpretations~\cite[\S3]{Broder84}.
\end{example}

The following theorem says that in a sense, the parametric family of
Hsu--Shiue Stirling number triangles is closed under the taking of
(row-wise) upper binomial transforms.

\begin{theorem}
  For all\/ $0\le k\le n$ and\/ $\delta$, the\/ $S_{n,k}$ satisfy
  \begin{equation*}
    \delta^{\underline k,b} S_{n,k}(a,-b;\,r+\delta) =
    \sum_{j=k}^n \binom{j}k\, \delta^{\underline j,b}S_{n,j}(a,b;\,r).
  \end{equation*}
\end{theorem}
\begin{proof}
  This follows from the UBT formula\/~{\rm(\ref{eq:secondc})}, if one
  views $\delta^{\underline k,b} S_{n,k}(a,b;r)$ as the GKP triangle\/
  $\left[\begin{smallarray}{cc|c}-a&b\,&r\\0&-b\,&\delta\end{smallarray}\right]_{n,k}$.
  (See\/~{\rm(\ref{eq:infact})}.)
\end{proof}

To prove the next theorem, recall from (\ref{eq:vertical}) that the
matrix $\mathcal{S}(a,b;r) = \left(S_{n,k}(a,b;r)\right)$ is an
exponential Riordan array:
\begin{equation}
\label{eq:shownas}
  \mathcal{S}(a,b;\,r) =
  \left[
    (1+az)^{r/a}, \frac{(1+az)^{b/a}-1}b
  \right],
\end{equation}
meaning that
\begin{equation}
\label{eq:meantthat}
  S_{n,k}(a,b;\,r) = \frac{n!}{k!}[z^n]
(1+az)^{r/a} \left[\frac{(1+az)^{b/a}-1}b\right]^k.
\end{equation}
(Taking the $a\to0$ limit, if desired, is straightforward.)  It is a
fundamental fact~\cite{Barry2016} that (exponential) Riordan arrays
form a group under matrix multiplication: if $d_i,h_i$ are formal
power series in~$z$ of respective orders~$0,1$, for $i=1,2$, then
$[d_1,h_1]\,[d_2,h_2] = \allowbreak [(d_2\circ h_1)d_1, \allowbreak
  h_2\circ h_1]$ and $[d_1,h_1]^{-1} =\allowbreak [1/(d_1\circ \bar
  h_1), \bar h_1]$, where $\bar h_1$~is the compositional inverse
of~$h_1$.

Of the following three identities, (i)~is presumably well known and
(ii)~appears in~\cite{Hsu98}; (iii)~extends a pair of identities of
Can and De\u{g}l\i\ \cite[eqs.\ (30)--(31)]{Can2014}, and appears to
be new.

\begin{theorem}
\label{thm:3Sidents}
\hangindent\leftmargini         
  {\rm(i)}~The parametric, infinite lower-triangular matrix\/
  $\mathcal{S}(a,b;r)$, i.e., $\left(S_{n,k}(a,b;r)\right)$, satisfies
  the product formula
  \begin{displaymath}
    \mathcal{S}(a,c;\,r_1+r_2) =     \mathcal{S}(a,b;\,r_1) \mathcal{S}(b,c;\,r_2);
  \end{displaymath}
and\/ $\mathcal{S}(a,b;r)^{-1} = \mathcal{S}(b,a;-r)$, as\/
$\mathcal{S}(a,a;0)=\mathcal{I}$, the identity matrix, for any\/~$a$.
\begin{itemize}
\item[{\rm(ii)}] For any non-negative\/ $n$ and\/ $k,k_1,k_2$ satisfying\/ $k=k_1+k_2$,
one has the convolution formula
\begin{displaymath}
  \frac{k!}{k_1!\,k_2!} S_{n,k}(a,b;r_1+r_2) = \sum \frac{n!}{n_1!\,n_2!} S_{n_1,k_1}(a,b;\,r_1)S_{n_2,k_2}(a,b;r_2),
\end{displaymath}
the sum being over non-negative pairs\/ $n_1,n_2$ satisfying\/ $n=n_1+n_2$, with
$n_1\ge k_1$ and $n_2\ge k_2$.
\item[{\rm(iii)}] For any non-negative\/ $k$ and\/ $n_1,n,k_2$ satisfying\/ $n_1=n+k_2$,
one has the asymmetric convolution formula
\begin{displaymath}
  \frac{n_1!}{n!\,k_2!} S_{n,k}(a,b;r_1+r_2) = \sum \frac{k_1!}{k!\,n_2!} S_{n_1,k_1}(a,b;\,r_1)S_{n_2,k_2}(b,a;\,r_2),
\end{displaymath}
the sum being over non-negative pairs\/ $k_1,n_2$ satisfying\/
$k_1=k+n_2$, with\/ $n_1\ge k_1$ and\/ $n_2\ge k_2$.  {\rm(}Note the
interchange of\/ $a,b$ in the summand.{\rm)}
\end{itemize}
\end{theorem}
\begin{proof}
  (i)~These facts are immediate corollaries of the connection
  formula~(\ref{eq:certainpolybases}).  They also have a Riordan-array
  interpretation: as is easily verified, they come from the
  just-stated formulas that express $[d_1,h_1]\,[d_2,h_2]$ and
  $[d_1,h_1]^{-1}$ as Riordan arrays, when $d_i,h_i$ depend on~$z$ as
  shown in~(\ref{eq:shownas}).

  (ii)~To prove this, substitute $(n,k)=(n_i,k_i)$ in the vertical EGF
  formula (\ref{eq:vertical}), and take the product of two copies
  of~it: one with $i=1$ and one with~$i=2$.  Then, equate the
  coefficients of like powers of~$z$ on the left and right sides.

  (iii)~Consider the infinite matrix ${B}(a,b;r) =
  \left(B_{n,w}(a,b;r)\right)$ defined by
  \begin{equation}
    \label{eq:Bdef}
  B_{n,w}(a,b;\,r) = {n!}[z^n]
(1+az)^{r/a} \left[\frac{bz}{(1+az)^{b/a}-1}\right]^w.
  \end{equation}
  This is an example of an improper Riordan array~\cite{Balas2005}.
  It is not a lower-triangular matrix indexed by $n,k\ge0$.  Rather,
  it is indexed by $n\ge0$, $w\in\mathbb{Z}$.  As a power series
  in~$z$, the quantity raised to the $w$'th power here is of
  order~$0$, not~$1$.

  If $r=r_1+r_2$, for any $k\ge0$ one has by elementary manipulations
  \begin{equation}
  \begin{aligned}
    &(1+az)^{r/a}\left[
      \frac{bz}{(1+az)^{b/a}-1}
      \right]^w \frac1{k!} \left[
      \frac{(1+az)^{b/a}-1}{b}
      \right]^k\\
    &\qquad {}=
    \sum_{n=k}^\infty
    \left[
      \sum_{j=k}^n \frac{1}{j!(n-j)!} B_{n-j,w}(a,b;\,r_1) S_{j,k}(a,b;\,r_2)
      \right]z^n.
  \end{aligned}
  \end{equation}
  If moreover $k\ge w$, this alternatively equals
  \begin{equation}
  \begin{aligned}
    \frac{z^w}{k!} (1+az)^{r/a}
    \left[
      \frac{(1+az)^{b/a}-1}{b}
      \right]^{k-w}\!\! &= \frac{z^w(k-w)!}{k!} \sum_{n=k-w}^\infty S_{n,k-w}(a,b;\,r)\frac{z^n}{n!}\\
    &=
    \sum_{n=k}^\infty\left[
      S_{n-w,k-w}(a,b;\,r)
      \frac{(k-w)!}{k!(n-w)!}
      \right] z^n.
  \end{aligned}
  \end{equation}
Equating the coefficients of like powers of~$z$ yields the identity
\begin{equation}
\label{eq:thisident}
  \frac{\binom{n}{k}}{\binom{n-w}{k-w}} S_{n-w,k-w}(a,b;\,r)
  =
  \sum_{j=k}^n \binom{n}{j}
  B_{n-j,w}(a,b;\,r_1) S_{j,k}(a,b;\,r_2),
\end{equation}
which holds when $n,k\ge\max(w,0)$.

In this, only the $w$'th column of ${B} = \left(B_{n,w}\right)$
appears.  Define a parametric, lower-diagonal Toeplitz matrix
$\mathcal{B}^{(w)}(a,b;r_1) =
\left(\mathcal{B}_{n,j}^{(w)}(a,b;r_1)\right)$, indexed by $n,j\ge0$,
which depends only on this column vector:
\begin{equation}
\mathcal{B}_{n,j}^{(w)}(a,b;r_1) = 
\begin{cases}
\binom{n}{j} B_{n-j,w}(a,b;\,r_1), & n\ge j,\\
0, & \textrm{otherwise}.
\end{cases}
\end{equation}
The identity (\ref{eq:thisident}) can be rewritten as
\begin{equation}
\label{eq:thatident}
  \frac{\binom{n}{k}}{\binom{n-w}{k-w}} S_{n-w,k-w}(a,b;\,r)
  =
  \sum_{j=k}^n
  \mathcal{B}^{(w)}_{n,j}(a,b;\,r_1) S_{j,k}(a,b;\,r_2),
\end{equation}
in which each side is an element of a matrix indexed by
$n,k\ge\max(w,0)$, and the right-hand side computes the product of two
such matrices.

It is known that if $g(z)=c_0+g_1z+\dots$ is a formal power series and
$M=M(g)$ is the lower-triangular Toeplitz matrix defined by
$M_{n,j}=g_{n-j}$, the map $g\mapsto M(g)$ is an algebra isomorphism.
It follows by examining (\ref{eq:Bdef}) that the inverse of
$\mathcal{B}^{(w)}(a,b;r_1)$ is $\mathcal{B}^{(-w)}(a,b;-r_1)$,
because negating $w,r$ in~(\ref{eq:Bdef}) replaces the power series
in~$z$ to which $[z^n]$~is applied by its reciprocal.

Taking this into account and computing the matrix inverse of both
sides of~(\ref{eq:thatident}) yields the inverted identity
\begin{equation}
  \frac{\binom{n}{k}}{\binom{n-w}{k-w}} S_{n-w,k-w}(b,a;\,-r)
  =
  \sum_{j=k}^n
  S_{n,j}(b,a;\,-r_2) \mathcal{B}^{(-w)}_{j,k}(a,b;\,-r_1),
\end{equation}
which again holds when $n,k\ge\max(w,0)$.  The case when $w\ge0$ is of
interest here.  Comparing (\ref{eq:meantthat}) and (\ref{eq:Bdef})
reveals that if $w\ge0$,
\begin{equation}
  B_{n,-w}(a,b;\,r) = {\binom{n+w}{w}}^{-1} S_{n+w,w}(a,b;\,r),
\end{equation}
so that the inverted identity can be rewritten as
\begin{equation}
  \binom{n}w S_{n-w,k-w}(b,a;\,-r) = \sum_{j=k}^n
  S_{n,j}(b,a;\,-r_2) \binom{j}{k-w}S_{j-k+w,w}(a,b;\,-r_1).
\end{equation}
By negating $r,r_1,r_2$ and interchanging $a,b$, one sees that this is
equivalent to the claimed asymmetric convolution formula.
\end{proof}

The following theorem lists formulas for certain $S_{n,k}(a,b;r)$
Hsu--Shiue triangles, parametrized by~$r$.  They may be known but seem
not to have not been assembled before in a single place.  The notation
used for a hypergeometric term in parts (iv),(v), and in the sequel,
is similar to that often used in the series for~${}_2F_1$:
\begin{equation}
  {}_2F_1\left(
  \myatop{A,\,B}{C}
  \biggm|
  w
  \right)
  =
  \sum_{k=0}^\infty
    \frac{A^{\overline k}\,B^{\overline k}}{1^{\overline k}\,C^{\overline k}}\,
    w^k
    \eqdef
  \sum_{k=0}^\infty
  \left[
    \myatop{A,\,B}{1,\,C}
    \right]^{\overline k} \!w^k.  
\end{equation}
That is,
\begin{equation}
  \left[
    \myatop{A_1,\dots,A_p}{C_1,\dots,C_q}
    \right]^{\overline k}
  \defeq
  \frac{(A_1)^{\overline k} \dotsm (A_p)^{\overline k}}
       {(C_1)^{\overline k} \dotsm (C_q)^{\overline k}}\,.
\end{equation}
By exploiting homogeneity and the identities of
Theorem~\ref{thm:3Sidents}, one can derive formulas for additional
triangles $S_{n,k}(a,b;r)$.

\begin{theorem}
  \label{thm:5points}
  For all\/ $n,k$ with\/ $n\ge k\ge0$, and for all\/ $r$,
  \begin{enumerate}
  \item[{\rm(i)}] $S_{n,k}(0,0;r) = \binom{n}{k} r^{n-k}$,
  \item[{\rm(ii)}] $S_{n,k}(1,1;r) = \binom{n}{k} r^{\underline{n-k}}$,
  \item[{\rm(iii)}] $S_{n,k}(-1,1;r) = \binom{n}{k} (n+r-1)^{\underline{n-k}}$,
  \item[{\rm(iv)}] $\begin{aligned} 
    S_{n,k}(1,2;r)
    &{}=\binom{n}{k}\frac{k!}{(2k-n)!} 2^{-(n-k)} {}_2F_1\left({\myatop{-r,\: -n+k}{-n+2k+1}}\biggm|2\right)\\
    &{}=2^{n-k}\left[ \myatop{-\frac{n}{2}, \: -\frac{n}2+\frac12}{1} \right]^{\overline {n-k}}  
    {}_2F_1\left({\myatop{-r,\: -n+k}{-n+2k+1}}\biggm|2\right),
  \end{aligned}$
  \item[{\rm(v)}] $\begin{aligned}
    S_{n,k}(-2,-1;r) 
    &{}=\binom{n}{k}\frac{(2n-k)!}{n!} 2^{-(n-k)} {}_2F_1\left({\myatop{r-1,\: -n+k}{-2n+k}}\biggm|2\right)\\
    &{}=\left(-\frac12\right)^{n-k}\left[ \myatop{-n, \: n+1}{1} \right]^{\overline {n-k}}  
    {}_2F_1\left({\myatop{r-1,\: -n+k}{-2n+k}}\biggm|2\right).
  \end{aligned}$
  \end{enumerate}
\end{theorem}
\begin{proof}
  A unified version of (i),(ii), namely $S_{n,k}(a,a;r) = \binom{n}{k}
  r^{\underline{n-k},a}$, appears as~\cite[eq.~(3.6)]{Maltenfort2020}.
  Also, (i),(ii),(iii) are valid because by examination, each
  satisfies the appropriate GKP recurrence, given in
  Definition~\ref{def:Sdef}.  Formulas equivalent to (iv),(v) were
  derived by Cheon, Jung, and
  Shapiro~\cite[eqs.~(11),(14)]{Cheon2013}, and each is restated here
  as a product of a (reversed) hypergeometric term and a terminating
  hypergeometric series.
\end{proof}

The ${}_2F_1(2)$ series in parts (iv) and~(v) must be interpreted with
care.  In~(iv), the lower hypergeometric parameter $-n+\nobreak
2k+\nobreak 1$ may be non-positive, causing a division by zero in the
terms of the series, but the division by $(2k-\nobreak n)!$
compensates for this.\footnote{The useful notation introduced by
  Olver~\cite[Chapter~15]{Olver2010} could be employed here.  If
  $C$~is the lower parameter of a ${}_2F_1$ function, ${}_2{\bf F}_1$
  signifies $\Gamma(C)^{-1}{}_2F_1$.  Unlike~${}_2F_1$, ${}_2{\bf
    F}_1$~is defined when $C$~is any non-positive integer, by taking a
  limit.}  It turns~out that $S_{n,k}(1,2;r)$ (when $r\in\mathbb{N}$)
is nonzero if and only if $0\le n-\nobreak k \le
\left\lfloor\frac{n+r}2\right\rfloor$.  When $(n,k)=(0,0)$, the
${}_2F_1(2)$ series in~(v) is troublesome also, and must be
interpreted as signifying unity.  It should be noted that the $r=0$
versions of the triangles $S_{n,k}(1,2;r)$ and $S_{n,k}(-2,-1;r)$ can
be left-trimmed, yielding the $r=1$ versions.

For all $r\in\mathbb{N}$, the triangle $S_{n,k}(1,1;r)$ in part~(ii)
has a combinatorial interpretation: its elements are rook numbers,
which count the number of ways of placing $n-\nobreak k$ non-attacking
rooks on an $r\times n$ chessboard.  Thus $S_{n,k}(1,1;r)$ is nonzero
only if $0\le n-\nobreak k \le r$.  The $r=0$ cases of the triangles
$S_{n,k}(a,b;r)$ in parts (iii),(iv),(v) are classical also, and can
be identified by examining the corresponding GKP recurrences.  The
elements $S_{n,k}(-1,1;0)$ are the (unsigned) Lah numbers $L_{n,k}$
(see \cite[$\texttt{A271703}$]{OEIS2022}), and the $S_{n,k}(1,2;0)$
and $S_{n,k}(-2,-1;0)$ are the second-kind Bessel numbers $B_{n,k}$
(see \cite[$\texttt{A122848}$]{OEIS2022}), resp.\ the (unsigned)
first-kind ones~$\hat b_{n,k}$ (see
\cite[$\texttt{A132062}$]{OEIS2022}).  The formulas for all three are
well known~\cite[p.~158]{Comtet74}.  They agree with the $r=0$ cases
of the formulas in (iii),(iv),(v).

In fact for all $r\in\mathbb{N}$, the generalized Stirling triangles
in (iii),(iv),(v) have combinatorial interpretations.  The elements
$S_{n,k}(-1,1;2r)\eqdef L_{n,k}^{(2r)}$ are the $r$\nobreakdash-Lah
numbers, which count restricted partitions of an $(n+\nobreak
r)$\nobreakdash-set into $k+\nobreak r$ lists, the restriction being
that $r$~distinguished elements of the set must be placed in distinct
lists~\cite{Nyul2015}.  The elements $S_{n,k}(1,2;r) \eqdef
B_{n,k}^{(r)}$ and $S_{n,k}(-2,-1;r) \eqdef \hat b_{n,k}^{(r)}$ are
the second-kind $r$\nobreakdash-Bessel numbers, resp.\ the (unsigned)
first-kind ones.  For an $(n+\nobreak r)$\nobreakdash-set,
$B_{n,k}^{(r)}$ counts partitions with each block having size
$1$~or~$2$, subject to the restriction that $r$~distinguished elements
must be placed in distinct blocks~\cite{Cheon2013}.  This
interpretation provides a proof that $S_{n,k}(1,2;r)$ (when
$r\in\mathbb{N}$) is nonzero if and only if $0\le n-\nobreak k \le
\left\lfloor\frac{n+r}2\right\rfloor$.

A notable feature of the hypergeometric series in parts (iv),(v) for
$S_{n,k}(1,2;r)$ and $S_{n,k}(-2,-1;r)$ is that when $r\in\mathbb{Z}$
with $r\ge0$, resp.\ $r\le1$, they terminate after $1+\nobreak
\min(r,n-\nobreak k)$ terms, resp.\ $1+\nobreak \min(1-r,n-\nobreak
k)$ terms.  That is, the number of terms does not grow with~$n$.  For
any such~$r$, this yields a `rank\nobreakdash-$0$' formula involving
no summation at~all, as in parts (i),(ii),(iii).

For example, $S_{n,k}(1,2;r)$ equals
\begin{equation}
\label{eq:subsumed}
  B_{n,k}^{(r)} = \binom{n}k \frac{k!}{(2k-n+r)!}  2^{-(n-k)}
  \times
  \begin{cases}
    1, & r=0,\\
    n+1, & r=1,\\
    n(n+1) + 2(k+1), &r=2,
  \end{cases}
\end{equation}
which holds when $2k-n+r\ge0$, i.e., $0\le n-\nobreak k \le
\left\lfloor\frac{n+r}2\right\rfloor$.  Similarly, $S_{n,k}(-2,-1;r)$
equals
\begin{equation}
  \hat b_{n,k}^{(r)} = \binom{n}k \frac{(2n-k)!}{n!(2n-k)^{\underline{1-r}}}2^{-(n-k)}
  \times
  \begin{cases}
    k(k+1)-2n, & r=-1,\\
    k, & r=0,\\
    1, &r=1,
  \end{cases}
\end{equation}
which holds when $0\le k\le n$, though if $2n-k+r\le0$ it must be
interpreted in a limiting sense: for instance, ${\hat b}_{0,0}^{(r)}$
always equals unity.  The curious identity
\begin{equation}
  \hat b_{n+1,k+1}^{(r)} = B^{(r)}_{2n-k, n},\qquad 0\le k\le n,
\end{equation}
holds not only when $r=0$, as has been noted~\cite{Han2008}, but also
when $r=1$.

By applying the coefficient extraction operator~$[z^n]$ to the
vertical univariate EGF~(\ref{eq:vertical}), one can derive the
alternative general formula
\begin{equation}
\label{eq:altgenform}
  \begin{aligned}
    B_{n,k}^{(r)} &=
    \frac{2^{-(n-k)}}{(n-k)!} {\binom{n+r}{r}}^{-1} (n+r)^{\underline{2n-2k}}
    \sum_{\ell=0}^r \binom{k+\ell}{\ell} \binom{n-k}{r-\ell} \\
    &{}=
    \frac{2^{-(n-k)}}{(n-k)!} {\binom{n+r}{r}}^{-1} (n+r)^{\underline{2n-2k}}
    \sum_{\ell=0}^r \binom{k+\lfloor\ell/2\rfloor}{\lfloor\ell/2\rfloor} \binom{n}{r-\ell},
  \end{aligned}
\end{equation}
which holds for all $r\in\mathbb{N}$ when $0\le k\le n$, and
subsumes~(\ref{eq:subsumed}).  Another approach to deriving formulas
for $B_{n,k}^{(r)}$ and~$\hat b_{n,k}^{(r)}$ would be to apply the
identities of Theorem~\ref{thm:2contigsK}, which allow the
parameter~$r$ in $S_{n,k}(a,b;r)$ to be incremented repeatedly by $a$
or by~$b$.
  
\smallskip
For any $r$ and $(a,b)$, such as the various choices in
Theorem~\ref{thm:5points}, the Hsu--Shiue array $S_{n,k}(a,b;r)$ can
be used to perform sequence transformations.  By
Theorem~\ref{thm:3Sidents}(i), the infinite lower-triangular matrix
$\mathcal{S}(a,b;r)$ has $\mathcal{S}(b,a;-r)$ as its inverse.  Hence
one has a `lower' transform pair
\begin{equation}
  v_k = \sum_{j=0}^k S_{k,j}(b,a;\,-r) \,u_j
    \,\Longleftrightarrow\,
  u_k = \sum_{j=0}^k S_{k,j}(a,b;\,r) \,v_j,
\end{equation}
and an `upper' transform pair
\begin{equation}
  v_k = \sum_{j=k}^\infty S_{j,k}(b,a;\,-r)\, u_j
    \,\Longleftrightarrow\,
  u_k = \sum_{j=k}^\infty S_{j,k}(a,b;\,r)\, v_j.
\end{equation}
One may need to require convergence in the latter, unless the
sequences terminate.

Such Hsu--Shiue--Stirling transforms include binomial transforms (for
which $(a,b;r)$ is $(0,0;1)$) and Stirling and
$r$\nobreakdash-Stirling transforms (for which it is $(0,1;0)$ or
$(0,1;r)$).  Other choices of $(a,b;r)$ seem not to have been much
investigated, such as the ones in parts (ii)--(v) of
Theorem~\ref{thm:5points}.  The choice $(a,b;r)=(1,1;r)$ in part~(ii)
is the most relevant here.  When $r\neq0$, it yields the pair of
`lower rook number transforms'
\begin{equation}
  v_k = \sum_{j=0}^k (-1)^{k-j} \binom{k}j r^{\underline{k-j},-1}\, u_j
    \,\Longleftrightarrow\,
  u_k = \sum_{j=0}^k \binom{k}j  r^{\overline{k-j},-1}\, v_j.
\end{equation}
which up to normalization is identical to the unusual
LBT~(\ref{eq:moregen}), an example of which appeared in
Theorem~\ref{thm:rephrased}.  The interplay between
Hsu--Shiue--Stirling transforms and the generalized Eulerian numbers
remains to be fully explored.

\subsection{Generalized Eulerian formulas}
\label{subsec:34}

In several cases, it is possible to derive a formula for the
generalized Eulerian numbers $E_{n,k}(a,b;c_0,c_\infty)$ not based on
a summation, or at least, not based on one in which the number of
terms grows with $n$ or~$k$.  Several such cases will now be explored.
By exploiting the previously mentioned reflection and homogeneity
properties
\begin{gather}
  \label{eq:reflectionlate}
  E_{n,k}(a,b;\,c_0,c_\infty) = E_{n,n-k}(-a,b;\, c_\infty,c_0),\\
  E_{n,k}(\lambda a,\lambda b;\, \lambda c_0, \lambda c_\infty) = \lambda^n E_{n,k}(a,b;\,c_0,c_\infty),
\end{gather}
additional ones can be generated.

\begin{theorem}
  \label{thm:simpleones}
  For all\/ $c_0,c_\infty$, one has
  \begin{enumerate}
  \item[{\rm(i)}] $E_{n,k}(0,0;\, c_0,c_\infty) =  \binom{n}k c_0^{{n-k}}c_{\infty}^{k} $,
  \item[{\rm(ii)}] $E_{n,k}(-1,1;\, c_0,c_\infty) =  \binom{n}k (c_0+k)^{\overline{n-k}}(c_{\infty})^{\underline k}$, and
  \item[{\rm(iii)}] $E_{n,k}(a,b;\,c_0,-c_0) = (-1)^k\binom{n}k (c_0)^{\underline{n},a}$, irrespective of\/~$b$.
  \end{enumerate}
\end{theorem}
\begin{proof}
  By examination, each of (i),(ii),(iii) satisfies the corresponding
  GKP recurrence, given in Definition~\ref{def:Edef}.  Note that
  formulas (i),(ii) can be unified: $E_{n,k}(-a,a;\, c_0,c_\infty)$
  equals $\binom{n}k
  (c_0+ka)^{\overline{n-k},a}(c_{\infty})^{\underline{k},a}$.  Another
  approach to (ii) and~(iii) is to note that when the parameter vector
  $(a,b;c_0,c_\infty)$ of~$E_{n,k}$ equals $(-1,1;c_0,c_\infty)$,
  resp.\ $(a,b;\allowbreak c_0,-c_0)$, the bivariate EGF
  (\ref{eq:speck}) reduces to
  \begin{displaymath}
    (1-z)^{-c_0-c_\infty} (1-z+tz)^{c_\infty},\qquad \textrm{resp.} \qquad (1+az-atz)^{c_0/a},
  \end{displaymath}
  from which (ii) and~(iii) follow by repeated differentiation.
  A~short approach to~(ii) is to note that $E_{n,k}(-1,1;\allowbreak
  c_0,c_\infty)$ and $(c_\infty)^{\underline {k}} \,S_{n,k}(-1,1;c_0)$ are
  identical, as they are both solutions of the GKP recurrence with
  parameter array
  $\left[\begin{smallarray}{cc|c}\alpha,&\beta&\gamma\\\alpha',&\beta'&\gamma'\end{smallarray}\right]
  =
  \left[\begin{smallarray}{cc|c}1,&1\,&c_0\\0,&-1\,&c_\infty\end{smallarray}\right]$.
  Then, one refers to Theorem~\ref{thm:5points}(iii).
\end{proof}

A more sophisticated case when the $E_{n,k}(a,b;c_0,c_\infty)$ may be
relatively easy to compute is the important `single progression' case
when $c_0+\nobreak c_\infty=b$, the first results on which were
apparently obtained by Carlitz~\cite[\S8]{Carlitz79}, whose
`degenerate' Eulerian numbers are of the innocuously normalized form
$E_{n,k}(\lambda,1;\allowbreak c_0,1-\nobreak c_0)$.  (When $c_0=1$,
these numbers have recently been placed in a combinatorial
setting~\cite{Herscovici2020}.)  Additional results are due to
Charalambides~\cite{Charalambides82}\footnote{In present notation, the
  composition number $A(m,k,r,s)$ in~\cite{Charalambides82} equals
  $E_{m,k}(1,s;\allowbreak r,s-\nobreak r)$ or
  $E_{m,m-k}(-1,s;\allowbreak,s-\nobreak r, r)$, divided by~$m!$.
  Some misprints in eqs.\ (2.26),(2.27) of~\cite{Charalambides82} are
  corrected in what follows.}, and Hsu and Shiue~\cite{Hsu99}.  As
mentioned in the introduction, the single-progression case includes
the usual Eulerian triangle $E_{n,k}(0,1;1,0)=\eulerian{n}k$, the
type\nobreakdash-$B$ Eulerian triangle
$E_{n,k}(0,2;1,1)={\eulerian{n}k}_B$, and others.

One
can write
\begin{equation}
  E_{n,k}(a,b;c_0,b-c_0) =
  \left[
    \begin{array}{cc|c}
      -a, & b & c_0\\
      a+b,& -b & b-c_0
    \end{array}
    \right]_{n,k},
\end{equation}
and if $b\neq0$,
\begin{equation}
  E_{n,k}(a,b;c_0,b-c_0) = b^n
  \left[
    \begin{array}{ccc}
      0, & 1, & \infty \\
      \hline
      a/b, & 1, & -a/b \\
      c_0/b, & \,-1\,, & 1-c_0/b
    \end{array}
    \right]_{n,k},
\end{equation}
in the new parametrization of Section~\ref{sec:characteristics}.
(Recall eqs.\ (\ref{eq:recalled1}),(\ref{eq:old2new}).)  One sees from
the parameter-pair in the second column of this tableau how the
single-progression case is special: not only is $r_1=1$, which is the
sign of the generalized Eulerian case (A\,III), but also~$g_1=-1$.

When $c_0+c_\infty=b$, the bivariate EGF (\ref{eq:speck}) reduces to
\begin{equation}
\label{eq:reducedspeck}
  \sum_{n=0}^\infty \sum_{k=0}^n E_{n,k}(a,b;\,c_0,b-c_0)
  \,t^k\,\frac{z^n}{n!} =
\frac{(1-t)(1+az-atz)^{c_0/a}}{1-t(1+az-atz)^{b/a}}
\end{equation}
(if $a\neq0$; taking the $a\to0$ limit is straightforward).  Moreover,
Theorems \ref{thm:ubtpair}(i) and~\ref{thm:addfactsE}(i) reduce to
\begin{equation}
    \label{eq:3waya}
    E_{n,k}(a,b;\, c_0, b-c_0) =
    \sum_{j=k}^n (-1)^{j-k} \binom{j}k b^{n-j} (n-j)!\,S_{n,n-j}(-a,b;\, b-c_0)
\end{equation}
and the rank\nobreakdash-$1$ formula
\begin{subequations}
  \begin{align}
    \label{eq:3wayb}
    E_{n,k}(a,b;\, c_0, b-c_0) 
    &=\sum_{j=0}^k (-1)^{k-j}\binom{n+1}{k-j} (bj+c_0)^{\underline{n},a}
    \\
    \label{eq:3wayc}
    &= \nabla_x^{n+1}
    \Bigl[
    \bigl( bx+c_0\bigr)^{\underline{n},a}\,
    {\bf1}_{0\le x\le k}
    \Bigr]  \Bigm|_{x=k},
  \end{align}
\end{subequations}
  where $\nabla_x$~is the backward first difference operator with
  respect to~$x$, defined by $\nabla_xf(x) =\allowbreak f(x)-\nobreak
  f(x-\nobreak 1)$, and $\mathbf{1}$ signifies an indicator function.
  The expression (\ref{eq:3wayc}), a~version of which was derived
  in~\cite{Charalambides82}, can be viewed as an unexpanded version of
  the preceding one.

  When $(a,b;c_0)=(0,1;1)$, resp.~$(0,1;0)$, (\ref{eq:3wayb})~and
  (\ref{eq:3wayc}) become known formulas for the Eulerian
  numbers~$\eulerian{n}k$, resp.\ the traditionally indexed numbers
  $\eulerian{n}{k-1}\defeq A_{n,k}$.  Similarly, when
  $(a,b;c_0)=(0,2;1)$, they become formulas for the
  type\nobreakdash-$B$ Eulerian numbers.  Each of the formulas
  (\ref{eq:3waya}),(\ref{eq:3wayb}),(\ref{eq:3wayc}) is useful, but
  certain restrictions on parameters may lead to alternative summation
  representations in which the number of terms does not grow with $n$
  or~$k$, as will now be seen.

  In the following theorem, parts (i) and~(ii) are analogues for the
  $E_{n,k}$ of the contiguity relations of
  Theorem~\ref{thm:2contigsK}, which applied to the~$S_{n,k}$.
  (Restricted versions of (i) and~(ii) were derived by
  Carlitz~\cite[\S8]{Carlitz79}.)  Note that though part~(ii) relates
  a pair of single-progression~$E_{n,k}$'s, (i)~holds more generally;
  and both can be iterated.  Parts (iii), (iv), and~(v) also relate
  parametric triangles $E_{n,k}$ which are not of the
  single-progression type.  They indicate how generalized Eulerian
  triangles can be left-trimmed (if~$c_0=0$), right-trimmed
  (if~$c_\infty=0$), and mid-trimmed (if $c_0+\nobreak c_\infty=0$).

  \begin{theorem}
    \label{thm:2contigsE}
    For all\/ $n\ge0$, the numbers\/ $E_{n,k}$ satisfy
    \begin{enumerate}
    \item[{\rm(i)}] $E_{n,k}(a,b; c_0+a, c_\infty-a) = E_{n,k} + an\left(E_{n-1,k}-E_{n-1,k-1}\right)$, the parameters of each number on the right-hand side being\/ $(a,b;c_0,c_\infty)$,  when\/ $0\le k\le n$;
    \item[{\rm(ii)}] $E_{n,k}(a,b;c_0+b,-c_0) = E_{n,k+1}(a,b;c_0,b-c_0) + (-1)^k (c_0)^{{\underline n},a}\binom{n+1}{k+1}$,\hfil\break when\/ $-1\le k\le n$;
    \item[{\rm(iii)}] $c_\infty\,E_{n,k}(a,b;b-a,c_\infty+a) = E_{n+1,k+1}(a,b;0,c_\infty)$,
      when $-1\le k\le n$;
    \item[{\rm(iv)}] $c_0\,E_{n,k}(a,b;c_0-a, a+b) = E_{n+1,k}(a,b;c_0,0)$, when\/
      $0\le k\le n$.
    \item[{\rm(v)}] $c\,[E_{n,k+1}-E_{n,k}](a,b;c-a, a-c) = E_{n+1,k+1}(a,b;c,-c)$, when\/ $-1\le k\le n$.
    \end{enumerate}
    In these, the convention that\/ $E_{n,k}=0$ if\/ $k<0$, $k>n$, or\/
    $n<0$, is adhered to.
  \end{theorem}
  \begin{proof}
    Identity~(i) follows readily from the bivariate
    EGF~(\ref{eq:speck}), and (ii) from its
    restriction~(\ref{eq:reducedspeck}) to the single-progression
    subcase, with the fact that the EGF of the term $(-1)^k
    (c_0)^{{\underline n},a}\binom{n+1}{k+1}$ equals
    $[(t-1)/t]\allowbreak (1+\nobreak az-\nobreak atz)^{c_0/a}$ taken
    into account.  Identities (iii) and~(iv) are specializations of
    Theorem~\ref{thm:trimming}, and (v)~follows from a fact indicated
    in Remark~\ref{rmk:mid}: the row polynomials $G_{n+1}(t)$,
    $n\ge0$, of any GKP triangle with $\gamma+\nobreak \gamma'=0$ and
    $\beta+\nobreak\beta'=0$ (with $\beta\beta'\neq0$) have
    $1-\nobreak t$ as a factor.  This factor can be divided~by,
    yielding (up~to a constant factor, here~$c$) the row polynomials
    $G_n^*(t)$ of a new, `mid-trimmed' GKP triangle.
  \end{proof}

  The following is an interesting formula for a triangle $E_{n,k}$ of
  single-progression type, having what amounts to a single discrete
  parameter: $1-\zeta+2p\in\mathbb{N}$.

  \begin{theorem}
    \label{thm:followsfrom}
    For all\/ $p\in\mathbb{N}$ and\/ $\zeta\in\{0,1\}$, 
    \begin{equation*}
      \begin{aligned}
        &E_{n,k}(-1,2;\,2-\zeta+2p, \zeta-2p) = \\
        &\qquad n!\binom{n+1}{2k+2p+1-\zeta}
        - (-1)^{k+p} \sum_{\ell=0}^{p-1} (-1)^{\ell} (2-\zeta+2\ell)^{\overline n}
        \binom{n+1}{k+p-\ell}.
      \end{aligned}
    \end{equation*}
  \end{theorem}
  \begin{proof}
    The $\zeta=0,1$ versions are proved independently, by induction
    on~$p$.  The inductive step uses Theorem~\ref{thm:2contigsE}(ii).
    The base ($p=0$) cases are respectively the triangles
    \begin{equation}
      E_{n,k}(-1,2;\,2-\zeta,\zeta) = n!\,\binom{n+1}{2k+1-\zeta},\qquad \zeta=0,1,
    \end{equation}
    which by examination satisfy the GKP recurrence of
    Definition~\ref{def:Edef}.  (These two number triangles,
    normalized by division by~$n!$, appear as $\texttt{A034839}$ and
    $\texttt{A034867}$ in the OEIS~\cite{OEIS2022}.)
  \end{proof}

  The number of terms in this summation formula does not grow with $n$
  or~$k$.  Besides applying Theorem~\ref{thm:2contigsE}(ii) to
  increment~$c_0$ and decrement~$c_\infty$ repeatedly (by~$b=2$), as
  was done in the proof, one could also apply
  Theorem~\ref{thm:2contigsE}(i) to decrement~$c_0$ and
  increment~$c_\infty$ repeatedly (by~$1$).  The number of terms in
  the resulting formula for any desired $E_{n,k}(-1,2;\allowbreak
  c_0,2-\nobreak c_0)$, $c_0\in\mathbb{Z}$, also will not grow with
  $n$ or~$k$.

  \smallskip
  Another notable feature of the single-progression case $c_0+\nobreak
  c_\infty=b$ is that in this case, the connection formula of
  Theorem~\ref{thm:addfactsE}(ii) reduces to the identity
  \begin{equation}
    \label{eq:consequential}
    b^n(x)^{\underline{n},a} = \frac1{n!} \sum_{k=0}^n
    E_{n,k}(a,b;\,c_0,b-c_0)
    \left[x-c_0+b(n-k)\right]^{\underline{n},b},
  \end{equation}
  or equivalently
  \begin{equation}
    \label{eq:rewrittenw}
    n!\,b^n (x)^{\underline n,a}=\sum_{k=0}^n E_{n,k}(a,b;\,b-c_\infty,c_\infty)
    (x+c_\infty - bk)^{\overline n,b}.
  \end{equation}
  This is a generalized Worpitzky identity, as the $a=0$ subcase makes
  clear.  When $(a,b;c_0,c_\infty)=(0,1;1,0)$, it reduces to the
  classical Worpitzky identity (\ref{eq:eulerian3}), and when
  $(a,b;c_0,c_\infty)=(0,1;0,1)$, to a slightly modified version.
  Also, when $(a,b;c_0,c_\infty)=(0,2;1,1)$, it reduces to the
  Worpitzky identity of type~$B$, which displays the
  type\nobreakdash-$B$ Eulerian numbers ${\eulerian{n}k}_B =
  E_{n,k}(0,2;1,1)$ as connection coefficients~\cite{Bagno2022}.

  The complementary subcase when $a\neq0$ has its own logic.  Setting
  $a=1$ in~(\ref{eq:consequential}) (which by homogeneity can be done
  without loss of generality), redefining the indeterminate~$x$, and
  applying the reflection property~(\ref{eq:reflectionlate}), produces
  the rather symmetric identity
\begin{equation}
\label{eq:applicable}
(b x+c_\infty)^{\underline n} = \sum_{k=0}^n \left[\frac1{n!}
  E_{n,k}(-1,b;\,b-c_\infty,c_\infty)\right]
(x+k)^{\underline n},
\end{equation}
which holds for arbitrary $b$ and~$c_\infty$.

For all $n\ge0$ and~$b$, let a matrix $\mathcal{A}^{(n,b)} =
(A^{(n,b)}_{k,j})$, $0\le k,j\le n$, which is not triangular, be
defined by
\begin{equation}
A_{k,j}^{(n,b)} = \frac1{n!} E_{n,j}(-1,b;\,b-k,k),
\end{equation}
Then (\ref{eq:applicable}), restricted to the case when
$c_\infty\in\{0,\dots n\}$, says that
\begin{equation}
\label{eq:followsfrom}
  (b x+k)^{\underline n} = \sum_{j=0}^n A_{k,j}^{(n,b)} (x+j)^{\underline n},
\qquad 0\le k\le n,
\end{equation}
which extends to
\begin{equation}
\label{eq:extendsto}
\Delta_x^r\left[
  (b x+k)^{\underline n}\right] = \sum_{j=0}^n A_{k,j}^{(n,b)}n^{\underline r}\, (x+j)^{\underline {n-r}},
\qquad 0\le k\le n,\quad r\in\mathbb{N}.
\end{equation}
Note that $\mathcal{A}^{(n,1)}$ equals $\mathcal{I}_{n+1}$, the
$(n+1)$-by-$(n+1)$ identity matrix, as follows from
Theorem~\ref{thm:simpleones}(ii).

Thus for all $n\ge0$ and $b\neq0$, the coefficients that connect the
two factorial bases $[(b x+\nobreak k)^{\underline n}]_{k=0}^n$ and
$[(x+\nobreak j)^{\underline n}]_{j=0}^n$ of the $(n+\nobreak
1)$\nobreakdash-dimensional space of polynomials of degree~${\le n}$
can be viewed as numbers in the $n$'th rows of certain
generalized Eulerian triangles of single-progression type, divided
by~$n!$.  Moreover, it follows from~(\ref{eq:followsfrom}) that for
all $n\ge0$, the map $b\mapsto \mathcal{A}^{(n,b)}$ is a homomorphism
from~$\mathbb{C}^*$, the multiplicative group of nonzero complex
numbers, to $GL(n+\nobreak 1,\mathbb{C})$.  This map and its image
deserve further study; empirically, one finds that the eigenvalues
of~$\mathcal{A}^{(n,b)}$ are $\{1,b,b^2,\dots,b^n\}$.

When $b$~is a positive integer, the preceding results make contact
with known identities, and $b=2$ is illustrative.  It follows from
Theorem~\ref{thm:followsfrom} that
\begin{equation}
  A_{k,j}^{(n,2)} = \binom{n+1}{2j-k+1}, \qquad 0\le k,j\le n.
\end{equation}
(It is assumed that $0\le 2j-k+1 \le n+1$, otherwise
$A_{k,j}^{(n,2)}$ vanishes.)  Hence (\ref{eq:followsfrom}) reduces
when $b=2$ to
\begin{equation}
  (2x+k)^{\underline n} =\sum_{j=0}^n
  \binom{n+1}{2j-k+1} (x+j)^{\underline n}, \qquad 0\le k\le n.
\end{equation}
(The summation includes only $j$ for which $0\le 2j-k+1 \le n+1$,
i.e., terms with $\left\lfloor \frac{k}2\right\rfloor \le j \le
\left\lfloor \frac{n+k}2\right\rfloor$.)  This identity exhibits the
coefficients of connection as binomial coefficients, and could be
proved directly.  When $b=3$ the coefficients of connection become
trinomial coefficients, and so~forth.

Because of the group homomorphism, for all $b\neq0$ one has
\begin{equation}
\label{eq:followed}
  \mathcal{A}^{(n,b)} \mathcal{A}^{(n,\frac1b)} = \mathcal{I}_{n+1}.
\end{equation}
The interesting pair of mutually inverse sequence transformations
\begin{equation}
  v_k = \frac1{n!}\sum_{j=0}^n E_{n,j}(-1,\tfrac1b;\,\tfrac1b-k,k) u_j
  \;\Longleftrightarrow\;
  u_k = \frac1{n!}\sum_{j=0}^n E_{n,j}(-1,b;\,b-k,k)  v_j
\end{equation}
follows from~(\ref{eq:followed}).  Such pairs were first derived
in~\cite{Charalambides82}.  However, deriving satisfactory
$n$\nobreakdash-dependent expressions for the elements of the matrix
$\mathcal{A}^{(n,b)}$ when $b\notin\mathbb{Z}$ is difficult.

\smallskip
Even outside the single-progression case, i.e., even when
$c_0+\nobreak c_\infty\neq b$, it may be possible to find a formula
for $E_{n,k}(a,b;c_0,c_\infty)$ that is not based on a sum, or at~most
involves one in which the number of terms does not grow with $n$
or~$k$.  Consider for instance the parametric triangle
$E_{n,k}(-1,2;c_0,0)$, for which the bivariate EGF
(see~(\ref{eq:speck})) is a manageable function of its arguments.
Versions of this triangle arise in several combinatorial contexts, and
formulas for its elements when $c_0=3$ and $c_0=1$ were given by Ma,
Ma, and Yeh~\cite{Ma2019}.

$E_{n,k}(-1,2;3,0)$ is the number of leaf-labeled rooted binary trees
with $n+\nobreak2$ leaves and $k+\nobreak1$ cherries, i.e., interior
vertices with exactly two descendant leaves.  (See
\cite[Table~6]{Rosenberg2019} and
\cite[$\texttt{A306364}$]{OEIS2022}.)  Also, the normalized version
$4^kE_{n,k}(-1,2;c_0,0)/(c_0)^{\overline n}$, when $c_0=3$, resp.~$1$,
has as its $n$'th row the $\gamma$\nobreakdash-vector of a simplicial
complex dual to the associahedron of type~$A_n$, resp.~$B_n$.  (See
\cite{Fomin2007} and
\cite[\texttt{A055151},\,\texttt{A105868}]{OEIS2022}.)  The element
$4^{n-k}E_{n,n-k}(-1,2;3,0)/(3)^{\underline n}$ of the reflected
triangle is the number of Motzkin paths of semi-length~$n$ with
$k$~steps, of types $U=\nobreak(1,1)$ and $H=\nobreak(1,0)$.  (See
\cite[$\texttt{A107131}$]{OEIS2022}.)

\begin{theorem}
  The triangles\/ $E_{n,k}(-1,2;c_0,0)$ and\/ $E_{n,k}(-1,2;c_0+1,1)$ have the
  hypergeometric term representations
  \begin{equation*}
    (c_0)^{\overline n}
    \left[ \myatop{-\frac{n}2,\: -\frac{n}2+\frac12}{1,\:\frac{c_0}2+\frac12}\right]^{\overline k},
    \quad\textrm{resp.}\quad
    (c_0+1)^{\overline n}
    \left[ \myatop{-\frac{n}2,\: -\frac{n}2-\frac12}{1,\:\frac{c_0}2+\frac12}\right]^{\overline k},
  \end{equation*}
  which are nonzero only if\/ $0\le k\le
  \left\lfloor\frac{n}2\right\rfloor$, resp.\ $0\le k\le
  \left\lfloor\frac{n+1}2\right\rfloor$.  Equivalently, their\/ n'th
  row polynomials have the hypergeometric representations
  \begin{equation*}
    (c_0)^{\overline n}\:
    {}_2F_1\biggl({\myatop{-\frac{n}2,\: -\frac{n}2+\frac12}{\frac{c_0}2+\frac{1}2}}
    \biggm|t\biggr),
    \quad\textrm{resp.}\quad
    (c_0+1)^{\overline n}\:
    {}_2F_1\biggl({\myatop{-\frac{n}2,\: -\frac{n}2-\frac12}{\frac{c_0}2+\frac{1}2}}
    \biggm|t\biggr).
  \end{equation*}
\label{thm:hyphypreps}
\end{theorem}
\begin{proof}
  The expressions for $E_{n,k}(-1,2;c_0,0)$ and
  $E_{n,k}(-1,2;c_0+1,1)$ can be obtained by repeatedly
  differentiating the bivariate EGF~(\ref{eq:speck}), or verified by
  by confirming that they satisfy the GKP recurrence of
  Definition~\ref{def:Edef}.  Also, the two are equivalent: the second
  comes from the first by a single application of
  Theorem~\ref{thm:2contigsE}(iv), the right-trimming identity.
\end{proof}

Beginning with $E_{n,k}(-1,2;c_0,0)$, one can apply
Theorem~\ref{thm:2contigsE}(i) repeatedly, so as to decrement~$c_0$
and increment~$c_\infty$ by any desired positive integer.  In this
way, for all $c_\infty\in\mathbb{N}$ one can obtain a formula for
$E_{n,k}(-1,2;\,c_0,c_\infty)$ in which the number of terms increases
with~$c_\infty$, but not with $n$ or~$k$.

Experimentation along this line resulted in the following unusual
conjecture, which involves the Bessel numbers
$B_{\nu,\kappa}=S_{\nu,\kappa}(1,2;0)$ and the $r$\nobreakdash-Bessel
numbers $B^{(r)}_{\nu,\kappa}=S_{\nu,\kappa}(1,2;r)$,
rank\nobreakdash-$0$ formulas for which appeared in
Theorem~\ref{thm:5points}(iv) and eq.~(\ref{eq:altgenform}).  A~proof
and a combinatorial interpretation are currently lacking.

\begin{conjecture}
  For all\/ $p\in\mathbb{N}$ and\/ $\zeta\in\{0,1\}$, and all\/~$c$,
  \begin{align*}
    &E_{n,k}(-1,2;\,c+2p+\zeta,2p+\zeta) = \\
    &\qquad
    \frac{(c+2p+\zeta)^{\overline n}}{(c+2p+2\zeta)^{{\overline p},2}}
    \sum_{\ell=0}^p
    \frac{(c+1)^{{\overline p},2}}{(c+1)^{{\overline {k+\ell}},2}}
    B_{2p+\zeta,2p+\zeta-\ell} B_{n,n-k}^{(2p+\zeta-2\ell)}.
  \end{align*}
\end{conjecture}

In addition to $E_{n,k}(-1,2;c_0,c_\infty)$, the parametric triangle
$E_{n,k}(-2,1;c_0,c_\infty)$ may be worth investigating, though not
many combinatorial applications of~it seem to be known.  It is not
difficult to derive the striking pair of identities
\begin{subequations}
\begin{align}
\label{eq:lesssat2}
  E_{n+1,k}(-2,1;\,c_0,0) &= c_0\,n!\,[t^k]P_{n}^{(c_0+n+1, -c_0-n-1)}(-t),\\
  E_{n,k}(-2,1;\,c_0+1,-1) &= n!\,[t^k]P_n^{(c_0+n, -c_0-n)}(-t),
\label{eq:lesssat1}
\end{align}
\end{subequations}
where $P_n^{(A,B)}(t)$ is the degree\nobreakdash-$n$ Jacobi
polynomial.  The two are equivalent, as the second is a right-trimmed
version of the first (with $c_0$ decremented by~$1$).  The second is
proved by verifying that the GKP row polynomial coming from its
left-hand side,
\begin{equation}
\label{eq:canident}
G_n(t)=\left[\begin{array}{cc|c}2,&1&c_0+1\\-1,&-1&-1\end{array}\right]_n(t),
\end{equation}
satisfies the same differential recurrence on~$n$ (see
Theorem~\ref{thm:pde}(ii)) as the known recurrence satisfied by
$n!\,P_n^{(c_0+n, -c_0-n)}(-t)$.  Alternatively a generating function
proof could be used, as the ordinary generating function of the
sequence $P_n^{(c_0+n, -c_0-n)}(t)$, $n\ge0$, is
known~\cite{Srivastava73}.

It follows from the Jacobi-polynomial representation
(\ref{eq:lesssat1}) that if $c_0=\nobreak\frac12$, the generalized
Eulerian row polynomial $G_n(t)$ of~(\ref{eq:canident}) can be
identified with $n!  \,P_n(-t)$, where $P_n$~is the $n$'th Boros--Moll
polynomial, originally introduced in the study of a quartic integral.
(See \cite{Boros2004} and \cite[$\texttt{A223549}$]{OEIS2022}.)  But
from a computational point of view, (\ref{eq:lesssat2}) and
(\ref{eq:lesssat1}) are weaker results than the row-polynomial
formulas of Theorem~\ref{thm:hyphypreps}: the coefficient of~$t^k$ in
$P_n^{(c_0+n,-c_0-n)}(t)$ can only be expressed as a sum, the number
of terms in~which grows with~$k$, as $n$~increases.

\section{Generalized Narayana triangles}
\label{sec:narayana}

In Section~\ref{sec:characteristics}, several cases when GKP
recurrences can be solved in closed form were introduced.  Case~(A)
was the generalized Stirling--Eulerian case, which gave rise in
Section~\ref{sec:stirlingeulerian} to the generalized Stirling and
Eulerian numbers, and many related identities.  Case~(B) will now be
examined.  As will be seen, the GKP triangles in case~(B) include many
triangles with combinatorial interpretations, including two
now-standard triangles of Narayana numbers~\cite{Petersen2015}.  This
is the reason for calling~(B) the generalized Narayana case.

In case (B\,I), when $(r_0,r_1,r_\infty) = (-\frac12,-\frac12,2)$,
Theorem~\ref{thm:BI} supplies an expression for the bivariate EGF
$G(t,z)$, based on a quadratic irrationality.  Much as in case~(A),
cases (B\,II) and~(B\,III) are obtained from (B\,I) by applying the
appropriate elements of the $S_3$\nobreakdash-group: the row-wise
sequence transformations $\textrm{RT}$ and
$\textrm{RT}\circ\textrm{UBT}\circ\textrm{RT}$.  The resulting vectors
$(r_0,r_1,r_\infty)$ are the permutations $(2,-\frac12,-\frac12)$ and
$(-\frac12,2,-\frac12)$.  Due to the analogy with the subcases of
case~(A), cases (B\,I), (B\,II), and (B\,III) will be called the
Stirling, reversed Stirling, and Eulerian subcases of the generalized
Narayana triangle.

It will be recalled that the new-style parameters $(r_0,r_1,r_\infty)$
and $(g_0,g_1,g_\infty)$, where the sums $r_0+r_1+r_\infty$ and
$g_0+g_1+g_\infty$ are constrained to equal $1$ and~$0$ respectively,
can be converted to the traditional GKP parameters
$(\alpha,\beta,\gamma;\allowbreak \alpha',\beta',\gamma')$ and vice
versa, by the formulas in (\ref{eq:old2new}) and~(\ref{eq:new2old}).
In the following, $3$\nobreakdash-parameter generalized Narayana
triangles $N_{n,k}^{\rm X}$ (where
$\textrm{X}=\allowbreak\textrm{S},\textrm{rS},\textrm{E}$, referring
to (B\,I), (B\,II), (B\,III)) are defined, in both the traditional and
new parametrizations.

\begin{definition}
\label{def:Bdefs}
\hfil\break
  (B\,I) The generalized Narayana triangle $N_{n,k}^{\rm S}(b;c_0,c_\infty)$, of Stirling type, is defined by
  \begin{equation*}
    N_{n,k}^{\rm S} = N_{n,k}^{\rm S}(b;c_0,c_\infty) \defeq
  \left[
    \begin{array}{cc|c}
      b/2, & b & c_0\\
      -b,& -b & c_\infty
    \end{array}
    \right]_{n,k},
  \end{equation*}
which if $b\neq0$ equals
\begin{equation*}
b^n
  \left[
    \begin{array}{ccc}
      0, & 1, & \infty \\
      \hline
      -\frac12, & -\frac12, & 2 \\
      c_0/b, & \,-(c_0+c_\infty)/b\,, & c_\infty/b
    \end{array}
    \right]_{n,k}.
\end{equation*}
  (B\,II) The generalized Narayana triangle $N_{n,k}^{\rm rS}(b;c_0,c_\infty)$, of reversed Stirling type, is defined by
  \begin{equation*}
    N_{n,k}^{\rm rS} = N_{n,k}^{\rm rS}(b;c_0,c_\infty) \defeq
  \left[
    \begin{array}{cc|c}
      -2b, & b & c_0\\
      3b/2,& -b & c_\infty
    \end{array}
    \right]_{n,k},
  \end{equation*}
which if $b\neq0$ equals
\begin{equation*}
b^n
  \left[
    \begin{array}{ccc}
      0, & 1, & \infty \\
      \hline
      2, & -\frac12, & -\frac12 \\
      c_0/b, & \,-(c_0+c_\infty)/b\,, & c_\infty/b
    \end{array}
    \right]_{n,k}.
\end{equation*}
  (B\,III) The generalized Narayana triangle $N_{n,k}^{\rm E}(b;c_0,c_\infty)$, of Eulerian type, is defined by
  \begin{equation*}
    N_{n,k}^{\rm E} = N_{n,k}^{\rm E}(b;c_0,c_\infty) \defeq
  \left[
    \begin{array}{cc|c}
      b/2, & b & c_0\\
      3b/2,& -b & c_\infty
    \end{array}
    \right]_{n,k},
  \end{equation*}
which if $b\neq0$ equals
\begin{equation*}
b^n
  \left[
    \begin{array}{ccc}
      0, & 1, & \infty \\
      \hline
      -\frac12, & 2, & -\frac12 \\
      c_0/b, & \,-(c_0+c_\infty)/b\,, & c_\infty/b
    \end{array}
    \right]_{n,k}.
\end{equation*}
\end{definition}

The GKP recurrences satisfied by these triangles are manifest in their
definition.  Their common homogeneity properties include
\begin{equation}
  N_{n,k}^{\rm X}(\lambda b;\,\lambda c_0,\lambda c_\infty) =
  \lambda^n N_{n,k}^{\rm X} (b;\,c_0,c_\infty),\qquad {\rm X}={\rm S},{\rm rS},{\rm E}.
\end{equation}
On account of this homogeneity, one can restrict without loss of
generality to a single nonzero value of~$b$.  To~facilitate comparison
with the standard Narayana triangles, the choice $b=2$ will be
convenient.  In what follows, the bivariate EGF's $G(t,z)$ of the
three types will be denoted by $N^{\rm X}(b;c_0,c_\infty;t,z)$, and
the corresponding row polynomials $G_n(t)$ by $N_n^{\rm
  X}(b;c_0,c_\infty;t)$.

\begin{theorem}
\label{thm:lessmanageable}
  The following EGF formulas hold in a neighborhood of\/ $(0,0)$.  
  \hfil\break 
  {\rm(B\,I)}, i.e., $(r_0,r_1,r_\infty)=(-\frac12,-\frac12,2)$, with\/ $b=2$:
  \begin{align*}
    &
    N^{\rm S}(2;\,c_0,c_\infty;\,t,z) =
  \left[
    \begin{array}{cc|c}
      1, & 2 & c_0\\
      -2,& -2 & c_\infty
    \end{array}
    \right](t,z)    
  =
  \left(\frac{s_+}{t_+}\right)^{c_0/2}
  \left(\frac{s_-}{t_-}\right)^{-(c_0+c_\infty)/2},
  \\[\jot]
  &
  \quad\quad
  s_\pm=\frac12\pm \frac{(2t-1)+z}{2\sqrt{1+2(2t-1)z+z^2}},
  \qquad
  t_\pm=\frac12\pm\left(t-\frac12\right).
  \end{align*}
  {\rm(B\,II)}, i.e., $(r_0,r_1,r_\infty)=(2,-\frac12,-\frac12)$, with\/ $b=2$:
  \begin{align*}
    &
    N^{\rm rS}(2;\,c_0,c_\infty;\,t,z) =
  \left[
    \begin{array}{cc|c}
      -4, & 2 & c_0\\
      3,& -2 & c_\infty
    \end{array}
    \right](t,z)    
  =
  \left(\frac{s_+}{t_+}\right)^{c_\infty/2}
  \left(\frac{s_-}{t_-}\right)^{-(c_0+c_\infty)/2},
  \\[\jot]
  &
  \quad\quad
  s_\pm=\frac12\pm \frac{(2-t)+t^2z}{2t\sqrt{1+2(2-t)z+t^2z^2}},
  \qquad
  t_\pm=\frac12\pm\frac{2-t}{2t}.
  \end{align*}
  {\rm(B\,III)}, i.e., $(r_0,r_1,r_\infty)=(-\frac12,2, -\frac12)$, with\/ $b=2$:
  \begin{align*}
    &
    N^{\rm E}(2;\,c_0,c_\infty;\,t,z) =
  \left[
    \begin{array}{cc|c}
      1, & 2 & c_0\\
      3,& -2 & c_\infty
    \end{array}
    \right](t,z)    
  =
  \left(\frac{s_+}{t_+}\right)^{c_0/2}
  \left(\frac{s_-}{t_-}\right)^{c_\infty/2},
  \\[\jot]
  &
  \quad\quad
  s_\pm=\frac12\pm \frac{(1+t)-(1-t)^2z}{2(t-1)\sqrt{1-2(1+t)z + (1-t)^2z^2}},
  \qquad
  t_\pm=\frac12\pm\frac{t+1}{2(t-1)}.
  \end{align*}
\end{theorem}
\begin{proof}
  The formula for $N^{\rm S}(2;c_0,c_\infty;t,z)$ is that of
  Theorem~\ref{thm:BI}, scaled by the factor $b=2$.  The
  $\textrm{X}=\textrm{rS},\textrm{E}$ formulas come by replacing
  $(t,z)$ by $(\frac1t,tz)$ and
  $\bigl(\frac{-t}{1-t},\allowbreak(1-\nobreak t)z\bigr)$
  respectively, as stated in Theorem~\ref{thm:38}.
\end{proof}

For general $c_0,c_\infty$, it is difficult to extract formulas for
the generalized Narayana numbers $N^{\rm X}_{n,k}(2;c_0,c_\infty)$
from the explicit EGF's in this theorem.  However, each of the three
EGF's satisfies a `contiguity relation.'  For instance, it follows by
elementary algebra that
\begin{equation}
\label{eq:unexplored}
\begin{split}
&\left[(2t-1)+z\right]N^{\rm S}(2;c_0,c_\infty;t,z)=\\
&\qquad 2t\,N^{\rm S}(2;c_0+1,c_\infty;t,z) - N^{\rm S}(2;c_0-1,c_\infty+2;t,z).
\end{split}
\end{equation}
This relates the parametric triangle $\genvert{n}k = N^{\rm
  S}_{n,k}(2;c_0,c_\infty)$ at three contiguous values of the
parameter-pair $(c_0,c_\infty)$.  The space of such relations and the
possibility of iterating them remain to be
explored.\footnote{Three-term relations resembling
  eq.~(\ref{eq:unexplored}) can be derived far more generally:
  in~fact, from any GKP recurrence in which one of the normalized
  parameters $r_0,r_1,r_\infty$ equals $1$~or~$2$.  An explicit EGF is
  not needed.  The derivation employs the PDE satisfied by the EGF,
  eq.~(\ref{eq:thepde}).}

By examination, there are three parametric restrictions (i.e.,
constraints on the pair $(c_0,c_\infty)$) which simplify each EGF
$N^{\rm X}(2;c_0,c_\infty;t,z)$ sufficiently that a closed-form
expression for the triangle entries $N^{\rm X}_{n,k}(2;c_0,c_\infty)$
can be obtained.  They will be denoted by (a),(b),(c).  The resulting
specialized EGF's are given in Theorem~\ref{thm:provides} below.

Recall that cases (B\,I),(B\,II),\allowbreak(B\,III) are really the
same, up~to permutation of the points $0,1,\infty$ of the projective
line: they constitute the case~(B), when the unordered set of
parameters $\{r_0,r_1,r_\infty\}$ is $\{-\frac12,-\frac12,2\}$.  The
restrictions (a),(b),(c) can be viewed as constraining the unordered
set of parameter-pairs $\{(r_0,g_0),\allowbreak (r_1,g_1), \allowbreak
(r_\infty,g_\infty)\}$.  For~(a), it must be of the form
$\{(-\frac12,g),\allowbreak (-\frac12,g),\allowbreak (2,-2g)\}$;
for~(b), of the form $\{(-\frac12,\frac12+h),\allowbreak
(-\frac12,\frac12-h),\allowbreak (2,-1)\}$; and for~(c), of the form
$\{(-\frac12,h),\allowbreak (-\frac12,-h),\allowbreak (2,0)\}$.

\begin{theorem}
\label{thm:provides}
  The following parametrically restricted EGF formulas hold in a
  neighborhood of\/ $(0,0)$.

  \noindent
  {\rm(B\,I)}, i.e., $(r_0,r_1,r_\infty)=(-\frac12,-\frac12,2)$, with\/ $b=2$:
  \begin{equation*}
  \begin{array}{lll}
  \mathrm{(a)} & g_0=g_1\mathrm{:} & N^{\rm S}(2;\,c,-2c;\,t,z) = S^{-c},    \vphantom{\left[\left(\frac{1-t}{t}\right)\frac{S-1+2t+z}{S+1-2t-z}\right]}\\
  \mathrm{(b)} & g_\infty=-1\mathrm{:} & N^{\rm S}(2;\,c,-2;\,t,z) = S^{-1}\left[\frac{-1+2t+z+S}{2t}\right]^{c-1},    \vphantom{\left[\left(\frac{1-t}{t}\right)\frac{S-1+2t+z}{S+1-2t-z}\right]^{(c-1)/2}}\\
  \mathrm{(c)} & g_\infty=0\mathrm{:} & N^{\rm S}(2;\,c-1,0;\,t,z) = \left[\left(\frac{t-1}{t}\right)\frac{1-2t-z-S}{1-2t-z+S}\right]^{(c-1)/2},
   \end{array}
  \end{equation*}
   where\/ $S\defeq\sqrt{1+2(2t-1)z+z^2}$.

  \smallskip
  \noindent
  {\rm(B\,II)}, i.e., $(r_0,r_1,r_\infty)=(2,-\frac12,-\frac12)$, with\/ $b=2$:
  \begin{equation*}
  \begin{array}{lll}
  \mathrm{(a)} & g_1=g_\infty\mathrm{:} & N^{\rm rS}(2;\,-2c,c;\,t,z) = S^{-c},  \vphantom{\left[(1-t)\frac{2-t+t^2z+tS}{2-t+t^2z-tS}\right]}\\
  \mathrm{(b)} & g_0=-1\mathrm{:} & N^{\rm rS}(2;\,-2,c;\,t,z) = S^{-1}\left[\frac{2-t+t^2z+tS}{2}\right]^{c-1},    \vphantom{\left[(1-t)\frac{2-t+t^2z+tS}{2-t+t^2z-tS}\right]^{(c-1)/2}}\\
  \mathrm{(c)} & g_0=0\mathrm{:} & N^{\rm rS}(2;\,0,c-1;\,t,z) = \left[(1-t)\frac{2-t+t^2z+tS}{2-t+t^2z-tS}\right]^{(c-1)/2},
   \end{array}
  \end{equation*}
   where\/ $S\defeq\sqrt{1+2(2-t)z+t^2z^2}$.

  \smallskip
  \noindent
  {\rm(B\,III)}, i.e., $(r_0,r_1,r_\infty)=(-\frac12,2,-\frac12)$, with\/ $b=2$:
  \begin{equation*}
  \begin{array}{lll}
  \mathrm{(a)} & g_0=g_\infty\mathrm{:} & N^{\rm E}(2;\,c,c;\,t,z) = S^{-c},  \vphantom{\left[(1-t)\frac{2-t+t^2z+tS}{2-t+t^2z-tS}\right]}\\
  \mathrm{(b)} & g_1=-1\mathrm{:} & N^{\rm E}(2;\,c,2-c;\,t,z) = S^{-1}\left[\frac{1+t-(1-t)^2z-(1-t)S}{2t}\right]^{c-1},    \vphantom{\left[((1-t)\frac{2-t+t^2z+tS}{2-t+t^2z-tS}\right]^{(c-1)/2}}\\
  \mathrm{(c)} & g_1=0\mathrm{:} & N^{\rm E}(2;\,c-1,1-c;\,t,z) = \left[\left(\frac1t\right)\frac{1+t-(1-t)^2z-(1-t)S}{1+t-(1-t)^2z+(1-t)S}\right]^{(c-1)/2},
   \end{array}
  \end{equation*}
   where\/ $S\defeq\sqrt{1-2(1+t)z+(1-t)^2z^2}$.
\end{theorem}

The EGF formulas of this theorem are more manageable than the general
ones of Theorem~\ref{thm:lessmanageable}.  With some effort, one can
extract from them the explicit formulas for $1$\nobreakdash-parameter
generalized Narayana numbers shown in Table~\ref{tab:2}.

Surprisingly, in each case the triangle entry $N_{n,k}^{\rm X}$ can be
represented as a hypergeometric term, parametrized by~$n,k$.  Because
of this, for any of the parametric triangles of Table~\ref{tab:2},
each row polynomial is a Gauss-hypergeometric polynomial, and could
optionally be rewritten in~terms of a Jacobi polynomial.  It is
straightforward to confirm each formula in the table, by verifying
that it satisfies the appropriate GKP recurrence.

A~few representations are omitted from the table because they are not
`pure,' in that they comprise more than a single hypergeometric term.
For instance, one can prove by induction that
\begin{sizeequation}{\small}
\label{eq:sizedtofit}
  N_{n,k}^{\rm E}(2;\,c-1,1-c) =
  (c-1)^{\overline n}
  \left\{
  \left[
    \myatop{-n+1,\,-n+c}{1,\,c}
    \right]^{\overline k}
   -
    \left[
    \myatop{-n+1,\,-n+c}{1,\,c}
    \right]^{\overline {k-1}}
  \right\},
\end{sizeequation}
when $1\le k\le n$.  This formula illustrates the operation of
mid-trimming.  The $n$'th row polynomial of the parametric triangle on
the left-hand side, for all $n\ge1$, is divisible by~$1\nobreak-t$;
whence the two terms on the right-hand side.

Many more trimming relationships could be mentioned.  In each of the
three sections of the table, triangle~(b) is a trimmed version of
triangle~(c); so the formulas given for the latter could be viewed as
redundant.  The trimming is respectively a right-, a left-, and a
mid-trimming, and (\ref{eq:sizedtofit}) is the identity resulting from
the last.

Moreover, in each of the three sections the $c=0$ and $c=2$ cases of
the (b)~triangle can be trimmed into the $c=3$ case of the
(a)~triangle.  The trimmings are respectively left- and mid-trimming;
right- and mid-trimming; and left- and right-trimming.  Also worth
noting is that in each section, the $c=1$ cases of the (a) and~(b)
triangles coincide.

A final observation regarding Table~\ref{tab:2} is that the triangles
in sections (B\,I), (B\,II), and (B\,III) are respectively invariant
under the involutive $S_3$\nobreakdash-group operations
$\textrm{UBT}$, $\textrm{RT}\circ\textrm{UBT}\circ\textrm{RT}$,
and~$\textrm{RT}$, each of which acts row-wise.  The last invariance
is manifest: the substitution $k\leftarrow n-\nobreak k$ leaves
(B\,III)(a) invariant and merely replaces $c$ by~$2-\nobreak c$ in
(B\,III)(b).  Up~to parametrization, $\textrm{RT}$~acts in cycle
notation as the permutation (I,II)(III), $\textrm{UBT}$~as
(II,III)(I), and $\textrm{RT}\circ\textrm{UBT}\circ\textrm{RT}$
as~(I,III)(II).  Although these operations permute the three sections
of the table, they do not affect the specialization letters
(a),(b),(c).

Many versions of the generalized Narayana triangles of
Table~\ref{tab:2} appear in the OEIS~\cite{OEIS2022}, in normalized
forms appropriate for combinatorial applications.  A~list is given in
Table~\ref{tab:3}, showing the value taken in each relevant OEIS entry
by the single parameter~$c$ of the corresponding Narayana triangle.
As indicated, in each OEIS entry the triangle elements~$\genvert{n}k$
are divided by a certain rising factorial, to reduce the triangle to
lowest terms.  For the (a)~specializations, this factor is
$(1)^{\overline n},(1)^{\overline n},(3)^{\overline n}$ when
$c=1,2,3$, and for the (b)~ones, it is $(2)^{\overline
  n},(1)^{\overline n},(2)^{\overline n}$ when $c=0,1,2$.  Also, when
${\textrm X}={\textrm S},{\textrm{rS}},{\textrm E}$, the triangle
elements $N^{\rm X}_{n,k}$ have respective signs
$(-)^k,(-)^{n-k},(+)$, and in the OEIS the negative signs are omitted.

For each generalized Narayana triangle in the OEIS, a compact
hypergeometric representation of its $n$'th row polynomial is shown in
the table.  (It~comes from the hypergeometric term representation
of~$\genvert{n}k$ given in Table~\ref{tab:2}, altered to agree with
the normalization used in the OEIS\null.)  The formulas listed for the
$c=0$ cases of {B\,I(b)} and~{B\,II(b)} require comment.  In~both,
$-2n$~is the lower parameter of the~${}_2F_1$, which would seem to
cause a division by zero if~$n=0$; but $-n$~being an upper parameter,
each ${}_2F_1$ is interpreted as unity if~$n=0$.

From a combinatorial point of view, the most important triangles in
Table~\ref{tab:3} may be the $c=1$ and $c=3$ cases of {B\,I(a)}
and~{B\,III(a)}.  Explicitly, they are
\begingroup
\begin{align}
  (-1)^k\left[(1)^{\overline n}\right]^{-1}N_{n,k}^{\rm S}(2;\,1,-2) & =\binom{n}k \binom{n+k}k,\\
  (-1)^k\left[(3)^{\overline n}\right]^{-1}N_{n,k}^{\rm S}(2;\,3,-6) & =\binom{n}k \binom{n+k+2}k\bigm/(k+1),\\
  \shortintertext{and}
  \left[(1)^{\overline n}\right]^{-1}N_{n,k}^{\rm E}(2;\,1,1) & ={\binom{n}k}^2,\label{eq:thirdoffour}\\
  \left[(3)^{\overline n}\right]^{-1}N_{n,k}^{\rm E}(2;\,3,3) & =\binom{n}k \binom{n+1}k\bigm/(k+1).\label{eq:fourthoffour}
\end{align}
\endgroup For each~$n$, the $n$'th rows of these four triangles are
respectively the \emph{f}\nobreakdash-vectors of simplicial complexes
dual to the associahedra of types $B_n$ and~$A_n$, and the
corresponding \emph{h}\nobreakdash-vectors~\cite{Fomin2007}.  As with
the permutohedra briefly encountered in~\S\ref{subsec:31}, the
\emph{f}\nobreakdash-vectors are mapped to the
\emph{h}\nobreakdash-vectors by what is essentially
$\textrm{UBT}\circ\textrm{RT}$, or equivalently
$\textrm{RT}\circ\textrm{UBT}\circ\textrm{RT}$, the
\emph{h}\nobreakdash-vectors of these simplicial polytopes being
reflection-invariant (the Dehn--Sommerville symmetry).  The invariance
of (\ref{eq:thirdoffour}) and~(\ref{eq:fourthoffour}) under
$k\leftarrow n-\nobreak k$ is evident.

The triangles (\ref{eq:thirdoffour}) and~(\ref{eq:fourthoffour}) are
known outside the combinatorics of polytopes: they are the
now-standard Narayana number triangles of types $B$ and~$A$\null.  In
the normalized triangle~(\ref{eq:fourthoffour}), the element
$\genvert{n}k$~counts the non-crossing partitions of an ordered
$(n+\nobreak1)$\nobreakdash-set into $k+\nobreak1$ blocks; in the
normalized triangle~(\ref{eq:thirdoffour}), it counts `signed' or
type\nobreakdash-$B$ non-crossing partitions~\cite{Chen2011,Fomin2007}.

\begin{landscape}
{
\setlength{\tabcolsep}{12pt}
\begin{table}
  \begin{center}
    \small
    \begin{tabular}{llll}
      &  $N^{\rm X}_{n,k}(2;\,c_0,c_\infty)$ & ${}_2F_1$ term representation & reversed representation
      \\
    \noalign{\smallskip}\hline\hline
    \noalign{\vskip2pt}
    (B\,I) & & & \\
    \noalign{\vskip-4pt}
    (a) & $N_{n,k}^{\rm S}(2;\,c,-2c)$
    & 
    $(c)^{\overline n}
    \left[
    \begin{array}{cc}
      -n, & n+c\\
      1, & \frac12+\frac{c}2
    \end{array}
    \right]^{\overline{k}\vphantom{\bigm|}}$
    &
    $(-2c)^{\underline n,4}
    \left[
    \begin{array}{cc}
      -n, & -n+\frac12-\frac{c}2\\
      1, & -2n+1-c
    \end{array}
    \right]^{\overline {n-k}}$
    \\
    (b) & $N_{n,k}^{\rm S}(2;\,c,-2)$
    & 
    $(c)^{\overline n}
    \left[
    \begin{array}{cc}
      -n, & n+1\\
      1, & c
    \end{array}
    \right]^{\overline{k}\vphantom{\bigm|}}$
    &
    $(-2)^{\underline n,4}
    \left[
    \begin{array}{cc}
      -n, & -n+1-c\\
      1, & -2n
    \end{array}
    \right]^{\overline {n-k}}$
    \\
    (c) & $N_{n,k}^{\rm S}(2;\,c-1,0)$
    & 
    $(c-1)^{\overline n}
    \left[
    \begin{array}{cc}
      -n+1, & n\\
      1, & c
    \end{array}
    \right]_{\vphantom{\underline p}}^{\overline{k}\vphantom{\bigm|}}
    $
    &
    {\hfil\textrm{---}}
    \\
    \hline
    \noalign{\vskip2pt}
    (B\,II) & & & \\
    \noalign{\vskip-4pt}
    (a) & $N_{n,k}^{\rm rS}(2;\,-2c,c)$
    & 
    $(-2c)^{\underline n,4}
    \left[
    \begin{array}{cc}
      -n, & -n+\frac12-\frac{c}2\\
      1, & -2n+1-c
    \end{array}
    \right]^{\overline{k}\vphantom{\bigm|}}$
    &
    $(c)^{\overline n}
    \left[
    \begin{array}{cc}
      -n, & n+c\\
      1, & \frac12+\frac{c}2
    \end{array}
    \right]^{\overline {n-k}}$
    \\
    (b) & $N_{n,k}^{\rm rS}(2;\,-2,c)$
    & 
    $(-2)^{\underline n,4}
    \left[
    \begin{array}{cc}
      -n, & -n+1-c\\
      1, & -2n
    \end{array}
    \right]^{\overline{k}\vphantom{\bigm|}}$
    &
    $(c)^{\overline n}
    \left[
    \begin{array}{cc}
      -n, & n+1\\
      1, & c
    \end{array}
    \right]^{\overline {n-k}}$
    \\
    (c) & $N_{n,k}^{\rm rS}(2;\,0,c-1)$
    & 
    {\hfil\textrm{---}}
    &
    $(c-1)^{\overline n}
    \left[
    \begin{array}{cc}
      -n+1, & n\\
      1, & c
    \end{array}
    \right]_{\vphantom{\underline p}}^{\overline{n-k}\vphantom{\bigm|}}$
    \\
    \hline
    \noalign{\vskip2pt}
    (B\,III) & & & \\
    \noalign{\vskip-4pt}
    (a) & $N_{n,k}^{\rm E}(2;\,c,c)$
    & 
    $(c)^{\overline n}
    \left[
    \begin{array}{cc}
      -n, & -n+\frac12-\frac{c}2\\
      1, & \frac12+\frac{c}2
    \end{array}
    \right]^{\overline{k}\vphantom{\bigm|}}$
    &
    $(c)^{\overline n}
    \left[
    \begin{array}{cc}
      -n, & -n+\frac12-\frac{c}2\\
      1, & \frac12+\frac{c}2
    \end{array}
    \right]^{\overline {n-k}}$
    \\
    (b) & $N_{n,k}^{\rm E}(2;\,c,2-c)$
    & 
    $(c)^{\overline n}
    \left[
    \begin{array}{cc}
      -n, & -n-1+c\\
      1, & c
    \end{array}
    \right]^{\overline{k}\vphantom{\bigm|}}$
    &
    $(2-c)^{\overline n}
    \left[
    \begin{array}{cc}
      -n, & -n+1-c\\
      1, & 2-c
    \end{array}
    \right]^{\overline {n-k}}$
    \\
    (c) & $N_{n,k}^{\rm E}(2;\,c-1,1-c)$
    & 
    $\vphantom{(c)^{\overline n}
    \left[
    \begin{array}{cc}
      -n, & n+1\\
      1, & c
    \end{array}
    \right]^{\overline {n-k}}}$
    {\hfil\textrm{---}}
    &
    {\hfil\textrm{---}}
    \end{tabular}
  \end{center}
   \caption{Certain 1-parameter generalized Narayana triangles with
     hypergeometric-term representations.  There are three each of the
     Stirling, reversed Stirling, and Eulerian kinds, i.e., three
     parametrically restricted versions of GKP cases (B\,I), (B\,II),
     and (B\,III).}
\label{tab:2}
\end{table}
}
\end{landscape}

\begin{landscape}
{
\setlength{\extrarowheight}{1.7pt}
\setlength{\tabcolsep}{12pt}
\begin{table}
  \begin{center}
    \begin{tabular}{@{}lllll}
      B\,I(a), $c=1$ & $(-1)^k[(1)^{\overline n}]^{-1} N_{n,k}^{\rm S}(2;1,-2)$ & ${}_2F_1(-n,n+1;1;-t)$ & \texttt{A063007} \\
      B\,I(a), $c=2$ & $(-1)^k[(1)^{\overline n}]^{-1} N_{n,k}^{\rm S}(2;2,-4)$ & $(n+1){}_2F_1(-n,n+2;\frac32;-t)$ & \texttt{A053124} \\
      B\,I(a), $c=3$ & $(-1)^k[(3)^{\overline n}]^{-1} N_{n,k}^{\rm S}(2;3,-6)$ & ${}_2F_1(-n,n+3;2;-t)$ & \texttt{A033282} \\
      \hdashline
      B\,I(b), $c=0$ & $(-1)^k[(2)^{\overline n}]^{-1} N_{n,k}^{\rm S}(2;0,-2)$ & $[(2)^{\overline n,4}/(2)^{\overline n}]t^{n}{}_2F_1(-n,-n+1;-2n;-\frac1t)$ & \texttt{A086810} \\
      B\,I(b), $c=1$ & $(-1)^k[(1)^{\overline n}]^{-1} N_{n,k}^{\rm S}(2;1,-2)$ & ${}_2F_1(-n,n+1;1;-t)$ & \texttt{A063007} \\
      B\,I(b), $c=2$ & $(-1)^k[(2)^{\overline n}]^{-1} N_{n,k}^{\rm S}(2;2,-2)$ & ${}_2F_1(-n,n+1;2;-t)$ & \texttt{A088617} \\
      \hline
      B\,II(a), $c=1$ & $(-1)^{n-k}[(1)^{\overline n}]^{-1} N_{n,k}^{\rm rS}(2;-2,1)$ & $t^n{}_2F_1(-n,n+1;1;-\frac1t)$ & \texttt{A104684} \\
      B\,II(a), $c=2$ & $(-1)^{n-k}[(1)^{\overline n}]^{-1} N_{n,k}^{\rm rS}(2;-4,2)$ & $(n+1)t^n{}_2F_1(-n,n+2;\frac32;-\frac1t)$ & \texttt{A053125} \\
      B\,II(a), $c=3$ & $(-1)^{n-k}[(3)^{\overline n}]^{-1} N_{n,k}^{\rm rS}(2;-6,3)$ & $t^n{}_2F_1(-n,n+3;2;-\frac1t)$ & \texttt{A126216} \\
      \hdashline
      B\,II(b), $c=0$ & $(-1)^{n-k}[(2)^{\overline n}]^{-1} N_{n,k}^{\rm rS}(2;-2,0)$ & $[(2)^{\overline n,4}/(2)^{\overline n}]{}_2F_1(-n,-n+1;-2n;t)$ & \texttt{A133336} \\
      B\,II(b), $c=1$ & $(-1)^{n-k}[(1)^{\overline n}]^{-1} N_{n,k}^{\rm rS}(2;-2,1)$ & $t^n{}_2F_1(-n,n+1;1;-\frac1t)$ & \texttt{A104684} \\
      B\,II(b), $c=2$ & $(-1)^{n-k}[(2)^{\overline n}]^{-1} N_{n,k}^{\rm rS}(2;-2,2)$ & $t^n{}_2F_1(-n,n+1;2;-\frac1t)$ & \texttt{A060693} \\
      \hline
      B\,III(a), $c=1$ & $[(1)^{\overline n}]^{-1} N_{n,k}^{\rm E}(2;1,1)$ & ${}_2F_1(-n,-n;1;t)$ & \texttt{A008459} \\
      B\,III(a), $c=2$ & $[(1)^{\overline n}]^{-1} N_{n,k}^{\rm E}(2;2,2)$ & $(n+1){}_2F_1(-n,-n-\frac12;\frac32;t)$ & \texttt{A091044} \\
      B\,III(a), $c=3$ & $[(3)^{\overline n}]^{-1} N_{n,k}^{\rm E}(2;3,3)$ & ${}_2F_1(-n,-n-1;2;t)$ & \texttt{A001263} \\
      \hdashline
      B\,III(b), $c=0$ & $[(2)^{\overline n}]^{-1} N_{n,k}^{\rm E}(2;0,2)$ & $t^n{}_2F_1(-n,-n+1;2,\frac1t)$ & \texttt{A090181} \\
      B\,III(b), $c=1$ & $[(1)^{\overline n}]^{-1} N_{n,k}^{\rm E}(2;1,1)$ & ${}_2F_1(-n,-n;1;t)$ & \texttt{A008459} \\
      B\,III(b), $c=2$ & $[(2)^{\overline n}]^{-1} N_{n,k}^{\rm E}(2;2,0)$ & ${}_2F_1(-n,-n+1;2;t)$ & \texttt{A131198}
    \end{tabular}
  \end{center}
  \caption{Generalized Narayana triangles in the OEIS~\cite{OEIS2022},
    with hypergeometric-polynomial representations of their row
    polynomials.\hfil\break  (N.B.: Triangles \texttt{A053124}, \texttt{A053125}
    in the OEIS are signed rather than signless, disagreeing with the
    convention adhered~to here.)\hfil\break  (N.B.: When $n=0$, each ${}_2F_1$
    equals unity, by convention if necessary.)}
\label{tab:3}  
\end{table}
}
\end{landscape}

\section{Generalized secant--tangent triangles}
\label{sec:sectang}

In addition to the generalized Stirling--Eulerian case (A) and the
generalized Narayana case~(B), there is a third case when a GKP
recurrence can be solved in closed form, or at~least the bivariate EGF
$G(t,z)$ can be computed explicitly by the method of characteristics.
As noted in Section~\ref{sec:characteristics}, this is the generalized
secant--tangent case, case~(C).

The GKP triangles in case~(C) include important ones with
combinatorial interpretations, but they are relatively few.  The
following treatment will briefly relate their EGF's to previous work.
The term `generalized secant--tangent' will be justified.  It comes
from an alternative method of generating the row polynomials $G_n(t)$,
$n\ge0$, of any GKP triangle.  (See Theorem~\ref{thm:cfg}.)  This
could be called the iterated derivation method, and is a specific
application of the context-free grammar approach taken by
Chen~\cite{Chen93} and Dumont~\cite{Dumont96} to exponential
structures in combinatorics.  As another application, some final
identities involving the generalized Eulerian numbers
$E_{n,k}(a,b;c_0,c_\infty)$ will be derived.

In subcase (C\,I), when $(r_0,r_1,r_\infty) = (\frac12,\frac12,0)$,
Theorem~\ref{thm:CI} supplies a transcendental expression for the
bivariate EGF $G(t,z)$.  Much as with cases (A) and~(B), subcases
(C\,II) and~(C\,III) are obtained from (C\,I) by applying the
appropriate elements of the $S_3$\nobreakdash-group: the row-wise
sequence transformations $\textrm{RT}$ and
$\textrm{RT}\circ\textrm{UBT}\circ\textrm{RT}$.  The resulting vectors
$(r_0,r_1,r_\infty)$ are the permutations $(0,\frac12,\frac12)$ and
$(\frac12,0,\frac12)$.  Due to the analogy with the subcases of
case~(A), cases (C\,I), (C\,II), and (C\,III) will be called the
Stirling, reversed Stirling, and Eulerian subcases of the generalized
secant--tangent triangle.

As always, the parameters $(r_0,r_1,r_\infty)$ and
$(g_0,g_1,g_\infty)$, where the sums $r_0+r_1+r_\infty$ and $g_0+g_1+g_\infty$
are constrained to equal $1$ and~$0$ respectively, can be converted to
the traditional GKP parameters $(\alpha,\beta,\gamma;\allowbreak
\alpha',\beta',\gamma')$ and vice versa, by the formulas in
(\ref{eq:old2new}) and~(\ref{eq:new2old}).  In the following,
$3$\nobreakdash-parameter generalized secant--tangent triangles
$W_{n,k}^{\rm X}$ (where
$\textrm{X}=\allowbreak\textrm{S},\textrm{rS},\textrm{E}$, referring
to (C\,I), (C\,II), (C\,III)) are defined, in both the traditional and
new parametrizations.  The bivariate EGF's $G(t,z)$ of the three types
will be denoted by $W^{\rm X}(b;c_0,c_\infty;t,z)$,

\begin{definition}
\label{def:Cdefs}
\hfil\break
  (C\,I) The generalized secant--tangent triangle $W_{n,k}^{\rm S}(b;c_0,c_\infty)$, of Stirling type, is defined by
  \begin{equation*}
    W_{n,k}^{\rm S} = W_{n,k}^{\rm S}(b;c_0,c_\infty) \defeq
  \left[
    \begin{array}{cc|c}
      -b/2, & b & c_0\\
      b,& -b & c_\infty
    \end{array}
    \right]_{n,k},
  \end{equation*}
which if $b\neq0$ equals
\begin{equation*}
b^n
  \left[
    \begin{array}{ccc}
      0, & 1, & \infty \\
      \hline
      \frac12, & \frac12, & 0 \\
      c_0/b, & \,-(c_0+c_\infty)/b\,, & c_\infty/b
    \end{array}
    \right]_{n,k}.
\end{equation*}
  (C\,II) The generalized secant--tangent triangle $W_{n,k}^{\rm rS}(b;c_0,c_\infty)$, of reversed Stirling type, is defined by
  \begin{equation*}
    W_{n,k}^{\rm rS} = W_{n,k}^{\rm rS}(b;c_0,c_\infty) \defeq
  \left[
    \begin{array}{cc|c}
      0, & b & c_0\\
      b/2,& -b & c_\infty
    \end{array}
    \right]_{n,k},
  \end{equation*}
which if $b\neq0$ equals
\begin{equation*}
b^n
  \left[
    \begin{array}{ccc}
      0, & 1, & \infty \\
      \hline
      0, & \frac12, & \frac12 \\
      c_0/b, & \,-(c_0+c_\infty)/b\,, & c_\infty/b
    \end{array}
    \right]_{n,k}.
\end{equation*}
  (C\,III) The generalized secant--tangent triangle $W_{n,k}^{\rm E}(b;c_0,c_\infty)$, of Eulerian type, is defined by
  \begin{equation*}
    W_{n,k}^{\rm E} = W_{n,k}^{\rm E}(b;c_0,c_\infty) \defeq
  \left[
    \begin{array}{cc|c}
      -b/2, & b & c_0\\
      b/2,& -b & c_\infty
    \end{array}
    \right]_{n,k},
  \end{equation*}
which if $b\neq0$ equals
\begin{equation*}
b^n
  \left[
    \begin{array}{ccc}
      0, & 1, & \infty \\
      \hline
      \frac12, & 0, & \frac12 \\
      c_0/b, & \,-(c_0+c_\infty)/b\,, & c_\infty/b
    \end{array}
    \right]_{n,k}.
\end{equation*}
\end{definition}

\smallskip
These parametric triangles have the common homogeneity property
\begin{equation}
  W_{n,k}^{\rm X}(\lambda b;\,\lambda c_0,\lambda c_\infty) =
  \lambda^n W_{n,k}^{\rm X} (b;\,c_0,c_\infty),\qquad {\rm X}={\rm S},{\rm rS},{\rm E}.
\end{equation}
To facilitate comparison with previous work, the choice $b=2$ will now
be made, without loss of generality.

\begin{theorem}
\label{thm:lessmanageableC}
  The following EGF formulas hold in a neighborhood of\/ $(0,0)$.  
  \hfil\break 
  {\rm(C\,I)}, i.e., $(r_0,r_1,r_\infty)=(\frac12,\frac12,0)$, with\/ $b=2$:
  \begin{align*}
    &
    W^{\rm S}(2;\,c_0,c_\infty;\,t,z) =
  \left[
    \begin{array}{cc|c}
      -1, & 2 & c_0\\
      2,& -2 & c_\infty
    \end{array}
    \right](t,z)    
  =
  \left(\frac{s_+}{t_+}\right)^{c_0/2}
  \left(\frac{s_-}{t_-}\right)^{-(c_0+c_\infty)/2},
  \\[\jot]
  &
  \quad
  s_\pm= \left[\sqrt{t_\pm}\cos(z\sqrt{t_+t_-}) \pm \sqrt{t_\mp}\sin(z\sqrt{t_+t_-})\right]^2,
  \quad
  t_\pm=\frac12\pm\left(t-\frac12\right).
  \end{align*}
  {\rm(C\,II)}, i.e., $(r_0,r_1,r_\infty)=(0,\frac12,\frac12)$, with\/ $b=2$:
  \begin{align*}
    &
    W^{\rm rS}(2;\,c_0,c_\infty;\,t,z) =
  \left[
    \begin{array}{cc|c}
      0, & 2 & c_0\\
      1,& -2 & c_\infty
    \end{array}
    \right](t,z)    
  =
  \left(\frac{s_+}{t_+}\right)^{c_\infty/2}
  \left(\frac{s_-}{t_-}\right)^{-(c_0+c_\infty)/2},
  \\[\jot]
  &
  \quad
  s_\pm= \left[\sqrt{t_\pm}\cos(z\sqrt{t_-/t_+}) \pm \sqrt{t_\mp}\sin(z\sqrt{t_-/t_+})\right]^2,
  \quad
  t_\pm=\frac12\pm\frac{2-t}{2t}.
  \end{align*}
  {\rm(C\,III)}, i.e., $(r_0,r_1,r_\infty)=(\frac12,0, \frac12)$, with\/ $b=2$:
  \begin{align*}
    &
    W^{\rm E}(2;\,c_0,c_\infty;\,t,z) =
  \left[
    \begin{array}{cc|c}
      -1, & 2 & c_0\\
      1,& -2 & c_\infty
    \end{array}
    \right](t,z)    
  =
  \left(\frac{s_+}{t_+}\right)^{c_0/2}
  \left(\frac{s_-}{t_-}\right)^{c_\infty/2},
  \\[\jot]
  &
  \quad
  s_\pm= \left[\sqrt{t_\pm}\cos(z\sqrt{t_+/t_-}) \pm \sqrt{t_\mp}\sin(z\sqrt{t_+/t_-})\right]^2,
  \quad
  t_\pm=\frac12\pm\frac{t+1}{2(t-1)}.
  \end{align*}
It should be noted that in each case, $s_++s_-$ and\/ $t_++t_-$ equal
unity.
\end{theorem}
\begin{proof}
  The formula for $W^{\rm S}(2;c_0,c_\infty;t,z)$ is that of
  Theorem~\ref{thm:CI}, scaled by the factor $b=2$.  The
  $\textrm{X}=\textrm{rS},\textrm{E}$ formulas come by replacing
  $(t,z)$ by $(\frac1t,tz)$ and
  $\bigl(\frac{-t}{1-t},\allowbreak(1-\nobreak t)z\bigr)$
  respectively, as stated in Theorem~\ref{thm:38}.
\end{proof}

It has long been known that case\nobreakdash-(C) GKP recurrences, in
particular case-(C\,II) ones, appear in the enumerative combinatorics
of peaks and valleys of permutations.  (See the papers of
Ma~\cite{Ma2012} and Zhuang~\cite{Zhuang2016}, with earlier work
extending from Andr\'e in the 1880's through Entringer in the 1960's
and Gessel in the 1970's~\cite{Gessel77} remaining relevant.)
$W^{\texttt{rS}}_{n,k}(2;1,0)$ counts the elements of the group~$S_n$
that have $k$~left (or~right) peaks
\cite[{\texttt{A008971}}]{OEIS2022}, and
$W^{\texttt{rS}}_{n,k}(2;2,0)$ the elements of~$S_{n+1}$ that have
$k$~peaks \cite[{\texttt{A008303}}]{OEIS2022}.  (In~\cite{Ma2012},
these are denoted by $W_{n,k}^l$ and~$W_{n+1,k}$.)  Also, the triangle
$W^{\texttt{rS}}_{n,k}(2;0,1)$ counts the elements of~$S_n$ that have
$k$~left--right peaks~\cite{Zhuang2016}, but left-trimming this
triangle reduces it to $W^{\texttt{rS}}_{n,k}(2;2,0)$.  In a separate
combinatorial application, the normalized triangle
$2^{-n}4^kW^{\texttt{rS}}_{n,k}(2;2,0)$,
resp.\ $4^kW^{\texttt{rS}}_{n,k}(2;1,0)$, has as its $n$'th row the
$\gamma$\nobreakdash-vector of a simplicial complex dual to the
permutohedron of type $A_n$, resp.~$B_n$.  (See~\cite{Fomin2007}
and~\cite[\texttt{A101280}]{OEIS2022}.)

\begin{proposition}
  The following EGF formulas hold in a neighborhood of\/ $(0,0)$.
  \begingroup
  \begin{align}
    &
    \begin{aligned}
      W^{\rm rS}_{n,k}(2;\,1,0;\,t,z) &=
      \left[
        \begin{array}{cc|c}
          0, & 2 & 1\\
          1,& -2 & 0
        \end{array}
        \right](t,z)\\[\jot]
      &=
        \frac{\sqrt{1-t}}{\sqrt{1-t}\cosh(z\sqrt{1-t}) - \sinh(z\sqrt{1-t})}
      ,
    \end{aligned}
    \\
    \shortintertext{and}
    &
    \begin{aligned}
      W^{\rm rS}_{n,k}(2;\,2,0;\,t,z) &=
      \left[
        \begin{array}{cc|c}
          0, & 2 & 2\\
          1,& -2 & 0
        \end{array}
        \right](t,z)
      \\[-\jot]
      &= 
      \left[
        \frac{\sqrt{1-t}}{\sqrt{1-t}\cosh(z\sqrt{1-t}) - \sinh(z\sqrt{1-t})}
        \right]^2
      \\[\jot]
      &=
      \frac{{\rm d}}{{\rm d}z}
      \left\{
      \left[
        \sqrt{1-t}\coth(z\sqrt{1-t})-1
        \right]^{-1}\right\}.
    \end{aligned}
  \end{align}
  \endgroup
\end{proposition}

These EGF formulas follow by some trigonometric manipulation from the
one for $W^{\rm rS}(2;\,c_0,c_\infty;\,t,z)$ in
Theorem~\ref{thm:lessmanageableC}.  They agree with those previously
known \cite{Gessel77,Ma2012,Zhuang2016}, but the present derivation
places them in analytic context as individual EGF's that belong to a
family of EGF's that can be computed by the method of characteristics.
It should be noted that for case\nobreakdash-(C) GKP recurrences, it
is the bivariate triangle EGF's and not the triangle
elements~$\genvert{n}k$ for which explicit, closed-form expressions
are currently available.

The iterated derivation method of solving a GKP recurrence with
parameters $(\alpha,\beta,\gamma;\allowbreak \alpha',\beta',\gamma')$,
which is an instance of a more general grammar-based method of
combinatorial enumeration, will now be summarized.  Its applicability
is not restricted to case\nobreakdash-(C) recurrences.

Let $D$ be a formal derivation satisfying Leibniz's rule, which acts
on any reasonable function of the variables or indeterminates~$x,y$;
such as a polynomial or a formal power series, though non-integral
powers are allowed.  The following is a known fact
(cf.\ \cite[Theorem~2.1]{Hao2015} and \cite[Lemma~8]{Ma2019}),
restated in present notation.  It is equivalent to the differential
recurrence of Theorem~\ref{thm:pde}(ii) and the iterated operator
formula of Theorem~\ref{thm:pde}(iii).

\begin{theorem}
\label{thm:cfg}
  If the variables\/ $x,y$ have the respective derivations
  \begin{equation}
    \label{eq:derivations}
D(x,y) = (x,y)\ast \left(x^{\alpha}y^{\alpha'}\!,\, x^{\alpha+\beta} y^{\alpha'+\beta'}\right),
  \end{equation}
where\/ $\ast$ signifies the elementwise product, then for all\/ $n\ge0$,
\begin{equation}
  D^{n} \bigl(x^{\gamma}y^{\gamma'}\bigr) = \bigl(x^{\gamma} y^{\gamma'}\bigr) \, \bigl(x^\alpha y^{\alpha'}\bigr)^n\,
    \sum_{k=0}^n \left[     \begin{array}{ll|l}
      \alpha, & \beta & \gamma \\
      \alpha', & \beta' & \gamma' \\
    \end{array}
      \right]_{n,k}
 \bigl(x^\beta y^{\beta'}\bigr)^k,
\end{equation}
or equivalently
\begin{equation*}
  D^n\bigl(x^{\gamma} y^{\gamma'}\bigr) = \bigl(x^\gamma y^{\gamma'}\bigr)  \bigl(x^\alpha y^{\alpha'}\bigr)^n \,G_n(x^\beta y^{\beta'}),
\end{equation*}
where\/ $G_n=G_n(t)$ is the n'th row polynomial of the GKP triangle
with parameter array
$\left[\begin{smallarray}{cc|c}\alpha,&\beta&\gamma\\\alpha',&\beta'&\gamma'\end{smallarray}\right]$.
\end{theorem}

It may be possible to realize $x,y$ as functions of an independent
variable~$w$, and $D$~as the derivative $D_w\defeq {\rm d}/{{\rm
    d}w}$.  If so, the following is immediate: it comes by
exponentiating~$D_w$.

\begin{corollary}
  For all\/ $\delta$, or equivalently as an equality between formal power series in\/~$\delta$,
  \begin{equation*}
    \frac{\bigl( x^\gamma y^{\gamma'}\bigr)(w+\delta)}
         {\bigl( x^\gamma y^{\gamma'}\bigr)(w)}
    =
    G\left(  
    \bigl(x^\beta y^{\beta'}\bigr)(w), 
    \delta 
    \bigl(x^\alpha y^{\alpha'}\bigr)(w)
    \right),
  \end{equation*}
where\/ $G=G(t,z)$ is the bivariate EGF of the GKP triangle with
parameter array
$\left[\begin{smallarray}{cc|c}\alpha,&\beta&\gamma\\\alpha',&\beta'&\gamma'\end{smallarray}\right]$.
\end{corollary}

There is also an easily checked homogeneity property: if the pair
$x=x(w)$, $y=y(w)$ satisfy the differential
equation~(\ref{eq:derivations}), then for any~$\lambda$, so do the
pair $\lambda^p x(\lambda w)$, $\lambda^q y(\lambda w)$, where
\begin{equation}
  p=\frac{\beta'}{\alpha\beta'-\alpha'\beta}, \qquad
  q=\frac{-\beta}{\alpha\beta'-\alpha'\beta},
\end{equation}
it being assumed that $\alpha\beta'-\alpha'\beta\neq0$.  If
$\beta'=-\beta$, these reduce to
$\lambda^{1/(\alpha+\alpha')}x(\lambda w)$,
$\lambda^{1/(\alpha+\alpha')}y(\lambda w)$.  Moreover, $w$~can be
shifted without affecting~(\ref{eq:derivations}), as it is an
autonomous equation.  Thus (\ref{eq:derivations}) has a
$2$\nobreakdash-parameter space of solutions.

For the three case-(C) GKP triangles $W^{\rm X}_{n,k}(b;c_0,c_\infty)$
of Definition~\ref{def:Cdefs}, solutions $x,y$
of~(\ref{eq:derivations}) can be found by inspection.  When $b=2$,
canonical ones are the following.

\setlength\tabcolsep{12.0pt}
\begin{center}
\begin{tabular}{llll}
  & $\genvert{n}k$ & $(\alpha,\beta;\,\alpha',\beta')$ & $x,y$\\
      \noalign{\smallskip}\hline
  (C\,I): & $W_{n,k}^{\rm S}$, $b=2$ & $(-1,2;\,2,-2)$ & $x=\tan w$, $y=\sec w$\\
  (C\,II): & $W_{n,k}^{\rm rS}$, $b=2$ & $(0,2;\,1,-2)$ & $x=\sec w$, $y=\tan w$\\
  (C\,III): & $W_{n,k}^{\rm E}$, $b=2$ & $(-1,2;\,1,-2)$ & $x=\cosh w$, $y=\sinh w$
\end{tabular}
\end{center}

\noindent
In each case the pair $\gamma,\gamma'$ equals $c_0,c_\infty$, and
Theorem~\ref{thm:cfg} yields the three elegant identities
\setlength\tabcolsep{3.0pt}

\begin{subequations}
\begin{align}
  &
  \left(\tan^{c_0}w\,\sec^{c_\infty}w\right)^{-1}  D_w^n \left(\tan^{c_0}w\,\sec^{c_\infty}w\right)\nonumber\\
  &
  \qquad\qquad =
  (\sec^n w) \sum_{k=0}^n W_{n,k}^{\rm S}(2;\,c_0,c_\infty)\, \csc^{n-2k}w,\label{eq:penulta}\\
  &
  \left(\sec^{c_0}w\,\tan^{c_\infty}w\right)^{-1}  D_w^n \left(\sec^{c_0}w\,\tan^{c_\infty}w\right)\nonumber\\
  &
  \qquad \qquad =
  (\tan^n w) \sum_{k=0}^n W_{n,k}^{\rm rS}(2;\,c_0,c_\infty)\, \csc^{2k}w,\label{eq:penultb}\\
\shortintertext{and}
  &
  \left(\cosh^{c_0}w\,\sinh^{c_\infty}w\right)^{-1}  D_w^n \left(\cosh^{c_0}w\,\sinh^{c_\infty}w\right)\nonumber\\
  &
  \qquad \qquad =
  (\tanh^n w) \sum_{k=0}^n W_{n,k}^{\rm E}(2;\,c_0,c_\infty)\, \coth^{2k}w.\label{eq:penultc}
\end{align}
\end{subequations}

Under the reflection operation (RT), the parametric triangles
$W_{n,k}^{\rm S}$, $W_{n,k}^{\rm rS}$ are interchanged, with
$c_0\leftrightarrow c_\infty$.  Also triangle $W_{n,k}^{\rm E}$ is
invariant, except that again, $c_0\leftrightarrow c_\infty$.  It
should be noted that in~(\ref{eq:penultc}), $\sinh,\cosh$ (and
likewise $\tanh,\coth$) could optionally be interchanged, because
$x(w)=\cosh w$, $y(w)=\sinh w$ are converted to $x(w)=\sinh w$,
$y(w)=\cosh w$ by a complex shift of~$w$.

As mentioned, the triangles $W_{n,k}^{\rm rS}(2;1,0)$ and
$W_{n,k}^{\rm rS}(2;2,0)$ are of combinatorial significance.
By~(\ref{eq:penultb}), their elements are the coefficients of
trigonometric polynomials obtained by repeated differentiation of
$\sec w$ and~$\sec^2 w$, respectively.  This fact is
known~\cite{Ma2012}, but one now sees that by a straightforward
generalization one can generate any desired case\nobreakdash-(C) GKP
triangle, with or without a combinatorial interpretation.  This
explains the term `generalized secant--tangent triangle.'

Case-(A) and case-(B) triangles can also be generated by the iterated
derivation method.  It could be applied as follows in the generalized
Stirling--Eulerian case~(A\null).  To~treat subcase (A\,I) similarly
to~(C\,I), one should generate not the Hsu--Shiue triangle
$S_{n,k}(a,b;r)$ but rather $(s)^{\overline k}S_{n,k}(a,b;r)$ for
arbitrary~$s$, because it is the latter which satisfies a GKP
recurrence with $\beta\beta'\neq0$.  (Recall eq.~(\ref{eq:infact}).)
The matching (A\,II) triangle is the reflection
$(r)^{\overline{n-k}}S_{n,n-k}(a,b;s)$, and the (A\,III) triangle is
the generalized Eulerian triangle $E_{n,k}(a,b;c_0,c_\infty)$.  The
GKP parameters of these triangles, and canonical solutions $x,y$
of~(\ref{eq:derivations}) when $a\neq0$ and~$b\neq0$ that can be found
by inspection, are as follows.

\setlength\tabcolsep{5.5pt}
\begin{center}
{
\small
\begin{tabular}{llll}
  & $\genvert{n}k$ & $(\alpha,\beta;\,\alpha',\beta')$, $(\gamma,\gamma')$ & $x,y$ \\
      \noalign{\smallskip}\hline
  (A\,I): & $s^{\overline{k}} \vphantom{\Bigl(A\Bigr)}  S_{n,k}(a,b;\,r)$ & $(-a,b;\,0,1)$ & $x=(1+aw)^{1/a}$\\
      & & $(r,s)$ & $y=b\bigl[1-(1+aw)^{b/a}\bigr]^{-1}$  \\[6pt]
  (A\,II): & $r^{\overline{n-k}}S_{n,n-k}(a,b;\,s)$ & $(1,-1;\,-a+b,-b)$ & $x=b\bigl[1-(1+aw)^{b/a}\bigr]^{-1}$ \\[\jot]
      & & $(r,s)$ & $y=(1+aw)^{1/a}$  \\[4pt] 
  (A\,III): & $E_{n,k}(a,b;\,c_0,c_\infty)$ & $(-a,b;\,a+b,-b)$ & $x= \bigl[1-(1-aw)^{-b/a}\bigr]^{-1/b}$ \\[2pt]
      & & $(c_0,c_\infty)$ & $y=\bigl[ (1-aw)^{b/a}-1\bigr]^{-1/b}$
\end{tabular}
}
\end{center}

\noindent
As in case~(C), the parametric (A\,I) and (A\,II) triangles are
interchanged by the reflection operation, with $r\leftrightarrow s$,
and so are the solutions~$x,y$.  The (A\,III) triangle is invariant,
with $c_0\leftrightarrow c_\infty$; and the solutions $x,y$ could
optionally be interchanged, as in case~(C\,III\null).

The three corresponding identities come at once from
Theorem~\ref{thm:cfg}.  The (A\,II) identity is merely a reflected
version of the (A\,I) one and is left to the reader.  The (A\,I) and
(A\,III) identities are respectively
\begin{subequations}
  \begin{sizealign}{\small}
    &
    \begin{aligned}
      &
      \left\{
      (1+aw)^{r/a} \bigl[1-(1+aw)^{b/a}\bigr]^{-s}
      \right\}^{-1}
      \!
      D_w^n
      \left\{
      (1+aw)^{r/a} \bigl[1-(1+aw)^{b/a}\bigr]^{-s}
      \right\}
      \\
      &
      \qquad
      =
      (1+aw)^{-n}
      \sum_{k=0}^n b^ks^{\overline{k}}\,S_{n,k}(a,b;\,r)
      \left\{
      (1+aw)^{b/a} \bigl[1-(1+aw)^{b/a}\bigr]^{-1}
      \right\}^{k}
      \label{eq:Sabc}
    \end{aligned}
    \\
    \intertext{\normalsize{and}}
    &
    \begin{aligned}
      &
      \left\{
      (1-aw)^{c_0/a}
      \bigl[(1-aw)^{b/a}-1\bigr]^{-(c_0+c_\infty)/b}
      \right\}^{-1}
      \\
      & \qquad\qquad\quad {}\times
      D_w^n
      \left\{
      (1-aw)^{c_0/a}
      \bigl[(1-aw)^{b/a}-1\bigr]^{-(c_0+c_\infty)/b}
      \right\}
      \\
      &
      \qquad
      =
      (1-aw)^{-n}\bigl[(1-aw)^{b/a}-1\bigr]^{-n}
      \sum_{k=0}^n E_{n,k}(a,b;\,c_0,c_\infty)\, (1-aw)^{kb/a}.
      \label{eq:Eabc}
    \end{aligned}
  \end{sizealign}
\end{subequations}
In the $a\to0$ limit (the $b\to0$ limit is not considered here), these
become
\begin{subequations}
  \begin{sizealign}{\small}
    &
    \begin{aligned}
      &
      \bigl[ {\rm e}^{rw} (1-{\rm e}^{bw})^{-s} \bigr]^{-1}
      \,D_w^n \bigl[ {\rm e}^{rw} (1-{\rm e}^{bw})^{-s} \bigr]
      \\
      &
      \qquad\qquad
      =
      \sum_{k=0}^n b^k s^{\overline k}\,S_{n,k}(0,b;\,r)
      \bigl[{\rm e}^{bw}(1-{\rm e}^{bw})^{-1}\bigr]^k
      \label{eq:S0bc}
    \end{aligned}
    \shortintertext{\normalsize{and}}
    &
    \begin{aligned}
      &
      \bigl[{\rm e}^{-c_0w} ({\rm e}^{-bw}-1)^{-(c_0+c_\infty)/b}\bigr]
      ^{-1}
      \,D_w^n
      \bigl[{\rm e}^{-c_0w} ({\rm e}^{-bw}-1)^{-(c_0+c_\infty)/b}\bigr]
      \\
      &
      \qquad\qquad
      =({\rm e}^{-bw} - 1)^{-n}
      \sum_{k=0}^n E_{n,k}(0,b;\,c_0,c_\infty) {\rm e}^{-kbw}.
      \label{eq:E0bc}
    \end{aligned}
  \end{sizealign}
\end{subequations}
\setlength\tabcolsep{3.0pt}

By applying (\ref{eq:Sabc}) and~(\ref{eq:S0bc}), the parameter~$s$
being arbitrary, one can compute Hsu--Shiue numbers $S_{n,k}(a,b;r)$
by repeated differentiation.  For instance, the De~Morgan numbers
$\textrm{Surj}(n,k)= k!\,\stirsub{n}k = (1)^{\overline k}
S_{n,k}(0,1;0)$, which count the number of maps from an
$n$\nobreakdash-set onto a $k$\nobreakdash-set (see
Example~\ref{ex:demorgan}), satisfy
\begin{equation}
  (1-{\rm e}^w)D_w^n\left[(1-{\rm e}^w)^{-1}
    \right] = \sum_{k=0}^n \textrm{Surj}(n,k)\left[
      {\rm e}^w (1-{\rm e}^w)^{-1}
      \right]^k.
\end{equation}
The generalized Eulerian numbers $E_{n,k}(a,b;c_0,c_\infty)$ can be
computed likewise from (\ref{eq:Eabc}) and~(\ref{eq:E0bc}), with the
latter applying in the Carlitz--Scoville $a=0$ case, examples of which
were mentioned in the introduction.







\begin{thebibliography}{67}
\providecommand{\natexlab}[1]{#1}
\providecommand{\url}[1]{\texttt{#1}}
\providecommand{\urlprefix}{URL }
\expandafter\ifx\csname urlstyle\endcsname\relax
  \providecommand{\doi}[1]{doi:\discretionary{}{}{}#1}\else
  \providecommand{\doi}[1]{doi:\discretionary{}{}{}\begingroup
  \urlstyle{rm}\url{#1}\endgroup}\fi
\providecommand{\bibinfo}[2]{#2}

\bibitem[{Bagno et~al.(2019)Bagno, Biagioli, and Garber}]{Bagno2019}
\bibinfo{author}{E.~Bagno}, \bibinfo{author}{R.~Biagioli},
  \bibinfo{author}{D.~Garber}, \bibinfo{title}{Some identities involving second
  kind {S}tirling numbers of types~{$B$} and~{$D$}},
  \bibinfo{journal}{Electron. J.~Combin.}
  \bibinfo{volume}{26}~(\bibinfo{number}{3}) (\bibinfo{year}{2019})
  \bibinfo{pages}{Paper No.~3.9, 20~pp.}, \bibinfo{note}{available on-line as
  arXiv:1901.07830 [math.CO]}.

\bibitem[{Bagno et~al.(2022)Bagno, Garber, and Novick}]{Bagno2022}
\bibinfo{author}{E.~Bagno}, \bibinfo{author}{D.~Garber},
  \bibinfo{author}{M.~Novick}, \bibinfo{title}{The {W}orpitzky Identity for the
  groups of signed and even-signed permutations},
  \bibinfo{journal}{J.~Algebraic Combin.}
  \bibinfo{volume}{55}~(\bibinfo{number}{2}) (\bibinfo{year}{2022})
  \bibinfo{pages}{413--428}, \bibinfo{note}{available on-line as
  arXiv:2004.03681 [math.CO]}.

\bibitem[{Balas(2015)}]{Balas2005}
\bibinfo{author}{P.~Balas}, \bibinfo{title}{Notes on generalized {R}iordan
  arrays}, \bibinfo{note}{preprint, currently available at {\tt
  http://oeis.org/A260492/a260492.pdf}}, \bibinfo{year}{2015}.

\bibitem[{Barbero~G. et~al.(2014)Barbero~G., Salas, and
  Villase{\~n}or}]{Barbero2014}
\bibinfo{author}{J.~F. Barbero~G.}, \bibinfo{author}{J.~Salas},
  \bibinfo{author}{E.~J.~S. Villase{\~n}or}, \bibinfo{title}{Bivariate
  generating functions for a class of linear recurrences: General structure},
  \bibinfo{journal}{J.~Combin. Theory Ser.~A} \bibinfo{volume}{125}
  (\bibinfo{year}{2014}) \bibinfo{pages}{146--165}.

\bibitem[{Barbero~G. et~al.(2015)Barbero~G., Salas, and
  Villase{\~n}or}]{Barbero2015}
\bibinfo{author}{J.~F. Barbero~G.}, \bibinfo{author}{J.~Salas},
  \bibinfo{author}{E.~J.~S. Villase{\~n}or}, \bibinfo{title}{Generalized
  {S}tirling permutations and forests: Higher-order {E}ulerian and {W}ard
  numbers}, \bibinfo{journal}{Electron. J.~Combin.}
  \bibinfo{volume}{22}~(\bibinfo{number}{3}) (\bibinfo{year}{2015})
  \bibinfo{pages}{Paper No.~3.37, 20~pp.}

\bibitem[{Barry(2016)}]{Barry2016}
\bibinfo{author}{P.~Barry}, \bibinfo{title}{Riordan Arrays: A Primer},
  \bibinfo{publisher}{Logic Press}, \bibinfo{address}{Kilcock, County Kildare,
  Ireland}, \bibinfo{year}{2016}.

\bibitem[{B{\'e}nyi et~al.(2022)B{\'e}nyi, Nkonkobe, and Shattuck}]{Benyi2002}
\bibinfo{author}{B.~B{\'e}nyi}, \bibinfo{author}{S.~Nkonkobe},
  \bibinfo{author}{M.~Shattuck}, \bibinfo{title}{Unfair distributions counted
  by the generalized {S}tirling numbers}, \bibinfo{journal}{Integers}
  \bibinfo{volume}{22}~(\bibinfo{number}{article id A19})
  (\bibinfo{year}{2022}) \bibinfo{pages}{28~pp.}

\bibitem[{Boros and Moll(2004)}]{Boros2004}
\bibinfo{author}{G.~Boros}, \bibinfo{author}{V.~H. Moll},
  \bibinfo{title}{Irresistible Integrals : Symbolics, Analysis, and Experiments
  in the Evaluation of Integrals}, \bibinfo{publisher}{Cambridge Univ. Press},
  \bibinfo{address}{Cambridge, UK}, \bibinfo{year}{2004}.

\bibitem[{Boyadzhiev(2018)}]{Boyadzhiev2018}
\bibinfo{author}{K.~N. Boyadzhiev}, \bibinfo{title}{Notes on the Binomial
  Transform. Theory and Table with Appendix on {S}tirling Transform},
  \bibinfo{publisher}{World Scientific}, \bibinfo{address}{Hackensack, NJ},
  \bibinfo{year}{2018}.

\bibitem[{Brenti(1994)}]{Brenti94}
\bibinfo{author}{F.~Brenti}, \bibinfo{title}{{$q$}-{E}ulerian polynomials
  arising from {C}oxeter groups}, \bibinfo{journal}{European J.~Combin.}
  \bibinfo{volume}{15}~(\bibinfo{number}{5}) (\bibinfo{year}{1994})
  \bibinfo{pages}{417--441}.

\bibitem[{Broder(1984)}]{Broder84}
\bibinfo{author}{A.~Z. Broder}, \bibinfo{title}{The {$r$}-{S}tirling numbers},
  \bibinfo{journal}{Discrete Math.} \bibinfo{volume}{49}~(\bibinfo{number}{3})
  (\bibinfo{year}{1984}) \bibinfo{pages}{241--259}.

\bibitem[{Can and Da{\u{g}}l{\i}(2014)}]{Can2014}
\bibinfo{author}{M.~Can}, \bibinfo{author}{M.~C. Da{\u{g}}l{\i}},
  \bibinfo{title}{Extended {B}ernoulli and {S}tirling matrices and related
  combinatorial identities}, \bibinfo{journal}{Linear Algebra Appl.}
  \bibinfo{volume}{444} (\bibinfo{year}{2014}) \bibinfo{pages}{114--131},
  \bibinfo{note}{available on-line as arXiv:1306.5888 [math.NT]}.

\bibitem[{Carlitz(1978)}]{Carlitz78}
\bibinfo{author}{L.~Carlitz}, \bibinfo{title}{Generalized {S}tirling and
  related numbers}, \bibinfo{journal}{Riv.~Mat. Univ. Parma~(4)}
  \bibinfo{volume}{4} (\bibinfo{year}{1978}) \bibinfo{pages}{79--99}.

\bibitem[{Carlitz(1979)}]{Carlitz79}
\bibinfo{author}{L.~Carlitz}, \bibinfo{title}{Degenerate {S}tirling,
  {B}ernoulli and {E}ulerian numbers}, \bibinfo{journal}{Utilitas Math.}
  \bibinfo{volume}{15} (\bibinfo{year}{1979}) \bibinfo{pages}{51--88}.

\bibitem[{Carlitz and Scoville(1974)}]{Carlitz74}
\bibinfo{author}{L.~Carlitz}, \bibinfo{author}{R.~Scoville},
  \bibinfo{title}{Generalized {E}ulerian numbers: Combinatorial applications},
  \bibinfo{journal}{J.~Reine Angew. Math.} \bibinfo{volume}{265}
  (\bibinfo{year}{1974}) \bibinfo{pages}{110--137}.

\bibitem[{Charalambides(1982)}]{Charalambides82}
\bibinfo{author}{Ch.~A. Charalambides}, \bibinfo{title}{On the enumeration of
  certain compositions and related sequences of numbers},
  \bibinfo{journal}{Fibonacci Quart.}
  \bibinfo{volume}{20}~(\bibinfo{number}{2}) (\bibinfo{year}{1982})
  \bibinfo{pages}{132--146}.

\bibitem[{Charalambides(1991)}]{Charalambides91}
\bibinfo{author}{Ch.~A. Charalambides}, \bibinfo{title}{On a generalized
  {E}ulerian distribution}, \bibinfo{journal}{Ann. Inst. Statist. Math.}
  \bibinfo{volume}{43}~(\bibinfo{number}{1}) (\bibinfo{year}{1991})
  \bibinfo{pages}{197--206}.

\bibitem[{Charalambides and Koutras(1983)}]{Charalambides83}
\bibinfo{author}{Ch.~A. Charalambides}, \bibinfo{author}{M.~Koutras},
  \bibinfo{title}{On the differences of the generalized factorials at an
  arbitrary point and their combinatorial applications},
  \bibinfo{journal}{Discrete Math.} \bibinfo{volume}{47} (\bibinfo{year}{1983})
  \bibinfo{pages}{183--201}.

\bibitem[{Chen(1993)}]{Chen93}
\bibinfo{author}{W.~Y.~C. Chen}, \bibinfo{title}{Context-free grammars,
  differential operators and formal power series}, \bibinfo{journal}{Theoret.
  Comput. Sci.} \bibinfo{volume}{117}~(\bibinfo{number}{1--2})
  (\bibinfo{year}{1993}) \bibinfo{pages}{113--129}.

\bibitem[{Chen et~al.(2011)Chen, Wang, and Zhao}]{Chen2011}
\bibinfo{author}{W.~Y.~C. Chen}, \bibinfo{author}{A.~Y.~Z. Wang},
  \bibinfo{author}{A.~F.~Y. Zhao}, \bibinfo{title}{Identities derived from
  noncrossing partitions of type~{$B$}}, \bibinfo{journal}{Electron.
  J.~Combin.} \bibinfo{volume}{18}~(\bibinfo{number}{article id 129})
  (\bibinfo{year}{2011}) \bibinfo{pages}{17~pp.}, \bibinfo{note}{available
  on-line as arXiv:0908.2291 [math.CO]}.

\bibitem[{Cheon et~al.(2013)Cheon, Jung, and Shapiro}]{Cheon2013}
\bibinfo{author}{G.-S. Cheon}, \bibinfo{author}{J.-H. Jung},
  \bibinfo{author}{L.~W. Shapiro}, \bibinfo{title}{Generalized {B}essel numbers
  and some combinatorial settings}, \bibinfo{journal}{Discrete Math.}
  \bibinfo{volume}{313}~(\bibinfo{number}{20}) (\bibinfo{year}{2013})
  \bibinfo{pages}{2127--2138}.

\bibitem[{Comtet(1974)}]{Comtet74}
\bibinfo{author}{L.~Comtet}, \bibinfo{title}{Advanced Combinatorics. The Art of
  Finite and Infinite Expansions}, \bibinfo{publisher}{Reidel},
  \bibinfo{address}{Boston/Dordrecht}, \bibinfo{year}{1974}.

\bibitem[{Corcino(2001)}]{Corcino2001}
\bibinfo{author}{R.~B. Corcino}, \bibinfo{title}{Some theorems on generalized
  {S}tirling numbers}, \bibinfo{journal}{Ars Combin.} \bibinfo{volume}{60}
  (\bibinfo{year}{2001}) \bibinfo{pages}{273--286}.

\bibitem[{Corcino et~al.(2001)Corcino, Hsu, and Tan}]{Corcino2001a}
\bibinfo{author}{R.~B. Corcino}, \bibinfo{author}{L.~C. Hsu},
  \bibinfo{author}{E.~L. Tan}, \bibinfo{title}{Combinatorial and statistical
  applications of generalized {Stirling} numbers}, \bibinfo{journal}{J.~Math.
  Res. Exposition} \bibinfo{volume}{21}~(\bibinfo{number}{3})
  (\bibinfo{year}{2001}) \bibinfo{pages}{337--343}.

\bibitem[{Dillon and Roselle(1968)}]{Dillon68}
\bibinfo{author}{J.~F. Dillon}, \bibinfo{author}{D.~P. Roselle},
  \bibinfo{title}{Eulerian numbers of higher order}, \bibinfo{journal}{Duke
  Math.~J.} \bibinfo{volume}{35} (\bibinfo{year}{1968})
  \bibinfo{pages}{247--256}.

\bibitem[{Dumont(1996)}]{Dumont96}
\bibinfo{author}{D.~Dumont}, \bibinfo{title}{Grammaires de {W}illiam {C}hen et
  D{\'e}rivations dans les arbres et arborescences}, \bibinfo{journal}{S{\'e}m.
  Lothar. Combin.} \bibinfo{volume}{37}~(\bibinfo{number}{article id B37a})
  (\bibinfo{year}{1996}) \bibinfo{pages}{21~pp.}

\bibitem[{Encinas and Jim{\'e}nez(2019)}]{Encinas2019}
\bibinfo{author}{A.~M. Encinas}, \bibinfo{author}{M.~J. Jim{\'e}nez},
  \bibinfo{title}{Triangular sequences, combinatorial recurrences and linear
  difference equations}, \bibinfo{journal}{Linear Algebra Appl.}
  \bibinfo{volume}{576} (\bibinfo{year}{2019}) \bibinfo{pages}{301--323}.

\bibitem[{Foata(2010)}]{Foata2010}
\bibinfo{author}{D.~Foata}, \bibinfo{title}{Eulerian polynomials: From
  {E}uler's time to the present}, in: \bibinfo{editor}{K.~Alladi},
  \bibinfo{editor}{J.~R. Klauder}, \bibinfo{editor}{C.~R. Rao} (Eds.),
  \bibinfo{booktitle}{The Legacy of {A}lladi {R}amakrishnan in the Mathematical
  Sciences}, \bibinfo{publisher}{Springer}, \bibinfo{address}{New York},
  \bibinfo{pages}{253--273}, \bibinfo{year}{2010}.

\bibitem[{Foata and Sch{\"u}tzenberger(1970)}]{Foata70}
\bibinfo{author}{D.~Foata}, \bibinfo{author}{M.-P. Sch{\"u}tzenberger},
  \bibinfo{title}{Th{\'e}orie g{\'e}ometrique des p{\^o}lynomes Eul{\'e}riens},
  no. \bibinfo{number}{138} in \bibinfo{series}{Lecture Notes in Mathematics},
  \bibinfo{publisher}{Springer-Verlag}, \bibinfo{address}{New York/Berlin},
  \bibinfo{year}{1970}.

\bibitem[{Fomin and Reading(2007)}]{Fomin2007}
\bibinfo{author}{S.~Fomin}, \bibinfo{author}{N.~Reading}, \bibinfo{title}{Root
  systems and generalized associahedra}, in: \bibinfo{editor}{E.~Miller},
  \bibinfo{editor}{V.~Reiner}, \bibinfo{editor}{B.~Sturmfels} (Eds.),
  \bibinfo{booktitle}{Geometric Combinatorics}, \bibinfo{publisher}{American
  Mathematical Society (AMS)}, \bibinfo{address}{Providence, RI},
  \bibinfo{pages}{63--131}, \bibinfo{note}{available on-line as
  arXiv:math/0505518 [math.CO]}, \bibinfo{year}{2007}.

\bibitem[{Galuzzi(1998)}]{Galuzzi98}
\bibinfo{author}{M.~Galuzzi}, \bibinfo{title}{A remark about the binomial
  transform}, \bibinfo{journal}{Fibonacci Quart.}
  \bibinfo{volume}{36}~(\bibinfo{number}{3}) (\bibinfo{year}{1998})
  \bibinfo{pages}{287--288}.

\bibitem[{Gessel(1977)}]{Gessel77}
\bibinfo{author}{I.~M. Gessel}, \bibinfo{title}{Generating Functions and
  Enumeration of Sequences}, \bibinfo{type}{{Ph.D.} dissertation},
  \bibinfo{school}{Massachusetts Institute of Technology},
  \bibinfo{year}{1977}.

\bibitem[{Graham et~al.(1994)Graham, Knuth, and Patashnik}]{Graham94}
\bibinfo{author}{R.~L. Graham}, \bibinfo{author}{D.~E. Knuth},
  \bibinfo{author}{O.~Patashnik}, \bibinfo{title}{Concrete Mathematics},
  \bibinfo{publisher}{Addison--Wesley}, \bibinfo{address}{Boston},
  \bibinfo{edition}{2nd} edn., \bibinfo{year}{1994}.

\bibitem[{Han and Seo(2008)}]{Han2008}
\bibinfo{author}{H.~Han}, \bibinfo{author}{S.~Seo},
  \bibinfo{title}{Combinatorial proofs of inverse relations and log-concavity
  for {B}essel numbers}, \bibinfo{journal}{European J.~Combin.}
  \bibinfo{volume}{29}~(\bibinfo{number}{7}) (\bibinfo{year}{2008})
  \bibinfo{pages}{1544--1554}, \bibinfo{note}{available on-line as
  arXiv:math/0406378 [math.CO]}.

\bibitem[{Hao et~al.(2015)Hao, Wang, and Yang}]{Hao2015}
\bibinfo{author}{R.~X.~J. Hao}, \bibinfo{author}{L.~X.~W. Wang},
  \bibinfo{author}{H.~R.~L. Yang}, \bibinfo{title}{Context-free grammars for
  triangular arrays}, \bibinfo{journal}{Acta Math. Sin. (Engl. Ser.)}
  \bibinfo{volume}{31}~(\bibinfo{number}{3}) (\bibinfo{year}{2015})
  \bibinfo{pages}{445--455}.

\bibitem[{Harris and Park(1994)}]{Harris94}
\bibinfo{author}{B.~Harris}, \bibinfo{author}{C.~J. Park}, \bibinfo{title}{A
  generalization of the {E}ulerian numbers with a probabilistic application},
  \bibinfo{journal}{Statist. Probab. Lett.}
  \bibinfo{volume}{20}~(\bibinfo{number}{1}) (\bibinfo{year}{1994})
  \bibinfo{pages}{37--47}.

\bibitem[{He(2014)}]{He2013}
\bibinfo{author}{T.-X. He}, \bibinfo{title}{Expression and computation of
  generalized {S}tirling numbers}, \bibinfo{journal}{J.~Combin. Math. Combin.
  Comput.} \bibinfo{volume}{86} (\bibinfo{year}{2014})
  \bibinfo{pages}{239--268}.

\bibitem[{He and Shiue(2020)}]{He2020}
\bibinfo{author}{T.-X. He}, \bibinfo{author}{P.~J.-S. Shiue}, \bibinfo{title}{A
  note on {E}ulerian numbers and {T}oeplitz matrices}, \bibinfo{journal}{Spec.
  Matrices} \bibinfo{volume}{8} (\bibinfo{year}{2020})
  \bibinfo{pages}{123--130}.

\bibitem[{Herscovici(2020)}]{Herscovici2020}
\bibinfo{author}{O.~Herscovici}, \bibinfo{title}{Generalized permutations
  related to the degenerate {E}ulerian numbers}, \bibinfo{note}{preprint,
  available on-line as arXiv:2007.13205 [math.CO]}, \bibinfo{year}{2020}.

\bibitem[{Hsu and Shiue(1998)}]{Hsu98}
\bibinfo{author}{L.~C. Hsu}, \bibinfo{author}{P.~J.-S. Shiue},
  \bibinfo{title}{A unified approach to generalized {S}tirling numbers},
  \bibinfo{journal}{Adv. in Appl. Math.}
  \bibinfo{volume}{20}~(\bibinfo{number}{3}) (\bibinfo{year}{1998})
  \bibinfo{pages}{366--384}.

\bibitem[{Hsu and Shiue(1999)}]{Hsu99}
\bibinfo{author}{L.~C. Hsu}, \bibinfo{author}{P.~J.-S. Shiue},
  \bibinfo{title}{On certain summation problems and generalizations of
  {E}ulerian polynomials and numbers}, \bibinfo{journal}{Discrete Math.}
  \bibinfo{volume}{204}~(\bibinfo{number}{1--3}) (\bibinfo{year}{1999})
  \bibinfo{pages}{237--247}.

\bibitem[{Janardan(1988)}]{Janardan88}
\bibinfo{author}{K.~G. Janardan}, \bibinfo{title}{Relationship between
  {M}orisita's model for estimating the environmental density and the
  generalized {E}ulerian numbers}, \bibinfo{journal}{Ann. Inst. Statist. Math.}
  \bibinfo{volume}{40}~(\bibinfo{number}{3}) (\bibinfo{year}{1988})
  \bibinfo{pages}{439--450}.

\bibitem[{Ma(2012)}]{Ma2012}
\bibinfo{author}{S.-M. Ma}, \bibinfo{title}{Derivative polynomials and
  enumeration of permutations by number of interior and left peaks},
  \bibinfo{journal}{Discrete Math.} \bibinfo{volume}{312}~(\bibinfo{number}{2})
  (\bibinfo{year}{2012}) \bibinfo{pages}{405--412}, \bibinfo{note}{available
  on-line as arXiv:1504.02372 [math.CO]}.

\bibitem[{Ma et~al.(2019)Ma, Ma, and Yeh}]{Ma2019}
\bibinfo{author}{S.-M. Ma}, \bibinfo{author}{J.~Ma}, \bibinfo{author}{Y.-N.
  Yeh}, \bibinfo{title}{{$\gamma$}-positivity and partial {$\gamma$}-positivity
  of descent-type polynomials}, \bibinfo{journal}{J.~Combin. Theory Ser.~A}
  \bibinfo{volume}{167} (\bibinfo{year}{2019}) \bibinfo{pages}{257--293},
  \bibinfo{note}{available on-line as arXiv:1802.02861 [math.CO]}.

\bibitem[{Ma and Mansour(2015)}]{Ma2015}
\bibinfo{author}{S.-M. Ma}, \bibinfo{author}{T.~Mansour}, \bibinfo{title}{The
  {$1/k$}-{E}ulerian polynomials and {$k$}-{S}tirling permutations},
  \bibinfo{journal}{Discrete Math.} \bibinfo{volume}{338}~(\bibinfo{number}{8})
  (\bibinfo{year}{2015}) \bibinfo{pages}{1468--1472}.

\bibitem[{Maltenfort(2020)}]{Maltenfort2020}
\bibinfo{author}{M.~Maltenfort}, \bibinfo{title}{New definitions of the
  generalized {S}tirling numbers}, \bibinfo{journal}{Aequationes Math.}
  \bibinfo{volume}{94}~(\bibinfo{number}{1}) (\bibinfo{year}{2020})
  \bibinfo{pages}{169--200}.

\bibitem[{Mansour and Schork(2016)}]{Mansour2016}
\bibinfo{author}{T.~Mansour}, \bibinfo{author}{M.~Schork},
  \bibinfo{title}{Commutation Relations, Normal Ordering, and {S}tirling
  Numbers}, \bibinfo{publisher}{CRC Press}, \bibinfo{address}{Boca Raton, FL},
  \bibinfo{year}{2016}.

\bibitem[{Mansour and Shattuck(2013)}]{Mansour2013}
\bibinfo{author}{T.~Mansour}, \bibinfo{author}{M.~Shattuck}, \bibinfo{title}{A
  combinatorial approach to a general two-term recurrence},
  \bibinfo{journal}{Discrete Appl. Math.}
  \bibinfo{volume}{161}~(\bibinfo{number}{13--14}) (\bibinfo{year}{2013})
  \bibinfo{pages}{2084--2094}.

\bibitem[{Miller(1966)}]{Miller66}
\bibinfo{author}{K.~S. Miller}, \bibinfo{title}{An Introduction to the Calculus
  of Finite Differences and Difference Equations}, \bibinfo{publisher}{Dover},
  \bibinfo{address}{New York}, \bibinfo{year}{1966}.

\bibitem[{Neuwirth(2001)}]{Neuwirth2001}
\bibinfo{author}{E.~Neuwirth}, \bibinfo{title}{Recursively defined
  combinatorial functions: Extending {G}alton's board},
  \bibinfo{journal}{Discrete Math.}
  \bibinfo{volume}{239}~(\bibinfo{number}{1--3}) (\bibinfo{year}{2001})
  \bibinfo{pages}{33--51}.

\bibitem[{Nyul and R{\'a}cz(2015)}]{Nyul2015}
\bibinfo{author}{G.~Nyul}, \bibinfo{author}{G.~R{\'a}cz}, \bibinfo{title}{The
  {$r$}-{L}ah numbers}, \bibinfo{journal}{Discrete Math.}
  \bibinfo{volume}{338}~(\bibinfo{number}{10}) (\bibinfo{year}{2015})
  \bibinfo{pages}{1660--1666}.

\bibitem[{{OEIS Foundation, Inc.}(2022)}]{OEIS2022}
\bibinfo{author}{{OEIS Foundation, Inc.}}, \bibinfo{title}{The {O}n-{L}ine
  {E}ncyclopedia of {I}nteger {S}equences}, \bibinfo{howpublished}{Published
  electronically at {\tt https://oeis.org}}, \bibinfo{year}{2022}.

\bibitem[{Olver et~al.(2010)Olver, Lozier, Boisvert, and Clark}]{Olver2010}
\bibinfo{editor}{F.~W.~J. Olver}, \bibinfo{editor}{D.~W. Lozier},
  \bibinfo{editor}{R.~F. Boisvert}, \bibinfo{editor}{C.~W. Clark} (Eds.),
  \bibinfo{title}{{NIST} Handbook of Mathematical Functions},
  \bibinfo{publisher}{U.S. Department of Commerce, National Institute of
  Standards and Technology}, \bibinfo{address}{Washington, DC},
  \bibinfo{year}{2010}.

\bibitem[{Petersen(2015)}]{Petersen2015}
\bibinfo{author}{T.~K. Petersen}, \bibinfo{title}{Eulerian Numbers},
  \bibinfo{publisher}{Birkh{\"a}user/Springer}, \bibinfo{address}{New York},
  \bibinfo{note}{with a foreword by {R}ichard {S}tanley}, \bibinfo{year}{2015}.

\bibitem[{Riordan(1958)}]{Riordan58}
\bibinfo{author}{J.~Riordan}, \bibinfo{title}{An Introduction to Combinatorial
  Analysis}, \bibinfo{publisher}{Wiley}, \bibinfo{address}{New York},
  \bibinfo{year}{1958}.

\bibitem[{Rosenberg(2019)}]{Rosenberg2019}
\bibinfo{author}{N.~A. Rosenberg}, \bibinfo{title}{Enumeration of lonely pairs
  of gene trees and species trees by means of antipodal cherries},
  \bibinfo{journal}{Adv. in Appl. Math.} \bibinfo{volume}{102}
  (\bibinfo{year}{2019}) \bibinfo{pages}{1--17}.

\bibitem[{Salas and Sokal(2021)}]{Salas2021}
\bibinfo{author}{J.~Salas}, \bibinfo{author}{A.~D. Sokal}, \bibinfo{title}{The
  {G}raham--{K}nuth--{P}atashnik recurrence: Symmetries and continued
  fractions}, \bibinfo{journal}{Electron. J.~Combin.}
  \bibinfo{volume}{28(2)}~(\bibinfo{number}{article id P2.18})
  (\bibinfo{year}{2021}) \bibinfo{pages}{72~pp.}, \bibinfo{note}{available
  on-line as arXiv:2008.03070 [math.CO]}.

\bibitem[{Singh~Chandel(1977)}]{SinghChandel77}
\bibinfo{author}{R.~C. Singh~Chandel}, \bibinfo{title}{Generalized {S}tirling
  numbers and polynomials}, \bibinfo{journal}{Publ. Inst. Math. (Beograd)
  (N.S.)} \bibinfo{volume}{22}~(\bibinfo{number}{36}) (\bibinfo{year}{1977})
  \bibinfo{pages}{43--48}.

\bibitem[{Spivey(2011)}]{Spivey2011}
\bibinfo{author}{M.~Z. Spivey}, \bibinfo{title}{On solutions to a general
  combinatorial recurrence}, \bibinfo{journal}{J.~Integer Seq.}
  \bibinfo{volume}{14(9)}~(\bibinfo{number}{article id 11.9.7})
  (\bibinfo{year}{2011}) \bibinfo{pages}{19~pp.}, \bibinfo{note}{available
  on-line as arXiv:1307.2010 [math.CO]}.

\bibitem[{Srivastava and Singhal(1973)}]{Srivastava73}
\bibinfo{author}{H.~M. Srivastava}, \bibinfo{author}{J.~P. Singhal},
  \bibinfo{title}{New generating functions for {J}acobi and related
  polynomials}, \bibinfo{journal}{J.~Math. Anal. Appl.}
  \bibinfo{volume}{41}~(\bibinfo{number}{3}) (\bibinfo{year}{1973})
  \bibinfo{pages}{748--752}.

\bibitem[{Stanton and Sprott(1962)}]{Stanton62}
\bibinfo{author}{R.~G. Stanton}, \bibinfo{author}{D.~A. Sprott},
  \bibinfo{title}{Some finite inversion formulae}, \bibinfo{journal}{Math.
  Gaz.} \bibinfo{volume}{46} (\bibinfo{year}{1962}) \bibinfo{pages}{197--202}.

\bibitem[{Tak{\'a}cs(1979)}]{Takacs79}
\bibinfo{author}{L.~Tak{\'a}cs}, \bibinfo{title}{A generalization of the
  {E}ulerian numbers}, \bibinfo{journal}{Publ. Math. Debrecen}
  \bibinfo{volume}{26}~(\bibinfo{number}{3--4}) (\bibinfo{year}{1979})
  \bibinfo{pages}{173--181}.

\bibitem[{Th{\'e}or{\^e}t(1994)}]{Theoret94}
\bibinfo{author}{P.~Th{\'e}or{\^e}t}, \bibinfo{title}{Hyperbinomiales: Double
  Suites Satisfaisant \`a des \'Equations aux Diff{\'e}rences Partielles de
  Dimension et D'Ordre Deux de la Forme
  {$H(n,k)=p(n,k)H(n-1,k)+q(n,k)H(n-1,k-1)$}}, \bibinfo{type}{{Ph.D.}
  dissertation}, \bibinfo{school}{University of Qu{\'e}bec at Montr{\'e}al},
  \bibinfo{year}{1994}.

\bibitem[{Th{\'e}or{\^e}t(1995)}]{Theoret95}
\bibinfo{author}{P.~Th{\'e}or{\^e}t}, \bibinfo{title}{Fonctions
  g{\'e}n{\'e}ratrices pour une classe d'{\'e}quations aux diff{\'e}rences
  partielles}, \bibinfo{journal}{Ann. Sci. Math. Qu{\'e}bec}
  \bibinfo{volume}{19}~(\bibinfo{number}{1}) (\bibinfo{year}{1995})
  \bibinfo{pages}{91--105}.

\bibitem[{Wilf(2004)}]{Wilf2004}
\bibinfo{author}{H.~S. Wilf}, \bibinfo{title}{The method of characteristics,
  and {``Problem 89''} of {G}raham, {K}nuth and {P}atashnik},
  \doi{\bibinfo{doi}{10.48550/arXiv.math/0406620}}, \bibinfo{note}{preprint,
  available on-line as arXiv:math/0406620 [math.CO]}, \bibinfo{year}{2004}.

\bibitem[{Xiong et~al.(2013)Xiong, Tsao, and Hall}]{Xiong2013}
\bibinfo{author}{T.~Xiong}, \bibinfo{author}{H.-P. Tsao},
  \bibinfo{author}{J.~I. Hall}, \bibinfo{title}{General {E}ulerian numbers and
  {E}ulerian polynomials}, \bibinfo{journal}{J.~Math.}
  \bibinfo{volume}{2013}~(\bibinfo{number}{article id 629132})
  (\bibinfo{year}{2013}) \bibinfo{pages}{9~pp.}, \bibinfo{note}{available
  on-line as arXiv:1207.0430 [math.CO]}.

\bibitem[{Zhuang(2016)}]{Zhuang2016}
\bibinfo{author}{Y.~Zhuang}, \bibinfo{title}{Counting permutations by runs},
  \bibinfo{journal}{J.~Combin. Theory Ser.~A} \bibinfo{volume}{142}
  (\bibinfo{year}{2016}) \bibinfo{pages}{147--176}, \bibinfo{note}{available
  on-line as arXiv:1505.02308 [math.CO]}.

\end{thebibliography}



\end{document}